\documentclass[12pt,reqno]{amsart}
\usepackage[letterpaper,margin=0.9in]{geometry}
\usepackage[T1]{fontenc}
\usepackage{lmodern,amsmath,amssymb,mathtools,mathrsfs,microtype}
\usepackage{xcolor,booktabs,tabularx,array,enumitem}
\usepackage[hidelinks]{hyperref}
\hypersetup{pdftitle={Finite conformal models and the Kac--van de Leur conjecture},pdfauthor={Haibo Chen, Mingqian Liu, and Yong Yang},pdfsubject={Exhaustion of simple superconformal Lie algebras; finite conformal models and twisted coefficients}}
\mathtoolsset{showonlyrefs}
\numberwithin{equation}{section}
\newtheorem{theorem}{Theorem}[section]
\newtheorem{proposition}[theorem]{Proposition}
\newtheorem{lemma}[theorem]{Lemma}
\newtheorem{corollary}[theorem]{Corollary}
\theoremstyle{definition}
\newtheorem{definition}[theorem]{Definition}
\theoremstyle{remark}
\newtheorem{remark}[theorem]{Remark}
\newcounter{furtherquestion}
\newenvironment{questionnote}[2]{\par\medskip\refstepcounter{furtherquestion}\label{#1}\noindent\textbf{Further question \thefurtherquestion\ (#2).}\ }{\par\medskip}
\newcommand{\C}{\mathbb C}
\newcommand{\Z}{\mathbb Z}

\newcommand{\Vir}{\mathfrak V}
\newcommand{\g}{\mathfrak g}

\newcommand{\ad}{\operatorname{ad}}
\newcommand{\Der}{\operatorname{Der}}
\newcommand{\Cder}{\operatorname{Cder}}
\newcommand{\Inder}{\operatorname{Inder}}
\newcommand{\Alg}{\operatorname{Alg}}

\newcommand{\Rad}{\operatorname{Rad}}
\newcommand{\Soc}{\operatorname{Soc}}
\newcommand{\Cur}{\operatorname{Cur}}
\newcommand{\im}{\operatorname{Im}}
\newcommand{\rank}{\operatorname{rank}}
\newcommand{\Span}{\operatorname{span}}

\newcommand{\diver}{\operatorname{div}}
\newcommand{\End}{\operatorname{End}}
\newcommand{\Aut}{\operatorname{Aut}}
\newcommand{\slie}{\mathfrak{sl}}

\newcommand{\dd}{\partial}
\newcommand{\wt}{\operatorname{wt}}
\newcommand{\Acal}{\mathcal A}
\newcommand{\Gcal}{\mathfrak G}

\newcommand{\Pcal}{\mathscr P}

\newcommand{\sm}{\mathrm s}
\newcommand{\nil}{\mathrm n}
\setlist[enumerate]{leftmargin=*,itemsep=3pt,topsep=5pt}
\allowdisplaybreaks[2]
\newcommand{\Div}{\operatorname{Div}}
\newcommand{\id}{\operatorname{id}}
\newcommand{\Prim}{\operatorname{Prim}}
\newcommand{\GL}{\operatorname{GL}}
\newcommand{\SO}{\operatorname{SO}}
\newcommand{\OO}{\operatorname{O}}
\newcommand{\supp}{\operatorname{supp}}
\newcommand{\s}{\mathfrak s}

\newcommand{\nn}[1]{\mathopen{\langle\!\langle}#1\mathclose{\rangle\!\rangle}}

\title[Finite conformal models and twisted coefficients]{Finite conformal models and the Kac--van de Leur conjecture}
\author[H. Chen]{Haibo Chen}
\address{School of Science, Jimei University, Xiamen, Fujian 361021, China}
\email{hypo1025@jmu.edu.cn}
\author[M. Liu]{Mingqian Liu}
\address{School of Mathematics and Statistics, Hubei Minzu University, Enshi 445000, China}
\email{mingqiangliu@163.com}
\author[Y. Yang]{Yong Yang}
\address{Department of Mathematics, Texas State University, San Marcos, Texas 78666, USA}
\email{yang@txstate.edu}
\date{November 26, 2026}
\subjclass[2020]{Primary 17B65, 17B69; Secondary 17B68}
\keywords{Superconformal Lie superalgebra, Lie conformal superalgebra, Virasoro algebra, polynomial commutator growth, twisted coefficient algebra}
\begin{document}
\begin{abstract}
We prove that every simple Lie superalgebra containing the centerless
Virasoro algebra, with diagonal $L_0$-action, uniformly bounded weight
multiplicities, and support in finitely many cosets of $\Z$, is isomorphic,
as a parity-graded Lie superalgebra, to a member of the Witt, special,
contact, or Cheng--Kac families. Thus the amended Kac--van de Leur
list is exhaustive. Every such algebra admits a finite ordinary
conformal model and is a quotient by the full center of a twisted
coefficient algebra of a finite simple conformal superalgebra.
For Lie superalgebras containing the centerless Virasoro algebra with
diagonal integral grading and finite-dimensional weight spaces, uniformly
bounded weight multiplicities, linear commutator growth, and polynomial
commutator growth are equivalent. We also give a finitely generated simple example with finite-dimensional
integral weight spaces of unbounded dimension, showing that some growth
restriction is essential.
\end{abstract}
\maketitle
\enlargethispage{2\baselineskip}
{\footnotesize\tableofcontents}
\clearpage
\section{Introduction}\label{sec:intro}
The Kac--van de Leur conjecture asks for an exhaustive list of
$\Z$-graded superconformal algebras
\cite[Conjecture~1]{KacVanDeLeur}. In that formulation, the grading is
the eigenspace decomposition for $\ad L_0$, the homogeneous spaces are
finite-dimensional, and the algebra satisfies the linear
commutator-growth condition SC3. The proposed Witt, special, and contact
list was enlarged by the Cheng--Kac algebra $CK(6)$ \cite{ChengKac}.

We prove exhaustion for simple Lie superalgebras in the setting of
\cite[Definition~2.1]{Billig}:
the $L_0$-action is diagonal, the weight multiplicities are uniformly
bounded, and the support is contained in finitely many cosets of $\Z$
in $\C$. The classification allows nonintegral weights and
diagonalizable twisting automorphisms of arbitrary order. A separate
growth theorem shows that, in the original integral setting, SC3
implies the required uniform bound. Thus uniform boundedness is part
of Billig's definition and follows from the original growth hypothesis
in the integral setting.

All vector spaces and algebras are over $\C$. Lie algebras may be
$\Z/2\Z$-graded. All ideals and homomorphisms respect parity, and $p(a)$
denotes the parity of a homogeneous element $a$. An even derivation
preserves parity. Let $\Vir$ be the centerless Virasoro algebra, with
basis $\{L_n:n\in\Z\}$ and bracket
\[
 [L_m,L_n]=(m-n)L_{m+n}.
\]
We index the $L_0$-eigenspaces by
\begin{equation}\label{eq:weight}
 L[k]=\{x\in L:[L_0,x]=-kx\}.
\end{equation}
We call $k$ the physical weight of an element of $L[k]$ and write
$\supp L=\{k\in\C:L[k]\ne0\}$ for the support. Thus
$[L[k],L[j]]\subseteq L[k+j]$.

\begin{definition}\label{def:superconformal}
A \emph{superconformal Lie algebra} is a Lie superalgebra $L$ containing
$\Vir$, decomposing into the eigenspaces~\eqref{eq:weight}, whose support
is contained in finitely many cosets of $\Z$ in $\C$, and for which
$\dim L[k]$ is bounded independently of $k$.
\end{definition}
This agrees with \cite[Definition~2.1]{Billig}, after changing the sign
of the weight index.

A Lie conformal superalgebra is a $\C[\dd]$-module $C$ with a polynomial
$\lambda$-bracket
\[
 [a_\lambda b]=\sum_{j\ge0}\frac{\lambda^j}{j!}a_{(j)}b
\]
satisfying sesquilinearity, super-skew-symmetry, and the conformal Jacobi
identity. In particular,
\begin{equation}\label{eq:sesqui}
 [(\dd a)_\lambda b]=-\lambda[a_\lambda b],\qquad
 [a_\lambda\dd b]=(\dd+\lambda)[a_\lambda b].
\end{equation}
It is \emph{finite} if it is finitely generated over $\C[\dd]$.
For a Lie superalgebra $\g$, its current conformal algebra is
$\Cur\g=\C[\dd]\otimes\g$, with $[a_\lambda b]=[a,b]$ on $\g$.

For an even diagonalizable conformal automorphism $\sigma$, write
$C^\alpha=\{a:\sigma a=e^{-2\pi i\alpha}a\}$, where
$\alpha\in\C/\Z$. The \emph{twisted coefficient algebra}
$\Alg(C,\sigma)$ is spanned by symbols $a[m]$, linear in
$a\in C^\alpha$, with $m\in\alpha+\Z$, subject to
\[
 (\dd a)[m]=-m a[m-1],\qquad
 [a[m],b[n]]=\sum_{j\ge0}\binom mj(a_{(j)}b)[m+n-j].
\]
The sum is finite because the $\lambda$-bracket is polynomial.
We write $Z(\g)$ for the full center of a Lie superalgebra $\g$.

The noncurrent finite simple conformal superalgebras are
\begin{equation}\label{eq:list}
 W_n,\quad S_{n,a},\quad\widetilde S_{2n},\quad
 K_n\ (n\ne4),\quad K'_4,\quad CK_6,
\end{equation}
with the parameter ranges in \cite[Theorem~1.1]{FK}. The remaining
finite simple conformal algebras are current algebras of
finite-dimensional simple Lie superalgebras.

To describe the resulting Lie superalgebras, let
\[
 R_N=\C[t,t^{-1}]\otimes\Lambda(\xi_1,\ldots,\xi_N),
 \qquad W(1,N)=\Der R_N.
\]
For a parity-homogeneous vector field
$D=f\partial_t+\sum_i g_i\partial_{\xi_i}$, define
\[
 \diver D=\partial_t f+\sum_i(-1)^{p(g_i)}\partial_{\xi_i}g_i.
\]
Let $S'_B(1,N)$ be the derived algebra of
\[
 \mathfrak S_B=
 \{D\in W(1,N):\diver D+B t^{-1}D(t)=0\}.
\]
This is the formal density condition $\diver(t^BD)=0$, also for
$B\notin\Z$. The notation $CK(6)$ denotes the standard untwisted
coefficient algebra of $CK_6$.

For $N\ge1$, let $\rho$ be an orthogonal reflection, set $M_N=K_N$
if $N\ne4$ and $M_4=K'_4$, and write
\[
 \mathfrak K_N^\epsilon
   =\Alg(M_N,\rho_\epsilon)/Z(\Alg(M_N,\rho_\epsilon)),
 \qquad \rho_0=\id,\quad\rho_1=\rho,
 \quad\epsilon\in\{0,1\}.
\]
These contact representatives are described in Section~\ref{sec:contact}.
In the notation of \cite[Theorem~4.1(a)]{KacVanDeLeur}, they represent
$K'(N;\epsilon)$ according to the parity of $\sum_i\epsilon_i$: the
identity form represents even parity and a reflection represents odd
parity. In rank four, the identity form uses the derived contact algebra.

\begin{theorem}
\label{thm:synthesis}
Let $L$ be a simple superconformal Lie algebra. Its underlying
parity-graded Lie superalgebra is isomorphic to at least one member
of the following list:
\[
\begin{gathered}
 W(1,N)\quad(N\ge0),\qquad
 S'_B(1,N)\quad(N\ge2,\ B\in\C),\\
 \mathfrak K_N^\epsilon\quad(N\ge1,\ \epsilon\in\{0,1\}),
 \qquad CK(6).
\end{gathered}
\]
Different entries may be isomorphic. The isomorphisms need not preserve
a prescribed Virasoro embedding or physical grading, and the twisting
automorphisms need not have finite order.
\end{theorem}

Theorem~\ref{thm:synthesis} applies to the full class of simple algebras
in Definition~\ref{def:superconformal}. Its proof uses the finite-coset
support and uniform bound directly and is independent of the
growth theorem below.

For a finite subset $F\subset L$, let $V_j(F)$ be the span of all iterated
brackets in elements of $F$ involving at most $j$ factors. We say that
$L$ has \emph{polynomial commutator growth} if for every finite $F$ there
are constants $C_F>0$ and $d_F\geq0$ such that
\[
 \dim V_j(F)\leq C_F j^{d_F}\qquad(j\geq1).
\]
The exponent is allowed to depend on $F$. The original condition
\cite[condition SC3]{KacVanDeLeur} is the special case $d_F=1$.

\begin{theorem}\label{thm:growth-rigidity}
Let $L$ be a Lie superalgebra containing the centerless Virasoro algebra,
and suppose that
\[
 L=\bigoplus_{k\in\Z}L[k],\qquad \dim L[k]<\infty\quad(k\in\Z).
\]
The following conditions are equivalent:
\begin{enumerate}[label=\textup{(\roman*)}]
\item the dimensions $\dim L[k]$ are uniformly bounded;
\item for every finite $F\subset L$, there is $C_F>0$ such that
      $\dim V_j(F)\leq C_Fj$ for all $j\geq1$;
\item $L$ has polynomial commutator growth.
\end{enumerate}
When these conditions hold,
\[
 \dim L[k]=\dim L[1]\qquad(k\ne0).
\]
\end{theorem}

The implication from polynomial growth to bounded multiplicities is the
new point. After restriction to a suitable Virasoro subalgebra, a
nontrivial simple highest- or lowest-weight Virasoro subquotient contains
linearly independent Poincar\'e--Birkhoff--Witt (PBW) vectors indexed by
partitions. Their lifts are iterated brackets from one fixed finite
subset of $L$, with lengths controlled by the partition sizes. The
resulting growth exceeds every polynomial. The classification
of simple Virasoro weight modules with finite-dimensional weight spaces
then leaves only intermediate-series subquotients, whose nonzero
integral weight spaces are one-dimensional, and trivial subquotients.
An induction on finite-dimensional generating spaces forces constant
nonzero weight multiplicities.
No simplicity assumption is needed in this growth theorem.

For an algebra satisfying the original conditions SC1--SC3 with its
integral $L_0$-grading and finite-dimensional homogeneous spaces,
Theorem~\ref{thm:growth-rigidity} gives uniform boundedness, so
Definition~\ref{def:superconformal} applies. Theorem~\ref{thm:synthesis} then
proves exhaustion for the original graded class; see
Corollary~\ref{gr:cor:original}. The classification theorem itself
has no integral-spectrum assumption.

The polynomial hypothesis cannot simply be omitted. A finitely generated
Lie algebra has \emph{exponential commutator growth} if, for one finite
generating set $S$, there are constants $c,C>0$ such that
\[
 cj\leq\log\bigl(1+\dim V_j(S)\bigr)\leq Cj
 \qquad\text{for all sufficiently large }j.
\]
Changing the finite generating set changes bracket lengths by at most
a constant factor, so this condition is independent of $S$.

\begin{theorem}\label{thm:unbounded-example}
There is a finitely generated simple complex Lie algebra $L$ of
exponential commutator growth containing the centerless Virasoro
algebra such that
\[
 L=\bigoplus_{k\in\Z}L[k],\qquad \dim L[k]<\infty\quad(k\in\Z),
\]
but the dimensions $\dim L[k]$ are unbounded. More precisely, if $p(n)$
denotes the partition number, then the example may be chosen so that
\[
 \dim L[0]=2,\qquad
 \dim L[\pm(3+n)]\geq 1+p(n)\qquad(n\geq0).
\]
\end{theorem}

Within the integral setting above, polynomial growth and uniform
boundedness therefore give the same class, while finite-dimensional
weight spaces alone allow further simple algebras. Theorem~\ref{thm:growth-rigidity}
is proved in Subsection~\ref{sec:growth}, and Theorem~\ref{thm:unbounded-example}
is proved in Section~\ref{sec:unbounded}.

\paragraph{Outline of the proof.}
The proof of Theorem~\ref{thm:synthesis} has three steps.
Sections~\ref{sec:residue}--\ref{sec:log} first construct a finite
conformal model from the finite family of local fields supplied by
\cite[Theorem~2.7]{Billig}. A nilpotent operator on the finite-dimensional
fibers of the covering Virasoro module descends to the nonzero weight
spaces of $L$ and extends uniquely to weight zero. This gives a derivation
of $L$. Allowing the
fields to depend polynomially on $\log z$ then gives a common translation
operator. Locality, the Virasoro action, and the uniform bound on weight
multiplicities imply that the resulting conformal algebra is finitely
generated over $\C[\dd]$.

Sections~\ref{sec:coeff}--\ref{sec:untwist} next recover $L$ from the
logarithmic coefficient algebra, with kernel equal to its full center.
Simplicity of $L$ forces the conformal model to be semisimple and the
twisting automorphism to act transitively on its minimal ideals. Current ideals are
excluded, so these ideals are copies of one noncurrent finite simple
conformal superalgebra. The derivation responsible for the logarithms
is inner; a finite exponential change of fields removes the logarithms.
Thus
\begin{equation}\label{eq:intro-reduction}
 L\simeq \Alg(M^{\oplus r},\sigma)/
 Z\bigl(\Alg(M^{\oplus r},\sigma)\bigr),
\end{equation}
where $M$ is noncurrent, finite, and simple, while $\sigma$ is
diagonalizable and permutes the $r$ copies transitively
(Theorem~\ref{thm:ordinary-model}).
No integrality of the $L_0$-weights is used in these two reductions.

Finally, Section~\ref{sec:cyclic} reduces the cyclic permutation of the
$r$ copies to an automorphism of $M$, keeping track of the change in mode
indices and the phase determined by conformal weights.
Sections~\ref{sec:endpoints}--\ref{sec:chengkac} identify the resulting
coefficient algebras for the Witt, special, contact, and Cheng--Kac
families, including the exceptional small ranks.
Section~\ref{sec:final-reduction} combines these results to prove
Theorem~\ref{thm:synthesis}.

The growth argument is independent of these reductions.
Subsection~\ref{sec:growth} proves that polynomial commutator growth
forces bounded weight multiplicities in the integral setting. This
allows the classification to be applied under the original
Kac--van de Leur hypotheses. Section~\ref{sec:unbounded} constructs a
finitely generated simple example with unbounded multiplicities and
exponential commutator growth, showing why finite-dimensional weight
spaces alone do not suffice.

\section{The residue derivation}\label{sec:residue}
A superconformal Lie algebra is \emph{pluperfect} if
\begin{equation}\label{eq:pluperfect}
 L[0]=\sum_{k\ne0}[L[k],L[-k]].
\end{equation}
Every simple superconformal algebra is centerless, since it contains the
nonabelian algebra $\Vir$, and is pluperfect by
\cite[Lemma 2.6(b)]{Billig}.

A formal distribution is a series
$a(z)=\sum_{m\in\alpha+\Z}a_mz^{-m-1}$. Two distributions are
\emph{local} if $(z-w)^N[a(z),b(w)]=0$ for some nonnegative integer
$N$. A family is \emph{spanning} if its coefficients span $L$.
By \cite[Theorem 2.7]{Billig}, every pluperfect superconformal algebra
has a finite spanning family of homogeneous mutually local distributions.

Set $A=\C[t,t^{-1}]$. An $A\Vir$-module is both an $A$-module and a $\Vir$-module, with the actions satisfying the Leibniz rule. The cover in \cite[Theorems 2.14 and 2.16]{Billig} has finite fibers
$U^\beta$, where $\beta$ ranges over finitely many classes in $\C/\Z$, and a homomorphism
\[
 \pi:\widehat L=\bigoplus_\beta\bigoplus_{k\in\beta+\Z}t^k\otimes U^\beta
 \longrightarrow L,\qquad \im\pi=[\Vir,L].
\]
For each fiber $U$ there are endomorphisms $e_j\in\End U$,
$j\geq0$, only finitely many of which are nonzero.
Writing $I$ for the identity on the relevant vector space, they satisfy
\begin{equation}\label{eq:jet}
\begin{gathered}
 [e_i,e_j]=(j-i)e_{i+j},\
 R(n)=\sum_{j\ge0}\frac{n^{j+1}}{(j+1)!}e_j,\qquad
 L_n(t^k\otimes u)=-t^{k+n}\otimes(kI+R(n))u.
\end{gathered}
\end{equation}
We call a module with this action a finite-jet module.
For an endomorphism $B$ of a finite-dimensional space,
write $B=B_{\sm}+B_{\nil}$ for its commuting semisimple and nilpotent
Jordan components. In particular, write $e_0=H_0+J$.
The construction in this section uses the chosen Virasoro embedding.
The operator $J$ first descends to every nonzero weight space. Its
compatibility with brackets of nonzero weights is established next, and
pluperfectness together with centerlessness determines the extension to
weight zero.

\subsection{Finite-dimensional invariance lemmas}
\begin{lemma}\label{lem:triangular}
For the operators in~\eqref{eq:jet},
\[
 [H_0,e_j]=je_j,\qquad [J,e_j]=0.
\]
For every finite sum $B=e_0+\sum_{j\ge1}c_je_j$, there is an invertible operator $S$, commuting with $J$, such that $BS=Se_0$. In particular $B_{\nil}=J$.
\end{lemma}
\begin{proof}
The commuting operators $\ad H_0$ and $\ad J$ on $\End U$ are respectively semisimple and nilpotent. Since $e_j$ is an eigenvector of $\ad e_0$ of eigenvalue $j$, it has $\ad H_0$-weight $j$ and is annihilated by $\ad J$.

Let $A_j=c_je_j$, $S_0=I$, and
\[
 S_d=-\frac1d\sum_{j=1}^d A_jS_{d-j}\qquad(d\ge1).
\]
Then $[e_0,S_d]=dS_d$ and $[J,S_d]=0$. Thus $S_d$ raises $H_0$-weight by $d$ and vanishes for sufficiently large $d$. The finite sum $S=\sum_dS_d$ is invertible, since $S-I$ strictly raises weight. The recurrence gives $BS=Se_0$. Since $S$ commutes with $J$, the nilpotent part of $B$ is $J$.
\end{proof}

\begin{lemma}\label{lem:commuting-difference}
Let $S,N\in\End E$, where $E$ is finite-dimensional, $S$ is semisimple, $N$ is nilpotent, and $[S,N]=0$. Suppose
\[
 A(x)=I+x^{-1}(S+N)+\sum_{j\ge2}A_jx^{-j}
\]
is the Laurent expansion of a rational matrix at infinity, and every $A_j$ commutes with $S$ and $N$. If a rational vector $y(x)$ satisfies $y(x+1)=A(x)y(x)$, then $Ny(x)=0$.
\end{lemma}
\begin{proof}
Every eigenspace of $S$ is invariant under the recurrence. Consider a nonzero component on which $S=sI$, and write
\[
 y(x)=x^m\sum_{i\ge0}y_ix^{-i},\qquad m\in\Z,\quad y_0\ne0.
\]
Comparison of the coefficient of $x^{m-1}$ gives
\[
 (m-s)y_0=Ny_0.
\]
Nilpotence implies $m=s$ and $Ny_0=0$. For $j\ge1$, comparison of the coefficient of $x^{m-j-1}$ gives
\begin{equation}\label{eq:coefficient-induction}
 (-jI-N)y_j
 =\sum_{a=2}^{j+1}A_a y_{j+1-a}
 -\sum_{i=0}^{j-1}\binom{m-i}{j+1-i}y_i.
\end{equation}
Assume $Ny_i=0$ for $i<j$. The right side is annihilated by $N$, because $N$ commutes with every $A_a$. Applying $N$ to~\eqref{eq:coefficient-induction} gives
\[
 (-jI-N)Ny_j=0.
\]
The operator $-jI-N$ is invertible, so $Ny_j=0$. Induction proves that every Laurent coefficient is annihilated by $N$. Hence the rational vector $Ny(x)$ is zero. Apply this argument to each eigenspace of $S$.
\end{proof}

\begin{lemma}\label{lem:rational}
Let $B\in\End U$, and let $W(x)$ be a rational map to the Grassmannian
$\operatorname{Gr}_d(U)$ of $d$-dimensional subspaces of $U$, satisfying
\[
 W(x+1)=(xI+B)W(x)
\]
at generic $x$. Then $B_{\nil}W(x)\subseteq W(x)$ at every regular point of the map.
\end{lemma}
\begin{proof}
The cases $d=0$ and $d=\dim U$ are immediate. We use ordinary exterior powers of the underlying vector space $U$; all operators in this argument are even. Write $B=B_{\sm}+N$. Choose a nonzero rational decomposable vector $w(x)\in\bigwedge^dU$ representing $W(x)$. There is $c(x)\in\C(x)^\times$ such that
\begin{equation}\label{eq:plucker-rec}
 w(x+1)=c(x)\bigwedge^d(xI+B)w(x).
\end{equation}
Writing $w(x)=x^m(w_0+O(x^{-1}))$, with $w_0\ne0$, shows from the leading terms that
\[
 c(x)x^d=1+c_1x^{-1}+O(x^{-2}).
\]
For $T\in\End U$, denote its additive action on $\bigwedge^dU$ by $T^{[d]}$. The matrix on the right of~\eqref{eq:plucker-rec} has expansion
\[
 I+x^{-1}\bigl((B_{\sm})^{[d]}+c_1I+N^{[d]}\bigr)+O(x^{-2}).
\]
Every coefficient commutes with $(B_{\sm})^{[d]}$ and $N^{[d]}$: the base matrix $I+x^{-1}B$ commutes with $\exp(tB_{\sm})$ and $\exp(tN)$, and differentiation of the induced exterior actions gives the asserted commutation. The operator $(B_{\sm})^{[d]}+c_1I$ is semisimple and $N^{[d]}$ is nilpotent. Lemma~\ref{lem:commuting-difference} therefore gives
\[
 N^{[d]}w(x)=0.
\]
For a nonzero decomposable vector representing a subspace $W$, this implies $N(W)\subseteq W$. Indeed, the component of its additive exterior action with one factor in $U/W$ is the map induced by $N:W\to U/W$; its vanishing forces that map to be zero. Finally, $N$-invariance is a closed condition on the Grassmannian, so the conclusion holds at every regular point of $W(x)$.
\end{proof}

\subsection{Descent on nonzero weight spaces}
A punctured invariant family means subspaces $K_k\subseteq U$, indexed by $k\in(\beta+\Z)\setminus\{0\}$, such that
\begin{equation}\label{eq:kernel-inc}
 (kI+R(n))K_k\subseteq K_{k+n}
 \quad\text{whenever }k\ne0\text{ and }k+n\ne0.
\end{equation}
No space $K_0$ is specified.

\begin{proposition}\label{prop:punctured-descent}
For every punctured invariant family satisfying~\eqref{eq:kernel-inc},
\begin{equation}\label{eq:Jkernel}
 JK_k\subseteq K_k\qquad(k\ne0).
\end{equation}
\end{proposition}
\begin{proof}
Fix an arbitrary nonzero $k_0\in\beta+\Z$. The two polynomial matrices
\[
 k_0I+R(n),\qquad (k_0+n)I+R(-n)
\]
are invertible for all but finitely many integers $n$, since both equal $k_0I$ at $n=0$. Their inclusions in opposite directions force equality of the two fiber dimensions and hence
\[
 K_{k_0+n}=(k_0I+R(n))K_{k_0}
\]
for all but finitely many integers $n$. Define the rational Grassmannian map
\begin{equation}\label{eq:anchored-family}
 W_{k_0}(x)=(k_0I+R(x-k_0))K_{k_0}.
\end{equation}
It is regular at $x=k_0$ and has value $K_{k_0}$ there. It agrees with the actual fibers on a cofinite subset of the lattice. The unit-step map $xI+R(1)$ is invertible generically. Hence, at infinitely many lattice points,
\[
 W_{k_0}(x+1)=(xI+R(1))W_{k_0}(x).
\]
Both sides are rational maps, so this is a rational identity. Lemma~\ref{lem:triangular} gives $R(1)_{\nil}=J$, and Lemma~\ref{lem:rational}, evaluated at the regular point $k_0$, proves~\eqref{eq:Jkernel}. Since $k_0$ was arbitrary, the assertion holds at every nonzero fiber.
\end{proof}

Let $M[k]=U/K_k$ and $N[k]=V/K'_k$ be punctured quotient families of finite-jet modules, with the same physical-mode coset. Their induced residues are denoted by $J_M,J_N$.

\begin{corollary}\label{cor:punctured-natural}
If maps $f_k:M[k]\to N[k]$ intertwine every Virasoro action whose source and target weights are nonzero, then
\begin{equation}\label{eq:punctured-natural}
 J_Nf_k=f_kJ_M\qquad(k\ne0).
\end{equation}
No surjectivity is assumed.
\end{corollary}
\begin{proof}
Let $q_k:U\to M[k]$ and $q'_k:V\to N[k]$ be the quotient maps. In $U\oplus V$, consider
\[
 \mathcal K_k=\{(u,v): f_k(q_k(u))=q'_k(v)\}.
\]
The intertwining identities make this a punctured invariant family for the direct sum of the two jet modules. Proposition~\ref{prop:punctured-descent} gives its invariance under $J_U\oplus J_V$. Applying that invariance to any $(u,v)\in\mathcal K_k$ proves the displayed identity. The same argument works for homogeneous maps of either parity, since the Virasoro action is even.
\end{proof}

Apply Proposition~\ref{prop:punctured-descent} to the kernel fibers of $\pi$. Since $\im\pi=[\Vir,L]$ contains every $L[k]$ with $k\ne0$, this gives an even nilpotent operator $D$ on all nonzero weight spaces, defined by
\[
 D\pi(t^k\otimes u)=\pi(t^k\otimes Ju).
\]
It intertwines all Virasoro actions that stay at nonzero weights. The one-dimensional jet module with $e_0=-1$ and $e_j=0$ for $j>0$ maps to $L$ by $1_k\mapsto L_k$. Its residue is zero. Corollary~\ref{cor:punctured-natural} therefore gives $D(L_k)=0$ for $k\ne0$.

\subsection{Compatibility with the bracket}
At a fixed total mode, locality expresses the bracket as a rational
function of one mode index. We separate its two possible poles from
its polynomial part. Punctured naturality controls the pole terms;
a finite tensor quotient controls the polynomial part.
Fix two finite-jet fibers $U,V$, with mode cosets
$\beta_U+\Z$ and $\beta_V+\Z$. Write
$u_p=\pi(t^p\otimes u)$ and $v_q=\pi(t^q\otimes v)$.
For a generalized $e_0$-eigenvector $u$ of eigenvalue $a$, the fields in
\cite[equation~(7), Proposition~4.1]{Billig} are
\[
 B_u(z)=
 \begin{cases}
  \displaystyle\sum_p u_pz^{-p},&a\ne1,\\[2mm]
  \displaystyle\sum_p pu_pz^{-p},&a=1,
 \end{cases}
\]
where the sum is over the mode coset of $u$. They are mutually local,
and
\begin{equation}\label{eq:weighted-local-fields}
 \sum_p pu_pz^{-p-1}=
 \begin{cases}
  -\dd_z B_u,&a\ne1,\\
  z^{-1}B_u,&a=1.
 \end{cases}
\end{equation}
Differentiation and multiplication by Laurent monomials preserve
locality, so these weighted fields are mutually local for arbitrary
labels, by linearity.

For a sequence $f(p)$, set $\Delta_p f(p)=f(p)-f(p+1)$, keeping
the total mode $k$ fixed. Finiteness of the fibers and the number of
mode cosets gives a single integer $N\ge3$ such that
\[
 \Delta_p^N\bigl(p(k-p)[u_p,v_{k-p}]\bigr)=0
\]
for every pair of labels and every total mode $k$. By \cite[Lemmas 3.2 and 3.3(a)]{Billig}, the expression in parentheses is a polynomial in $p$ of degree less than $N$.

Denote this unique polynomial continuation by $P_k(p)\in\operatorname{Hom}(U\otimes V,L[k])[p]$. For $k\ne0$, define
\[
 F_k(p)=\frac{P_k(p)}{p(k-p)}.
\]
Its partial-fraction decomposition has the form
\begin{equation}\label{eq:pole-decomp}
 F_k(p)=H_k(p)+\frac{A_k}{p}+\frac{B_k}{k-p},
\end{equation}
where $A_k,B_k:U\otimes V\to L[k]$ and $H_k(p)$ is polynomial. The degrees of all $H_k$ are less than $N$, uniformly in $k$. At lattice points with $p(k-p)\ne0$, $F_k(p)(u,v)=[u_p,v_{k-p}]$.

For $s=k+n$, with $k,s\ne0$, the super Jacobi identity with the even element $L_n$ gives the following rational identity in $p$:
\begin{align}
 L_nF_k(p)={}&-F_s(p+n)((pI+R_U(n))\otimes I)\notag\\
 &-F_s(p)(I\otimes((k-p)I+R_V(n))).\label{eq:bracket-cov}
\end{align}
It holds as a rational identity because it holds at all but finitely many lattice points.

Fix $s\ne0$ and let $n$ vary through integers, with $k=s-n$. Away from finitely many pole collisions, comparison at the new pole $p=-n$ gives
\[
 A_s((R_U(n)-nI)\otimes I)=0.
\]
The coefficient map $A_s$ is fixed in this comparison. Polynomial comparison in $n$ therefore yields
\begin{equation}\label{eq:Aann}
 A_s((e_0^U-I)\otimes I)=0,\qquad
 A_s(e_j^U\otimes I)=0\quad(j\ge1).
\end{equation}
In particular, $A_s(J_U\otimes I)=0$: on every generalized $e_0$-eigenspace except that of eigenvalue $1$, $e_0-I$ is invertible, while on the generalized $1$-eigenspace it equals $J_U$. The new pole $p=s$ gives, in the same way,
\[
 B_s(I\otimes(e_0^V-I))=0,\qquad
 B_s(I\otimes e_j^V)=0\quad(j\ge1),\qquad
 B_s(I\otimes J_V)=0.
\]
At the old poles $p=0,k$, respectively, one obtains
\begin{align}
 L_nA_k&=-A_s(I\otimes(kI+R_V(n))),\label{eq:Aintertwine}\\
 L_nB_k&=-B_s((kI+R_U(n))\otimes I).\label{eq:Bintertwine}
\end{align}
For $n=0$ these formulas are also immediate from the physical weights.

For fixed $u\in U$, the maps $v\mapsto A_k(u,v)$ intertwine the jet module with fiber $V$, now placed on the total mode coset, with the punctured quotient family $L[k]$. Such a jet module exists on any mode coset by the same formula~\eqref{eq:jet}. Corollary~\ref{cor:punctured-natural} and~\eqref{eq:Aann} give
\[
 DA_k=A_k(J_U\otimes I+I\otimes J_V).
\]
The same conclusion holds for $B_k$.

Each pole term satisfies~\eqref{eq:bracket-cov} separately. For the first pole, the right-hand side is
\begin{align*}
 &-\frac{A_s((pI+R_U(n))\otimes I)}{p+n}
 -\frac{A_s(I\otimes((k-p)I+R_V(n)))}p\\
 &\hspace{12mm}=-A_s-\frac{A_s(I\otimes(kI+R_V(n)))}p+A_s
 =\frac{L_nA_k}{p}.
\end{align*}
Here~\eqref{eq:Aann} gives $A_s(R_U(n)\otimes I)=nA_s$. Similarly, the second pole contributes
\begin{align*}
 &-\frac{B_s((pI+R_U(n))\otimes I)}{k-p}
 -\frac{B_s(I\otimes((k-p)I+R_V(n)))}{s-p}\\
 &\hspace{12mm}=-\frac{B_s((kI+R_U(n))\otimes I)}{k-p}+B_s-B_s
 =\frac{L_nB_k}{k-p}.
\end{align*}
Consequently $H_k$ satisfies~\eqref{eq:bracket-cov} with $F$ replaced by $H$.

Choose a single integer $d\ge N$. The tensor product of the two full jet modules can be realized with variables $t_1,t_2$. Let $q=t_1/t_2$, and quotient its scalar factor by $(q-1)^d$. This ideal is invariant under the diagonal Virasoro action: its differential part in $q$ is a multiple of $q(q^n-1)\partial_q$. With $x=\log q$ in the quotient, the finite fiber is
\[
 B_d\otimes U\otimes V,\qquad B_d=\C[x]/(x^d).
\]
Fix a representative $p_0$ of the first mode coset, so $p=p_0+j$ with $j\in\Z$. The factor $e^{p_0x}$ is a finite polynomial in the truncated ring $B_d$. Its jet polynomial is
\begin{equation}\label{eq:tensor}
 R_{U,V,d}(n)=(e^{nx}-1)(\partial_x+p_0)
   +e^{nx}R_U(n)\otimes I+I\otimes R_V(n).
\end{equation}
The total Laurent coordinate is $t_2$. Although $\partial_x$ alone is not a derivation of $B_d$, its product with $e^{nx}-1$ is well defined. This construction is a quotient of the diagonal Virasoro action, so it satisfies the Virasoro relations.

Its Euler operator is
\[
 e_0^{\rm ten}=x\partial_x+p_0x+e_0^U\otimes I+I\otimes e_0^V.
\]
The multiplication operator $e^{p_0x}$ on $B_d$ satisfies
\[
 e^{p_0x}(x\partial_x+p_0x)e^{-p_0x}=x\partial_x.
\]
The operator on the right is semisimple and acts on a separate tensor factor. This conjugation commutes with the two fiber actions. Hence
\begin{equation}\label{eq:tensor-nil}
 (e_0^{\rm ten})_{\nil}=J_U\otimes I+I\otimes J_V.
\end{equation}
For every $k\ne0$, define the linear map
\[
 \Phi_k:B_d\otimes U\otimes V\longrightarrow L[k]
\]
by
\[
 \Phi_k(e^{jx}\otimes u\otimes v)=H_k(p_0+j)(u,v)\qquad(j\in\Z).
\]
It exists uniquely: if $H_k(p_0+j)=\sum_{a=0}^{d-1}j^aH_{k,a}$, then
$\Phi_k(x^a\otimes u\otimes v)=a!H_{k,a}(u,v)$. Thus the same finite jet module works for every $k$.

For $p=p_0+j$, the required intertwining identity follows from the identity in $B_d\otimes U\otimes V$:
\begin{align*}
 &(kI+R_{U,V,d}(n))(e^{jx}\otimes u\otimes v)\\
 &\quad=e^{(j+n)x}\otimes((pI+R_U(n))u)\otimes v
 +e^{jx}\otimes u\otimes((k-p)I+R_V(n))v.
\end{align*}
The possible top-degree discrepancy in differentiating the truncated exponential is killed by $e^{nx}-1$. Applying $\Phi_{k+n}$ and the covariance identity for $H$ proves that $\Phi_k$ is a punctured intertwining family. Naturality and~\eqref{eq:tensor-nil} give
\[
 D\Phi_k=\Phi_k(J_U\otimes I+I\otimes J_V).
\]
Together with the two pole terms, this proves
\begin{equation}\label{eq:nonzero-Leibniz}
 D[x,y]=[Dx,y]+[x,Dy]
\end{equation}
when both weights and their sum are nonzero.

\subsection{Extension through weight zero}
\begin{lemma}\label{lem:tail}
Suppose $L$ is pluperfect. If $w\in L[0]$ centralizes all but finitely many weight spaces, then $w\in Z(L)$.
\end{lemma}
\begin{proof}
For $x\in L[b]$, $b\ne0$, choose a sufficiently large positive integer $n$ so that $w$ centralizes $L_n,L_{-n}$ and the spaces containing $[L_n,x]$ and $[L_{-n},x]$. Then
\[
 [w,[L_n,[L_{-n},x]]]=[w,[L_{-n},[L_n,x]]]=0.
\]
Their difference is $[w,[2nL_0,x]]=-2nb[w,x]$. Thus $[w,x]=0$. Pluperfectness gives $L[0]$.
\end{proof}

\begin{theorem}\label{thm:residue}
Let $L$ be a simple superconformal Lie algebra, with its chosen Virasoro subalgebra and finite $A\Vir$-cover. The nilpotent parts $J$ of the fiber Euler operators induce a unique even nilpotent derivation $D$ of $L$ such that
\[
 D\pi(t^k\otimes u)=\pi(t^k\otimes Ju)\qquad(k\ne0),
 \qquad D(L_n)=0\quad(n\in\Z).
\]
If $J^r=0$ on every fiber, then $D^{2r-1}=0$ on $L$. The same conclusion holds for a centerless pluperfect superconformal algebra with the stated finite-cover and locality properties.
\end{theorem}
\begin{proof}
The preceding arguments define $D$ on the nonzero weight spaces and prove
\eqref{eq:nonzero-Leibniz} when both weights and their sum are nonzero.
For $z\in L[0]$, choose a finite expression
\[
 z=\sum_i[x_i,y_i],\qquad x_i\in L[a_i],\quad y_i\in L[-a_i],\quad a_i\ne0,
\]
and define
\begin{equation}\label{eq:Dzero}
 Dz=\sum_i\bigl([Dx_i,y_i]+[x_i,Dy_i]\bigr).
\end{equation}
If the expression for $z$ is zero, let $w$ be the right side. For homogeneous $u\in L[b]$ with $b\ne0,\pm a_i$, expand $[[x_i,y_i],u]$ by the super Jacobi identity and apply~\eqref{eq:nonzero-Leibniz} to its two nested brackets. All intermediate weights are nonzero. Summing gives $[w,u]=0$. Lemma~\ref{lem:tail} makes $w$ central, so $w=0$ because $L$ is centerless. Thus~\eqref{eq:Dzero} is well defined. It also proves the derivation identity for opposite nonzero weights. Since $[L_1,L_{-1}]=2L_0$, it gives $D(L_0)=0$.

For fixed $z\in L[0]$, the same Jacobi expansion gives
\begin{equation}\label{eq:generic-zero-defect}
 D[z,u]=[Dz,u]+[z,Du]
\end{equation}
for all nonzero weights $b$ except a finite set depending on the chosen expression for $z$.

For $z\in L[0]$, define the Virasoro-equivariance defect
\[
 E_n(z)=D[L_n,z]-[L_n,Dz].
\]
Equation~\eqref{eq:generic-zero-defect} gives $E_n(z)=0$ for sufficiently large positive and negative $n$. If $m,n,m+n\ne0$, the Jacobi identity and punctured equivariance give
\[
 (m-n)E_{m+n}(z)=[L_m,E_n(z)]-[L_n,E_m(z)].
\]
Fix $m\ne0$ and let $n$ be sufficiently large in both directions. It follows that $[L_n,E_m(z)]=0$. Applying this with $n$ and $-n$ and using $[L_n,L_{-n}]=2nL_0$ forces $E_m(z)=0$, since it has physical weight $m$. The case $m=0$ is immediate. Together with~\eqref{eq:Dzero} when a nonzero weight is mapped to weight zero, this proves $[D,\ad L_n]=0$ on all of $L$.

For arbitrary homogeneous $a,b$, write
\[
 \mathcal E(a,b)=D[a,b]-[Da,b]-[a,Db].
\]
The Virasoro equivariance just proved implies
\[
 [L_n,\mathcal E(a,b)]
 =\mathcal E([L_n,a],b)+\mathcal E(a,[L_n,b]).
\]
Take $a=z\in L[0]$ and $b=u\in L[b]$, $b\ne0$. For sufficiently large positive and negative $n$, the first term on the right vanishes by~\eqref{eq:nonzero-Leibniz}, and the second by~\eqref{eq:generic-zero-defect}. Thus $\mathcal E(z,u)$ is annihilated by $L_n,L_{-n}$ and has nonzero weight. It must be zero.

For $z,w\in L[0]$, expanding the super Jacobi identity for $z,w,u$, with $u\in L[b]$, $b\ne0$, and applying the already established cases of the derivation identity gives
\[
 [\mathcal E(z,w),u]=0.
\]
Pluperfectness and centerlessness then give $\mathcal E(z,w)=0$. Thus $D$ is a derivation on all of $L$.

Choose $r$ such that $J^r=0$ on every fiber. Then $D^r=0$ on nonzero weights. For opposite nonzero weights, the iterated derivation rule gives
\[
 D^{2r-1}[x,y]
 =\sum_{j=0}^{2r-1}\binom{2r-1}{j}[D^jx,D^{2r-1-j}y]=0.
\]
Every summand is zero. Pluperfectness proves the same bound at weight zero.
It also proves uniqueness: a derivation with the prescribed action on
the nonzero weight spaces is determined on their opposite-weight brackets.
\end{proof}

\section{Logarithmic fields and finite generation}\label{sec:log}
An ordinary derivation of a Lie conformal superalgebra commutes with
$\dd$ and satisfies the graded Leibniz rule for every $j$-product. It is
\emph{inner} if it equals $h_{(0)}$ for some state $h$.

Let $D$ be the residue derivation of Theorem~\ref{thm:residue}, and set
\[
 T=\ad L_{-1},\qquad D_z=\dd_z+z^{-1}\dd_\zeta,
\]
where $\zeta$ is a formal logarithm of $z$. The derivation $D$ is even,
preserves physical weights, and satisfies $D^{s_0}=0$ for some $s_0$,
$D|_{\Vir}=0$, and $D(u_k)=(Ju)_k$ at every nonzero cover mode.
In particular, $D$ commutes with $T$. All fields below have coefficients in $L$.
The variables $z$ and $\zeta$ are formally independent.
For mutually local fields, $a_{(n)}b$ denotes the nonnegative
field product of \cite[Definition~3.3]{Bakalov}.

\subsection{Local covariant fields}
For $u\in U^\beta$, write
\[
 u_k=\pi(t^k\otimes u),\qquad
 F_u(z)=\sum_{k\in\beta+\Z}u_kz^{-k-1},\qquad
 G_u(z)=\sum_{k\in\beta+\Z}ku_kz^{-k-1}.
\]
The zero coefficient of $G_u$ vanishes. Hence Theorem~\ref{thm:residue} gives
\begin{equation}\label{eq:DG}
 DG_u=G_{Ju}.
\end{equation}

Write
\[
 B_-=e_0-\frac12e_1+\frac16e_2-\cdots.
\]
Lemma~\ref{lem:triangular} gives $(B_-)_{\nil}=J$, and the jet formula gives
\begin{equation}\label{eq:translation-F}
 TF_u=\dd_zF_u+z^{-1}F_{B_-u},
\end{equation}
\begin{equation}\label{eq:translation-G}
 z(T-\dd_z)G_u=G_{(B_--I)u}+F_{(B_--I)u}.
\end{equation}
The fields $G_u$ are mutually local by
\eqref{eq:weighted-local-fields}. Let
$\ell(z)=\sum_nL_nz^{-n-2}$. At fixed total physical mode $s$, the
jet formula gives
\[
 [L_n,(s-n)u_{s-n}]
 =-(s-n)\bigl(((s-n)I+R(n))u\bigr)_s.
\]
This is polynomial in $n$, of uniformly bounded degree, so $\ell$
is local with every $G_u$. Since
$\ell_{(0)}G_u=[L_{-1},G_u]=TG_u$, Dong's lemma also gives locality
of the fields $TG_u$ with one another and with the $G_u$.

Define
\[
 E_u=z(T-\dd_z)G_u-G_{(B_--I)u}.
\]
Then $E_u=F_{(B_--I)u}$. Since $D$ commutes with $T$, and $J$ commutes with $B_-$,
\begin{equation}\label{eq:DE}
 DE_u=E_{Ju}.
\end{equation}
Thus $DF_v=F_{Jv}$ for $v\in(B_--I)U$, including the zero mode.
On a generalized $B_-$-eigenspace write
\[
 B_-=\omega I+J.
\]
If $\omega\ne1$, the operator $B_--I$ is invertible, so every $F_u$ in that block is local and satisfies
\[
 DF_u=F_{Ju}.
\]
If $\omega=1$, equation~\eqref{eq:DE} gives the same identity for $F_{Ju},F_{J^2u},\ldots$.

For $\omega\ne1$ define
\begin{equation}\label{eq:Xfield}
 X_u(z,\zeta)=z^\omega\sum_{j\ge0}\frac{\zeta^j}{j!}F_{J^ju}(z).
\end{equation}
For $\omega=1$, begin with the same formal expression and use its total derivative
\begin{equation}\label{eq:Yfield}
 Y_u=D_zX_u
 =\sum_{j\ge0}\frac{\zeta^j}{j!}\left(F_{J^{j+1}u}-G_{J^ju}\right).
\end{equation}
All sums are finite. Equations~\eqref{eq:translation-F}--\eqref{eq:DE} give, on the respective blocks,
\begin{equation}\label{eq:common-residue}
 TX_u=D_zX_u,\quad \dd_\zeta X_u=DX_u \qquad (\omega\ne1),
\end{equation}
and
\begin{equation}\label{eq:common-residue-exceptional}
 TY_u=D_zY_u,\quad \dd_\zeta Y_u=DY_u \qquad (\omega=1).
\end{equation}
For $\omega=1$ we use $Y_u$: its $F$-terms are indexed by vectors in
$\operatorname{im}(B_--I)$, and its $G$-terms have vanishing zero coefficient.
Thus neither requires descent through the zero fiber of the cover.
The logarithm-free coefficient of $Y_u$ at physical mode $k$ is
$((J-kI)u)_k$, which spans the same image as $u_k$ whenever $k\ne0$.
Hence finitely many fields of the forms~\eqref{eq:Xfield}--\eqref{eq:Yfield},
together with the Virasoro field, span every nonzero weight space through
their coefficients.

\subsection{The logarithmic conformal closure}
For a faithful restricted realization, set
\[
 \mathfrak p=\bigoplus_{\operatorname{Re}k\geq0}L[k],
 \qquad W=U(L)\otimes_{U(\mathfrak p)}\C,
\]
where $\mathfrak p$ acts trivially on $\C$. The subspace $\mathfrak p$ is
a subalgebra. The Poincar\'e--Birkhoff--Witt theorem identifies $W$, as a
vector space, with the enveloping algebra of the sum of the negative
real-weight spaces. In particular, $W\ne0$ and $L_{-1}$ acts nontrivially.
The kernel of this representation is an ideal, so simplicity makes it
faithful. Write $\rho:L\to\End W$ for this representation and
$\mathbf1=1\otimes1$ for its vacuum vector.
For a monomial $x_1\cdots x_t\mathbf1$ of negative-weight
elements, a mode of sufficiently large real weight annihilates that
monomial: commute it successively to the right, observing that its
weight remains nonnegative after any subset of the finitely many
negative shifts. Thus every homogeneous field above is a restricted
logarithmic quantum field on $W$. Since $D$ preserves $\mathfrak p$, its
induced operator $\widehat D$ on $W$ satisfies
$[\widehat D,\rho(x)]=\rho(Dx)$.

Let $C$ be the smallest linear space containing the fields
\eqref{eq:Xfield}--\eqref{eq:Yfield} and $\ell$, and closed under $D_z$,
$\dd_\zeta$, and the products $a_{(n)}b$ for $n\geq0$. The identities
following \cite[Definition 3.3]{Bakalov} give sesquilinearity and the
derivation rules for $D_z$ and $\dd_\zeta$. Dong's lemma
\cite[Lemma 3.4]{Bakalov} preserves locality. To obtain the Jacobi
identity, place the initial fields and their finitely many $\dd_\zeta$-derivatives
inside the larger local space obtained by adjoining the identity and closing
under $D_\zeta=z\dd_z+\dd_\zeta$ and all integer products, as in
\cite[Theorem 3.7]{Bakalov}. That theorem makes the resulting vertex algebra
closed under the translation $D_z$. It is also $\dd_\zeta$-stable:
the chosen generators are stable, $\dd_\zeta$ commutes with $D_\zeta$,
and it is a derivation of every field product. Thus this larger space
contains the $D_z$-, $\dd_\zeta$-, and
nonnegative-product closure defining $C$. Hence the conformal Jacobi identity
restricts to $C$. The super version uses supercommutators throughout. Only this
singular-product subspace is used below.

\begin{lemma}\label{lem:finite-boundary}
Let $a,b$ be polynomial logarithmic $L$-valued local fields, and choose
$N\geq1$ with $(z-w)^N[a(z,\zeta),b(w,\eta)]=0$. For one exponent
component write
\[
 a(z,\zeta)=\sum_{i\in\Z}a_i(\zeta)z^{-i-\alpha_0},
 \qquad a_-=\sum_{i\geq1}a_i(\zeta)z^{-i-\alpha_0},
\]
where $-1<\operatorname{Re}\alpha_0\leq0$. Then
\begin{equation}\label{eq:finite-boundary}
 (z-w)^N[a_-(z,\zeta),b(w,\eta)]
 =\sum_{r=1-N}^{0}z^{-r-\alpha_0}Q_r(w;\zeta,\eta),
\end{equation}
where
\begin{equation}\label{eq:boundary-coefficients}
 Q_r(w;\zeta,\eta)=
 \sum_{j=0}^{r+N-1}(-1)^j\binom Nj
 w^j[a_{r+N-j}(\zeta),b(w,\eta)].
\end{equation}
Every coefficient of every nonnegative product $a_{(n)}b$ is a
finite sum of brackets of coefficients in $L$.
\end{lemma}
\begin{proof}
The coefficient of $z^{-r-\alpha_0}$ on the left of
\eqref{eq:finite-boundary} has index $i=r+N-j$.
For $r\leq-N$ there are no terms with $i\geq1$.
For $r\geq1$ all $j=0,\ldots,N$ occur, and their sum is the
corresponding coefficient of the zero full commutator. The remaining
indices are exactly $1-N\leq r\leq0$, giving
\eqref{eq:boundary-coefficients}.

For $0\leq n<N$, \cite[Proposition 3.6]{Bakalov} expresses $a_{(n)}b$
as
\begin{equation}\label{eq:singular-boundary-product}
 \left.
 \frac{D_z^{\,N-1-n}}{(N-1-n)!}
 \left((z-w)^N[a_-(z,\zeta),b(w,\eta)]\right)
 \right|_{z=w,\ \zeta=\eta}.
\end{equation}
There are only $N$ powers of $z$ in \eqref{eq:finite-boundary}.
For a polynomial $p(\zeta)$ and $q\geq0$,
\[
 \frac{D_z^q}{q!}\bigl(z^\gamma p(\zeta)\bigr)
 =z^{\gamma-q}\binom{\gamma+\dd_\zeta}{q}p(\zeta).
\]
Consequently differentiation and diagonal substitution leave only
finitely many summands in each coefficient. The brackets in
\eqref{eq:boundary-coefficients} lie in $L$, so the resulting
coefficients lie in $L$ as well. Products with $n\geq N$ vanish.
Finite sums of exponent components give the general case.
\end{proof}

\begin{proposition}\label{prop:log-closure}
The space $C$ is a Lie conformal superalgebra of $L$-valued polynomial
logarithmic fields. With $\dd=D_z$ and $\delta=\dd_\zeta$, every
$a\in C$ satisfies
\begin{equation}\label{eq:common-all}
 Ta=D_za,\qquad \delta a=Da.
\end{equation}
If $D^s=0$ on $L$, then
\begin{equation}\label{eq:log-degree}
 \delta^s C=0.
\end{equation}
Its physical-mode and conformal-weight gradings are direct, and
\begin{equation}\label{eq:Hdiag}
 \ell_{(0)}=\dd,\qquad
 \mathcal H=\ell_{(1)}+\delta
            =\ad L_0-z\dd_z
\end{equation}
is diagonalizable. Moreover $\delta\ell=0$.
\end{proposition}
\begin{proof}
Lemma~\ref{lem:finite-boundary} proves $L$-valuedness at every stage
of the singular-product closure; derivatives preserve it as well.
Coefficientwise $D$ is a derivation of the brackets in
\eqref{eq:boundary-coefficients}, commutes with the total derivatives,
and commutes with the cutoff and diagonal substitution. Hence it is
a derivation of all the nonnegative field products. The same holds
for coefficientwise $T$. On the field side, $D_z$ and $\dd_\zeta$ are
derivations of these products and commute with each other.
Thus $T-D_z$ and $D-\dd_\zeta$ annihilate the closure because they
annihilate its generators. This proves \eqref{eq:common-all};
applying it repeatedly proves \eqref{eq:log-degree}.
In particular,
\[
 a(z,\zeta)=e^{\zeta D}a(z,0).
\]

For $\Delta\in\C$ and $\beta\in\C/\Z$, let $C_{\Delta,\beta}$ be the
subspace of $C$ consisting of fields of the form
\[
 a(z,\zeta)=\sum_{k\in\beta+\Z}a_k(\zeta)z^{-k-\Delta},
 \qquad a_k(\zeta)\in L[k][\zeta].
\]
We call $\Delta$ the conformal weight and $\beta$ the physical-mode coset,
and write $C_\Delta=\sum_\beta C_{\Delta,\beta}$.
The derivative $D_z$ raises $\Delta$ by one, and $\dd_\zeta$
preserves both labels. In \eqref{eq:singular-boundary-product},
multiplication by $(z-w)^N$ raises total degree in $z,w$ by $N$,
and the indicated derivative lowers it by $N-1-n$. Brackets add
physical weights. Therefore
\begin{equation}\label{eq:field-grading}
 (C_{\Delta,\beta})_{(n)}C_{\Gamma,\gamma}
 \subseteq C_{\Delta+\Gamma-n-1,\beta+\gamma}.
\end{equation}
The initial $X_u$ have $\Delta=1-\omega$, the $Y_u$ have
$\Delta=1$, and $\ell$ has $\Delta=2$. Every element of $C$ is a
finite sum of homogeneous fields. The physical cosets are direct
because those of $L$ are direct; the $\Delta$-grading is direct
because $\ad L_0-z\dd_z$ has eigenvalue $\Delta$ on such a field.

The untwisted Virasoro commutator formula, applied to the field
$\ell(z)=\sum_nL_nz^{-n-2}$, gives
\[
 \ell_{(0)}a=[L_{-1},a],\qquad
 \ell_{(1)}a=[L_0,a]-z[L_{-1},a].
\]
Together with $Ta=D_za$ this gives \eqref{eq:Hdiag}. The equality
$D(L_n)=0$ proves $\delta\ell=0$.
\end{proof}

\subsection{Criteria for finite generation}
\begin{lemma}\label{lem:onesided}
Let $L$ be centerless, and let $\mathcal S$ be a family of fields whose
coefficients generate $L$ as a Lie superalgebra. Suppose that a polynomial
logarithmic $L$-valued field $a(z,\zeta)$ is local with every member of
$\mathcal S$. If its exponents lie in one $\Z$-coset and are bounded above or bounded below in the integer indexing of that coset, then $a=0$.
\end{lemma}
\begin{proof}
Fix a coefficient of a power of $\zeta$ and write its $z$-expansion as
$\sum_i a_i z^{-i-\alpha_0}$, with $i$ in a subset of $\Z$. The locality
relation with $b(w,\eta)\in\mathcal S$ holds coefficientwise in the
independent variables $\zeta,\eta$. If the indices have a largest value
$i_0$, the coefficient of $z^{-i_0-\alpha_0}$ in
$(z-w)^N[a(z),b(w,\eta)]$ is
$(-w)^N[a_{i_0},b(w,\eta)]$. If they have a least value $i_0$, the
coefficient of $z^{-i_0-\alpha_0+N}$ is
$[a_{i_0},b(w,\eta)]$. Thus the extreme coefficient commutes with every
coefficient of $b$. Removing it and proceeding through the indices shows
that every coefficient of $a$ commutes with every coefficient of $b$.
Repeating this for $b\in\mathcal S$ and for each logarithmic coefficient
proves that all coefficients of $a$ are central in $L$.
\end{proof}

\begin{lemma}\label{lem:finitegen}
Let $C$ be a Lie conformal superalgebra with a direct grading
\[
 C=\bigoplus_{\Delta\in\C}C_\Delta,
 \qquad \dd C_\Delta\subseteq C_{\Delta+1}.
\]
Assume that the support is contained in finitely many $\Z$-cosets,
the dimensions $\dim C_\Delta$ are uniformly bounded, and $\dd$ is
injective. If some $\ell\in C$ satisfies $\ell_{(0)}=\dd$, then $C$ is a
finite free $\C[\dd]$-module. In each weight coset, its module rank is
the eventual value of $\dim C_{\Delta+n}$ as $n\to+\infty$.
\end{lemma}
\begin{proof}
We begin with
\begin{equation}\label{eq:translation-separated}
 \bigcap_{r\geq0}\dd^r C=0.
\end{equation}
If $0\ne x=\dd^r y_r$ for every $r$, then sesquilinearity gives
\[
 [\ell_\lambda x]=(\dd+\lambda)^r[\ell_\lambda y_r].
\]
The polynomial $[\ell_\lambda y_r]$ is nonzero, since its constant
coefficient is $\dd y_r\ne0$. Multiplication by $(\dd+\lambda)^r$
raises its $\lambda$-degree by exactly $r$, since the leading
coefficient is unchanged.
This contradicts the finite $\lambda$-degree of $[\ell_\lambda x]$.

Fix a weight coset with nonzero support. If its support were unbounded
below, then, for any fixed weight $\Delta$ in that coset, every space
$C_{\Delta-r}$ would be nonzero. Indeed, a nonzero vector at a still
lower weight can be translated injectively to $\Delta-r$. The spaces
\[
 0\ne\dd^rC_{\Delta-r}\subseteq C_\Delta\qquad(r\geq0)
\]
form a descending chain. Finite dimensionality forces this chain to
stabilize at a nonzero space, contrary to
\eqref{eq:translation-separated}. Thus each coset is bounded below.

In each nonempty weight coset, choose $\Delta_0$ in the support
such that $C_{\Delta_0+n}=0$ for $n<0$. Choose complements $V_n$
respecting parity, so that
\[
 C_{\Delta_0+n}=\dd C_{\Delta_0+n-1}\oplus V_n
 \qquad(n\geq0).
\]
Injectivity of $\dd$ gives, by induction,
\[
 C_{\Delta_0+n}=\bigoplus_{j=0}^n\dd^{\,n-j}V_j.
\]
Since
\[
 \dim V_n=\dim C_{\Delta_0+n}-\dim C_{\Delta_0+n-1},
\]
the bounded nondecreasing sequence of homogeneous dimensions makes
all but finitely many $V_n$ zero. Bases of the nonzero $V_n$ form a
free $\C[\dd]$-basis for this coset. Their total number is its
eventual homogeneous dimension. There are only finitely many cosets.
\end{proof}

\begin{lemma}\label{lem:translation-tails}
For each physical-mode coset $\beta$, the map
\[
 T:L[k]\longrightarrow L[k-1]
\]
is an isomorphism for all but finitely many $k\in\beta+\Z$.
\end{lemma}
\begin{proof}
On the cover, $T$ is represented by $B_--kI$, which is invertible
outside the finite spectrum of $B_-$. For $k\ne0$, write
$K_k=\ker(\pi:t^k\otimes U^\beta\to L[k])$.
The anchored-family argument in Proposition~\ref{prop:punctured-descent}
shows that $\dim K_k$ is constant outside a finite subset of the
lattice: compare every generic fiber with one fixed nonzero fiber
using the forward and reverse invertible jet matrices. Thus, outside
a finite set, the inclusion
\[
 (B_--kI)K_k\subseteq K_{k-1}
\]
is an equality, since its two sides have the same dimension.
The induced map on the two quotients is an isomorphism. Both weights
are nonzero here, so the cover maps onto the full spaces $L[k]$ and
$L[k-1]$. 
\end{proof}

\subsection{The finite logarithmic model}
For each physical-mode coset $\beta$, choose $k_\beta\in\beta+\Z$ such
that $T$ is injective at $k_\beta+j$ for every $j\geq1$. Lemma~\ref{lem:translation-tails}
permits this choice independently of the conformal weight. Let $s$ be
the number of physical-mode cosets, and choose $M$ with $\dim L[k]\leq M$.

\begin{theorem}\label{thm:logfinite}
Let $L$ be a simple superconformal Lie algebra, with the residue derivation
$D$ of Theorem~\ref{thm:residue}. Then the construction above gives a
finite free Lie conformal superalgebra $C$ of $L$-valued logarithmic fields
whose coefficients span $L$. Its rank satisfies
\[
 \rank_{\C[\dd]}C\leq\sum_\beta\dim L[k_\beta]\leq sM.
\]
The realization is faithful, $\dd=D_z$, and $\delta=\dd_\zeta$ is a
nilpotent ordinary derivation with $\delta a=Da$. Moreover
$\mathcal H=\ell_{(1)}+\delta$ is diagonalizable,
$\ell_{(0)}=\dd$, and $\delta\ell=0$.
\end{theorem}
\begin{proof}
Proposition~\ref{prop:log-closure} gives the conformal identities, covariance,
and common residue.
For the dimension estimate, we use the Virasoro block decomposition.
The tensor module $V(a,\beta)$ has basis $v_k$, $k\in\beta+\Z$,
and action $L_n v_k=-(k+an)v_{k+n}$. Its simple composition factors
have block label $(a+\Z,\beta)$; the trivial module and the exceptional
factors of $V(0,0)$ and $V(1,0)$ have label $(0,0)$.
For $\alpha,\beta\in\C/\Z$, let $L[\alpha,\beta]$ be the sum of
Virasoro submodules of $L$ whose simple composition factors have
label $(\alpha,\beta)$. By \cite[Lemma~2.10 and Theorem~2.11]{Billig},
these spaces form a direct Lie-superalgebra grading.
Only finitely many blocks of $L$ occur:
the finite cover contains finitely many blocks and maps onto every
nonzero physical-weight space, while a Virasoro module concentrated
at physical weight zero belongs to the block $(0,0)$.
The cover may be taken blockwise, as observed after
\cite[Theorem~2.14]{Billig}. The operator $D$ preserves every block:
this follows from $J$ on nonzero fibers and from the opposite-weight
bracket formula defining $D$ at weight zero.

The triangular conjugacy for $B_-$ uses only operators that shift
Euler eigenvalues by integers. Thus a generalized $B_-$-eigenvector
of eigenvalue $\omega$ remains in Euler block $\alpha=\omega+\Z$.
For an initial field $X_u$ its conformal weight is $1-\omega$; for
the exceptional field $Y_u$ it is $1$. In both cases the relation
between the block label and the conformal weight is
\begin{equation}\label{eq:block-weight-relation}
 \alpha=-\Delta+\Z.
\end{equation}
The Virasoro field has weight $2$ and Euler block $0$, and satisfies
the same relation. Derivatives increase $\Delta$ by one and preserve
the block label. For a nonnegative product,
\[
 \Delta(a_{(n)}b)=\Delta(a)+\Delta(b)-n-1,
\]
while its coefficients lie in the sum of the two blocks of its arguments.
Consequently \eqref{eq:block-weight-relation} holds throughout the
singular-product closure, so only finitely many conformal-weight
cosets occur.

A homogeneous field has the expansion
\begin{equation}\label{eq:homfield}
 a(z,\zeta)=e^{\zeta D}\sum_{k\in\beta+\Z}a_kz^{-k-\Delta},
 \qquad a_k\in L[-\Delta+\Z,\beta]\cap L[k].
\end{equation}
Translation covariance is exactly
\begin{equation}\label{eq:coeff-rec}
 T a_k=\bigl((1-\Delta-k)I+D\bigr)a_{k-1}.
\end{equation}
If $a_{k_\beta}=0$, this recurrence and injectivity of $T$ give
$a_{k_\beta+j}=0$ for every $j\geq0$. The common residue identity
then gives the same vanishing for every logarithmic coefficient.
The field is one-sided, so Lemma~\ref{lem:onesided} gives $a=0$.
Its locality is with the initial fields, whose coefficients contain
all nonzero physical-weight spaces; these generate $L$ by
pluperfectness. Thus evaluation gives an injection
\begin{equation}\label{eq:evaluation-inj}
 C_{\Delta,\beta}\longrightarrow
 L[-\Delta+\Z,\beta]\cap L[k_\beta],\qquad
 a\longmapsto a_{k_\beta}.
\end{equation}
In particular,
\begin{equation}\label{eq:dim-bound}
 \dim C_\Delta\leq\sum_\beta\dim L[k_\beta]\leq sM.
\end{equation}

If $D_za=0$, then $a$ is constant. On a term
$z^q p(\zeta)$ the equation is
$(q+\dd_\zeta)p=0$. For $q\ne0$ this operator is invertible on
polynomials, while for $q=0$ its kernel consists of constants.
The resulting constant coefficient is central by
Lemma~\ref{lem:onesided}, and hence zero. Thus $D_z$ is injective,
and Lemma~\ref{lem:finitegen} proves finite generation.

In fact the sharper blockwise estimate in
\eqref{eq:evaluation-inj} gives a rank bound. The translation
preserves $\beta$, so the rank contributed by a pair
$(\alpha,\beta)$ is its eventual homogeneous dimension, at most
$\dim(L[\alpha,\beta]\cap L[k_\beta])$. Summing the direct block
decomposition gives
\begin{equation}\label{eq:finite-rank-bound}
 \rank_{\C[\dd]} C
 \leq\sum_{\alpha,\beta}
       \dim\bigl(L[\alpha,\beta]\cap L[k_\beta]\bigr)
 =\sum_\beta\dim L[k_\beta]\leq sM.
\end{equation}

For the coefficient span, let $a\langle m\rangle$
denote the coefficient of $z^{-m-1}$ in $a(z,0)$, with $m$ in its
exponent coset. The formal monodromy operator is
$\phi=e^{2\pi iD_\zeta}$. On each exponent coset its scalar phase is
constant, and $\dd_\zeta$ is nilpotent on each polynomial-logarithmic
field. Thus the ambient monodromy is locally finite. On $C$ the nilpotent
factor is $e^{2\pi i\delta}$. The logarithmic commutator formula
\cite[(5.9)]{Bakalov}, with $\mathcal N=-\delta$, reads
\begin{equation}\label{eq:actual-log-modes}
 [a\langle m\rangle,b\langle n\rangle]
 =\sum_{j\geq0}
 \left(\left(\binom{m-\delta}{j}a\right)_{(j)}b\right)
       \langle m+n-j\rangle.
\end{equation}
The sum is finite: only finitely many
$\delta^i a$ occur, and each is local with $b$. Since $C$ is
$\delta$-stable, its logarithm-free coefficient span is a Lie
subsuperalgebra of $L$. The initial fields already contain every
nonzero physical-weight space in that span. Pluperfectness therefore
gives $L[0]$, and the span is all of $L$. The common residue identity
shows that the other logarithmic coefficients lie in the same span.
Since $C$ is a space of $L$-valued fields, its realization is faithful.
\end{proof}

\section{Coefficient algebras and the socle}\label{sec:coeff}
Let $C$ be the faithful finite logarithmic model of
Theorem~\ref{thm:logfinite}, and let $\phi$ denote the monodromy of its
logarithmic realization. Write its Jordan decomposition as
\[
 \phi=\sigma e^{2\pi i\delta},
\]
where $\sigma$ is diagonalizable, $\delta$ is nilpotent, and they commute.
On $C_{\Delta,\beta}$, the automorphism $\sigma$ acts by
$e^{-2\pi i(\Delta+\beta)}$. The homogeneous product rule shows that
$\sigma$ preserves all conformal products. Thus $\sigma\ell=\ell$, and
$\sigma$ commutes with $\mathcal H$ and $\delta$.

\subsection{The logarithmic coefficient algebra and its central kernel}
For a homogeneous field in $C_{\Delta,\beta}$, the physical and coefficient
indices satisfy $m=k+\Delta-1$, so its coefficient modes lie in
$\Delta+\beta+\Z$. Let
\[
 \Gamma=\bigl\langle\Z,\ \Delta+\beta:
                      C_{\Delta,\beta}\ne0\bigr\rangle_{\Z}\subseteq\C,
\]
where a representative of each $\beta\in\C/\Z$ is chosen.
The group $\Gamma$ is independent of these choices, since it contains $\Z$.
On the group algebra of formal monomials
\[
 R_\Gamma=\C[z^\gamma:\gamma\in\Gamma][\zeta],\qquad
 d_z(z^\gamma)=\gamma z^{\gamma-1},\qquad d_z\zeta=z^{-1},
\]
form the differential scalar extension $R_\Gamma\otimes C$, with
translation $\dd_R=\dd+d_z$ and products
\begin{equation}\label{eq:scalar-extension-products}
 (fa)_{(n)}(gb)=\sum_{j\geq0}\frac{d_z^jf}{j!}\,g\,
                         (a_{(n+j)}b).
\end{equation}
The scalar-extension identities give a Lie conformal superalgebra.
The operator $\nabla=\dd_\zeta+\delta$ is an ordinary derivation and
commutes with $\dd_R$. Its horizontal subspace in the sector determined
by $\sigma$ is
\[
 \mathcal A_\delta=
 \Span\{z^m e^{-\zeta\delta}a:
          \sigma a=e^{-2\pi i m}a\}.
\]
It is closed under the products and under $\dd_R$. All exponentials are
polynomial, since $\delta$ is nilpotent. Define
\[
 \Gcal(C;\sigma,\delta)=\mathcal A_\delta/\dd_R\mathcal A_\delta,
\]
with the Lie bracket induced by the zero product. For
$\sigma a=e^{-2\pi i\alpha}a$ and $m\in\alpha+\Z$, write $a[m]$ for
the class of $z^m e^{-\zeta\delta}a$. Direct differentiation gives
\[
 d_z^j\bigl(z^m e^{-\zeta\delta}a\bigr)
 =z^{m-j}e^{-\zeta\delta}
       (m-\delta)(m-1-\delta)\cdots(m-j+1-\delta)a.
\]
Consequently the mode relations are
\begin{equation}\label{eq:log-translation}
 (\dd a)[m]=-((m-\delta)a)[m-1]
\end{equation}
and
\begin{equation}\label{eq:log-bracket}
 [a[m],b[n]]=
 \sum_{j\ge0}\left(\left(\binom{m-\delta}{j}a\right)_{(j)}b\right)[m+n-j].
\end{equation}
The sum is finite. Evaluation of modes gives a surjective homomorphism
\[
 \rho:\Gcal(C;\sigma,\delta)\longrightarrow L.
\]

\begin{lemma}\label{lem:separation}
Write $Y(a,z,\zeta)$ for the field associated with $a\in C$ in the
faithful logarithmic realization, and set
$X(a,z)=Y(a,z,\zeta)|_{\zeta=0}$. If
\[
 \sum_{s\in S}z^sX(a_s,z)=0
\]
for a finite set of distinct complex numbers $S$, then every $a_s=0$.
\end{lemma}
\begin{proof}
On logarithm-free fields define
\[
 \mathcal E_0=zT-z\dd_z-D.
\]
The common residue identity gives $\mathcal E_0X(a,z)=0$, hence
\[
 \mathcal E_0(z^sX(a,z))=-s z^sX(a,z).
\]
The summands lie in distinct eigenspaces. Since $X$ is injective, the claim follows.
\end{proof}

\begin{proposition}\label{prop:central}
The kernel of $\rho$ is the center of $\Gcal(C;\sigma,\delta)$. In particular
\[
 L\simeq \Gcal(C;\sigma,\delta)/Z(\Gcal(C;\sigma,\delta)).
\]
\end{proposition}
\begin{proof}
Let $v=\sum_i a_i[m_i]\in\ker\rho$. Commuting its image with $X(b,z)$ and using~\eqref{eq:log-bracket} gives a finite relation of the form
\[
 0=\sum_{i,j}z^{m_i-j}X\!\left(\left(\binom{m_i-\delta}{j}a_i\right)_{(j)}b,z\right).
\]
After grouping terms with the same value of $m_i-j$, we apply
Lemma~\ref{lem:separation}. The resulting identities are the grouped
coefficients of $[v,b[n]]$ in the universal algebra, so $v$ is central.
Conversely, a central element maps into $Z(L)=0$.
\end{proof}

\subsection{Semisimplicity and the socle}
For a subspace $I\subseteq C$, let $\operatorname{coeff} I$ be the
span in $L$ of the logarithm-free coefficients of its fields.
If $I$ is a conformal ideal stable under $\sigma$ and $\delta$, then
\eqref{eq:log-bracket} makes $\operatorname{coeff} I$ an ideal of $L$.
The conformal derived series is $I^{(0)}=I$ and
$I^{(r+1)}=\Span_{\C[\dd]}\{a_{(j)}b:a,b\in I^{(r)},\ j\ge0\}$.
An ideal is solvable if this series reaches zero. The solvable radical
$\Rad C$ is the largest solvable ideal; a conformal algebra is
semisimple when its radical is zero. The radical is invariant under
even automorphisms.
Since $\delta$ is even and nilpotent, $e^{t\delta}$ is an automorphism
for every $t\in\C$. Differentiation at $t=0$ shows that the radical is
$\delta$-stable; it is also $\sigma$-stable. Each term of its conformal
derived series is stable under both operators. If $I^{(r)}$ denotes
such a term, \eqref{eq:log-bracket} gives
\[
 [\operatorname{coeff} I^{(r)},\operatorname{coeff} I^{(r)}]
 \subseteq\operatorname{coeff} I^{(r+1)}.
\]
Thus the coefficient span of $\Rad C$ is a solvable ideal of $L$. If $\Rad C\ne0$, faithfulness gives a nonzero such ideal, contradicting simplicity. Thus $C$ is semisimple.

The socle $\Soc C$ is the sum of the minimal nonzero conformal ideals.
Write
\[
 J=\Soc C=M_1\oplus\dots\oplus M_r.
\]
The finite semisimple structure theory in \cite[Theorem 5.1 and Lemma 5.2]{FKR} gives perfect minimal ideals and says that distinct minimal ideals commute. By \cite[Theorem 4.7]{FKR}, each $M_i$ is either a current algebra tensored with a finite exterior algebra or a simple noncurrent conformal algebra. Every ordinary derivation preserves each perfect ideal: differentiating its generating products yields brackets with one argument in that ideal. Hence $\delta$ preserves each $M_i$, whereas $\sigma$ may permute them.

For one $\sigma$-orbit $\mathcal O$, the sum $J_{\mathcal O}$ is stable under $\sigma$ and $\delta$. Faithfulness makes its coefficient span a nonzero ideal of $L$, hence all of $L$. Two different orbits would give two commuting coefficient images equal to $L$, impossible because $L$ is nonabelian. Thus $\sigma$ acts transitively on the $M_i$, and all are isomorphic.

The faithful fields of $J$ therefore have coefficient span $L$.
For $a\in J$, the ambient translation $T$ still gives
$(zT-z\dd_z-D)X(a,z)=0$. Lemma~\ref{lem:separation} and the proof of
Proposition~\ref{prop:central}, applied directly to these fields, yield
\[
 \ker\bigl(\Gcal(J;\sigma,\delta)\longrightarrow L\bigr)
 =Z\bigl(\Gcal(J;\sigma,\delta)\bigr).
\]
This argument does not require an embedding of the coefficient algebra
of $J$ into that of $C$. We retain $C$ for its Virasoro state and the
conformal derivations used in untwisting.

\subsection{Excluding current minimal ideals}\label{ft:sec:current}
\begin{proposition}\label{ft:prop:current}
Assume a finite conformal algebra has a faithful logarithmic realization in a nonabelian simple Lie superalgebra $L$, with commuting operators $\sigma,\delta$, where $\sigma$ is a diagonalizable automorphism and $\delta$ is an even nilpotent ordinary derivation. Suppose a $\sigma$-orbit sum $J=M_1\oplus\cdots\oplus M_r$ of perfect minimal ideals has coefficient span $L$, and the reconstruction map
\[
 \Gcal(J;\sigma,\delta)\longrightarrow L
\]
has central kernel. Then the ideals in this orbit cannot be of current type
$(\Cur\s)\otimes\Lambda(q)$, where $\s$ is a finite-dimensional simple
Lie superalgebra and $\Lambda(q)$ is the exterior algebra on $q\geq0$
odd generators.
\end{proposition}
\begin{proof}
The derivation $\delta$ preserves each $M_i$, and $\sigma$ permutes them. Suppose first that $q>0$, and let $\Lambda(q)_+$ be
the augmentation ideal. Then
\[
 I_i=\C[\partial]\otimes\s\otimes\Lambda(q)_+
\]
is nilpotent and is the solvable radical of $M_i$, since
$M_i/I_i\simeq\Cur\s$ is simple. Thus $e^{t\delta}$ preserves $I_i$
for every $t\in\C$, and differentiation gives
$\delta I_i\subseteq I_i$. Also $\sigma$ permutes the $I_i$.
Consequently $I=\bigoplus_i I_i$ is a nonzero solvable conformal
ideal of $J$ stable under both operators. The derived-series argument
used for $\Rad C$ shows that its coefficient span is a solvable ideal
of $L$. Faithfulness makes this ideal nonzero, contrary to simplicity.

When $q=0$, on one summand $\Cur\s$, an even ordinary derivation has the form
\begin{equation}\label{ft:eq:current-der}
 \delta=c\partial+d,
 \qquad c\in\C,\quad d\in\Der(\s)_{\overline0},
\end{equation}
as in \cite[Remark 3.12(6)]{FKR}. On the finite-dimensional constant space, write
\[
 \delta(a)=\sum_{j=0}^m\partial^jA_j(a).
\]
If $m\ge2$, comparison of the coefficient of $\partial^{m-1}\lambda$ in the derivation identity gives $[a,A_m(b)]=0$ for all $a,b\in\s$. Centerlessness gives $A_m=0$, a contradiction. Thus $m\le1$. Comparing the remaining coefficients gives $A_0\in\Der\s$ and
\[
 A_1[a,b]=[A_1a,b]=[a,A_1b].
\]
The even centroid of $\s$ is $\C$, so $A_1=c\id$. This proves~\eqref{ft:eq:current-der}.

Nilpotence forces $c=0$: the highest $\partial$-degree term of $(c\partial+d)^k(a)$ is $c^k\partial^ka$. Therefore $\delta$ preserves the constant space in every summand.

Write $J=\Cur\g$, where $\g$ is a direct sum of the simple finite-dimensional algebras in the orbit. Every ideal of $\g$ has zero center. By \cite[Theorem 3.4]{KLP}, every conformal automorphism of $\Cur\g$ preserves $\g$. In particular, $\sigma$ also preserves the constant generators.

For constant $a,b$, every state $\binom{m-\delta}{j}a$ is still constant. Consequently all terms in~\eqref{eq:log-bracket} with $j>0$ vanish, and
\[
 [a[m],b[n]]=[a,b][m+n].
\]
The relation
\[
 (\partial a)[m]=-((m-\delta)a)[m-1]
\]
eliminates derivative states without imposing relations among the
constant-generator modes. Indeed, it is triangular in the degree in
$\partial$, with leading coefficient one, and gives the normal form
\[
 (\partial^r a)[m]=(-1)^r
 \bigl((m-\delta)(m-1-\delta)\cdots(m-r+1-\delta)a\bigr)[m-r].
\]
Thus $\Gcal(J;\sigma,\delta)$ is the usual twisted loop algebra
\[
 \mathcal L(\g,\sigma)=
 \bigoplus_{\alpha\in\C/\Z}t^\alpha\C[t,t^{-1}]\otimes\g^\alpha,
 \qquad
 \g^\alpha=\{a:\sigma a=e^{-2\pi i\alpha}a\}.
\]
It is centerless: bracketing a central finite sum with homogeneous
monomials separates its distinct exponents and forces each coefficient
into $Z(\g)=0$. Choosing one representative of each eigenvalue coset makes
$\mathcal L(\g,\sigma)$ a free $A_0=\C[t,t^{-1}]$-module of rank
$\dim\g$. Its bracket is $A_0$-bilinear:
\[
 [f(t)t^\alpha\otimes x,\,g(t)t^\beta\otimes y]
 =f(t)g(t)t^{\alpha+\beta}\otimes[x,y].
\]
This description is independent of the chosen representatives of the
cosets. Hence $(t-1)\mathcal L(\g,\sigma)$ is a nonzero proper Lie
ideal, and evaluation at $t=1$ shows that the quotient has vector-space
dimension $\dim\g$.
Hence the loop algebra is not simple. A central-kernel surjection from this centerless algebra to $L$ would be an isomorphism, a contradiction.
\end{proof}

Thus
\begin{equation}\label{eq:socle-cartantype}
 J=M_1\oplus\dots\oplus M_r\simeq M^{\oplus r},
\end{equation}
where $M$ is one of the noncurrent simple conformal superalgebras in~\eqref{eq:list}, and $\sigma$ permutes the copies transitively. 

\section{Innerness of the residue and untwisting}\label{sec:untwist}
A conformal derivation of a Lie conformal superalgebra $C$ is a
polynomial family $\psi_\lambda:C\to C[\lambda]$ satisfying
\[
 \psi_\lambda(\dd a)=(\dd+\lambda)\psi_\lambda(a)
\]
and, for homogeneous $\psi,a,b$,
\[
 \psi_\lambda([a_\mu b])
 =[(\psi_\lambda a)_{\lambda+\mu}b]
  +(-1)^{p(\psi)p(a)}[a_\mu(\psi_\lambda b)].
\]
Write
\[
 \psi_\lambda=\sum_{j\geq0}\frac{\lambda^j}{j!}\psi_{(j)}.
\]
The inner conformal derivation associated with $h\in C$ is
$(\ad h)_\lambda a=[h_\lambda a]$.
We write $\Cder C$ for the conformal derivation algebra, with
$(\dd\psi)_\lambda=-\lambda\psi_\lambda$ and conformal commutator
\[
 ([\psi_\lambda\chi])_\mu
 =\psi_\lambda\chi_{\mu-\lambda}
  -(-1)^{p(\psi)p(\chi)}\chi_{\mu-\lambda}\psi_\lambda.
\]

Let $V=\ell_{(1)}$ and $\mathcal H=V+\delta$. The operator $\mathcal H$
is semisimple, $\delta$ is nilpotent, and they commute. Hence
$V_{\sm}=\mathcal H$ and $V_{\nil}=-\delta$ on the logarithmic field algebra.

\subsection{Constituents excluded by translation}\label{ft:sec:translation}
The conformal translation is denoted by $\partial$; the odd exterior derivatives are $\partial_i=\partial/\partial\xi_i$. We use the conventions of \cite[Examples 3.3--3.6]{FKR}.

Write $\Lambda_N=\Lambda(\xi_1,\ldots,\xi_N)$,
$W(N)=\Der\Lambda_N$, and $E=\sum_i\xi_i\partial_i$ for the exterior
Euler operator. Then
\[
 W_N=\C[\partial]\otimes\bigl(W(N)\oplus\Lambda_N\bigr).
\]
Its brackets on generators are
\begin{align*}
 [X_\lambda Y]&=[X,Y],\\
 [X_\lambda f]&=X(f)-\lambda(-1)^{p(X)p(f)}fX,\\
 [f_\lambda g]&=-\partial(fg)-2\lambda fg.
\end{align*}
In particular, $(-1)_{(0)}=\partial$ on $W_N$. For an even $X+f\in W(N)\oplus\Lambda_N$ and an odd vector field $Q\in W(N)$, skew-symmetry gives
\begin{equation}\label{ft:eq:zero}
 (X+f)_{(0)}Q=[X,Q]-Q(f)-\partial(fQ).
\end{equation}
The vector-field and function summands are independent over $\C[\partial]$.

For $A=\sum_iP_i(\partial,\xi)\partial_i+f(\partial,\xi)$, define
\[
 \Div A=\sum_i(-1)^{p(P_i)}\partial_iP_i-\partial f,
 \qquad \Div_a A=\Div A+af.
\]
Thus
\[
 S_{N,a}=\ker\Div_a.
\]
When $N$ is even, set $\nu=\xi_1\cdots\xi_N$, $F=1+\nu$, and $u=F^{-1}=1-\nu$. The second family is
\[
 \widetilde S_N=\{A\in W_N:\Div(FA)=0\}=uS_{N,0}.
\]

\begin{proposition}\label{ft:prop:no-translation}
Let $N\ge2$. There is no even state $s\in S_{N,a}$ with $a\ne0$ such that $s_{(0)}=\partial$ on $S_{N,a}$. For even $N\ge2$, there is no even state $s\in\widetilde S_N$ such that $s_{(0)}=\partial$ on $\widetilde S_N$.
\end{proposition}
\begin{proof}
Suppose first that $s\in S_{N,a}$ has the stated zero mode. Write
\[
 s+1=X+f+\partial b,
 \qquad X\in W(N)_{\overline0},\quad f\in(\Lambda_N)_{\overline0},\quad b\in(W_N)_{\overline0}.
\]
Since $(\partial b)_{(0)}=0$, the zero mode of $X+f$ annihilates every $\partial_i\in S_{N,a}$. Formula~\eqref{ft:eq:zero}, by its coefficient of $\partial$, gives $f\partial_i=0$, and hence $f=0$. It then gives $[X,\partial_i]=0$ for all $i$. Writing $X=\sum_jP_j\partial_j$, the coefficients $P_j$ are odd and all their exterior derivatives vanish. They are therefore zero. Consequently
\[
 s=-1+\partial b.
\]
Applying $\Div_a$, which commutes with $\partial$, gives
\[
 0=\partial-a+\partial\Div_a b.
\]
Evaluation at $\partial=0$ gives $a=0$, a contradiction.

Now suppose that $s\in\widetilde S_N$ has zero mode $\partial$. Again write $s+1=X+f+\partial b$. The odd vector fields
\[
 Q_i=u\partial_i
\]
belong to $\widetilde S_N$, because $FQ_i=\partial_i$. Formula~\eqref{ft:eq:zero} first gives $fQ_i=0$. As $u$ is a unit, $f=0$. Thus $[X,Q_i]=0$ for all $i$.

Write $X=\sum_jP_j\partial_j$ with odd $P_j\in\Lambda_N$. Direct comparison of the vector-field coefficients gives
\begin{equation}\label{ft:eq:centralizer}
 u\partial_iP_j=\delta_{ij}X(u).
\end{equation}
For $i\ne j$, this says that $P_j$ is independent of $\xi_i$. Hence $P_j=c_j\xi_j$. For $i=j$, the constant term of~\eqref{ft:eq:centralizer} gives $c_i=0$, because $X(u)=-(\sum_jc_j)\nu$ has zero constant term. Therefore $X=0$, and again $s=-1+\partial b$.

The defining condition for $\widetilde S_N$ now gives
\[
 0=\Div(Fs)=\partial F+\partial\Div(Fb).
\]
Since multiplication by $\partial$ is injective on $\C[\partial]\otimes\Lambda_N$,
\[
 \Div(Fb)=-F.
\]
Evaluate at $\partial=0$ and take the coefficient of $\nu$. The coefficient of $\nu$ in any exterior divergence is zero: an odd derivative lowers exterior degree. The coefficient on the right is $-1$. This contradiction proves the second assertion.
\end{proof}

\begin{corollary}\label{ft:cor:exclude}
Suppose $C$ is a Lie conformal superalgebra containing an even element $\ell$ with $\ell_{(0)}=\partial$. A conformal ideal $M$ of $C$ isomorphic to $S_{N,a}$ with $a\ne0$, or to $\widetilde S_N$, cannot occur.
\end{corollary}
\begin{proof}
The restriction of $\operatorname{ad}\ell$ to $M$ is an even conformal derivation. For each of the two stated families, every conformal derivation is inner by \cite[Proposition 3.11]{FKR}. Therefore there is an even $s\in M$ such that
\[
 [\ell_\lambda x]=[s_\lambda x]\qquad(x\in M).
\]
Taking the zero product gives $s_{(0)}=\partial$ on $M$, contradicting Proposition~\ref{ft:prop:no-translation}.
\end{proof}

\begin{remark}
This exclusion concerns conformal ideals in a model with
$\ell_{(0)}=\partial$. It does not exclude the corresponding
superconformal Lie algebras or identify their other realizations.
\end{remark}

\subsection{Innerness of the residue}
Let $\Der M$ denote the Lie superalgebra of ordinary derivations, and write
\[
 \Inder M=\{a_{(0)}:a\in M\},\qquad
 \mathfrak o(M)=\Der M/\Inder M.
\]
Inner ordinary derivations are zero modes; inner conformal derivations
are the families $[a_\lambda\,\cdot\,]$.
For $\nu=\xi_1\cdots\xi_N$, we use $E_{(0)}$ and $\nu_{(0)}$ also for
their outer classes. The following descriptions of the underlying vector
superspaces are stated in \cite[Remark~3.12]{FKR}, where the derivation
calculations are left to the reader:
\begin{align*}
 \mathfrak o(S_{N,a})&=\C E_{(0)} &&(a\ne0),\\
 \mathfrak o(S_{N,0})&=\C E_{(0)}\oplus\C\nu_{(0)} &&(N>2),\\
 \mathfrak o(S_{2,0})&\simeq\slie_2,\\
 \mathfrak o(\widetilde S_N)&=\C\nu_{(0)},\\
 \mathfrak o(K'_4)&=\C\nu_{(0)}.
\end{align*}
The class of $E_{(0)}$ is even, and that of $\nu_{(0)}$ has parity $N\bmod2$.
The identification with $\slie_2$ is also an isomorphism of Lie algebras.
For $W_N$, $K_N$ with $N\ne4$, and $CK_6$, every ordinary derivation is inner.

A grading operator on a conformal algebra is an even endomorphism
$\mathcal H$ satisfying
\[
 [\mathcal H,\dd]=\dd,\qquad
 [\mathcal H,a_{(j)}]
   =(\mathcal Ha)_{(j)}-(j+1)a_{(j)}.
\]

The only remaining cases for the innerness of $\delta$ are
$S_{N,0}$ and $K'_4$. The following lemma describes the grading action
on their ordinary derivations.

\begin{lemma}
\label{lem:outer-grading}
Let $M$ be a finite conformal algebra with a diagonalizable grading
operator $\mathcal H$. Then $\ad\mathcal H$ is diagonalizable on
$\Der M$ and on $\Der M/\Inder M$. If an ordinary derivation $\delta$
satisfies $\delta^s=0$ on $M$, then
$(\ad\delta)^{2s-1}=0$ on both spaces.
\end{lemma}
\begin{proof}
Choose finite homogeneous $\C[\dd]$-generators $v_i$ with
$\mathcal Hv_i=\Delta_i v_i$. Evaluation embeds
$\End_{\C[\dd]}M$ into $\bigoplus_i M$ by $d\mapsto(dv_i)_i$.
It intertwines $\ad\mathcal H$ with
$\bigoplus_i(\mathcal H-\Delta_i I)$, which is diagonalizable.
The derivations form an invariant subspace. Moreover
\[
 [\mathcal H,h_{(0)}]=(\mathcal Hh)_{(0)}-h_{(0)},
\]
so inner derivations form an invariant subspace as well. Restriction
and quotient preserve diagonalizability. Finally, the binomial
expansion of $(\ad\delta)^{2s-1}(d)$ has a factor $\delta^s$ on one
side of $d$ in every term, proving nilpotence.
\end{proof}

\begin{proposition}\label{prop:inner-S0-K4}
The restriction of $\delta$ to a constituent $S_{N,0}$ or $K'_4$ is inner.
\end{proposition}
\begin{proof}
First take $M=S_{N,0}$. By \cite[Proposition~3.11]{FKR}, every
conformal derivation of $M$ is inner. The conformal action of the
ambient Virasoro field on $M$ therefore equals the action of an even
state $s\in M$. Its zero mode is $\dd$. Comparing $s$ with
\[
 s_*=-1+\frac1N\dd E
\]
inside $W_N$ gives
\[
 s=s_*+\dd b,
\]
with $b\in S_{N,0}$: modulo $\dd W_N$, apply the zero mode of $s-s_*$ to all constant odd derivations $\dd_i$; the coefficient of $\dd$ first eliminates the function part, and the remaining vector field has constant coefficients, whose even part is zero. Divergence then places the primitive in $S_{N,0}$.

Hence $V|_M=(s_*)_{(1)}-b_{(0)}$.
We show that $\ad V$ vanishes on $\mathfrak o(M)$.
Set $V_*=(s_*)_{(1)}$. In the ambient $W_N$,
\[
 [s_*{}_{\lambda}E]=(\dd+\lambda)E,
 \qquad [s_*{}_{\lambda}\nu]=(\dd+\lambda)\nu.
\]
For the second equality use $[E_\lambda\nu]=N\nu$ and $\nu E=0$.
The mode commutator formula gives
$[V_*,E_{(0)}]=[V_*,\nu_{(0)}]=0$.
For $N>2$ these classes span the outer algebra. For $N=2$ they are,
after rescaling, a Cartan element and one root vector of the outer
$\slie_2$. An even derivation killing both kills the opposite root
vector as well: the Cartan relation makes its image a scalar multiple
of that vector, and the bracket with the first root vector makes the
scalar zero. The inner term $b_{(0)}$ also acts trivially on the quotient. Thus $\ad V$ is zero on $\mathfrak o(M)$. Since
\[
 \ad V=\ad\mathcal H-\ad\delta
\]
is the difference of commuting semisimple and nilpotent operators by Lemma~\ref{lem:outer-grading}, both induced operators vanish. The class of $\delta$ is therefore central in $\mathfrak o(M)$. This center is zero: for $N>2$, $[E_{(0)},\nu_{(0)}]=N\nu_{(0)}$, and for $N=2$ the outer algebra is $\slie_2$. Hence $\delta|_M$ is inner.

Now take $M=K'_4$, and write $\nu=\xi_1\xi_2\xi_3\xi_4$.
The natural bracket action identifies
\[
 K_4\simeq\Cder K'_4.
\]
Indeed, in Billig's contact convention $P=I-\frac12E$, the bracket
of an exterior monomial $a$ of degree $r$ with $1$ is
\[
 [a_\lambda 1]
 =\left(1-\frac r2\right)\partial a
  +\left(2-\frac r2\right)\lambda a.
\]
The two coefficients do not vanish simultaneously. Conformal
sesquilinearity and the free exterior basis therefore make the
bracket action faithful on $K_4$. Its image contains the inner
conformal derivations and the derivation induced by $\nu$, which
generate $\Cder K'_4$ over $\C[\partial]$ by
\cite[Proposition~3.11]{FKR}. Moreover,
\[
 \partial(\operatorname{ad}\nu)
   =\operatorname{ad}(\partial\nu)
\]
is inner, since $\partial\nu\in K'_4$. Thus the outer conformal
quotient is one-dimensional and is annihilated by $\partial$.
Under this identification, the conformal action of $\ell$ on $K'_4$
is the action of an even state $s\in K_4$. Thus
$s_{(0)}=\dd$ and $V|_{K'_4}=s_{(1)}$.
We claim that
\[
 s=1+\dd b\qquad(b\in K_4).
\]
Indeed, if $a=s-1$ modulo $\dd K_4$, then $a_{(0)}1=0$ gives $P(a)=0$,
so $a$ has exterior degree $2$; applying $a_{(0)}$ to the four
generators $\xi_i$ gives $\partial_i a=0$, hence $a=0$.
It follows that $V|_{K'_4}=V_*-b_{(0)}$, where $V_*=1_{(1)}$.
The top monomial $\nu$ has conformal weight $0$ under $V_*$, so
\[
 [V_*,\nu_{(0)}]=-\nu_{(0)}.
\]
The zero mode $b_{(0)}$ may itself have an outer component, but its
adjoint action on the ordinary outer quotient $\mathfrak o(K'_4)$
is zero. Since $[V,\delta]=0$, a nonzero outer component of $\delta$
would have eigenvalue $-1$ under $\ad V$, a contradiction. Thus
$\delta|_{K'_4}$ is inner.
\end{proof}

\begin{theorem}\label{thm:compatible-inner}
For the finite semisimple model constructed above, there is an even
$h\in\Soc C$ such that
\begin{equation}\label{eq:innerdelta}
 \delta=h_{(0)}
\end{equation}
on the same ambient algebra $C$, and
\begin{equation}\label{eq:hchoice}
 \mathcal Hh=h,\qquad \sigma h=h.
\end{equation}
\end{theorem}
\begin{proof}
The current constituents have been excluded. Corollary~\ref{ft:cor:exclude}
excludes $S_{N,a}$ with $a\ne0$ and $\widetilde S_N$ in this model.
For all other noncurrent types, the ordinary derivation list and
Proposition~\ref{prop:inner-S0-K4} show that
$\delta|_{M_i}=(h_i)_{(0)}$ for an even $h_i\in M_i$.
Let $h=\sum_i h_i$ and $\varepsilon=\delta-h_{(0)}$.
For $c\in C$, $m\in\Soc C$,
\[
 0=\varepsilon[c_\lambda m]=[(\varepsilon c)_\lambda m].
\]
Thus $\varepsilon c$ centralizes the socle, whose centralizer is zero
by the proof of \cite[Theorem 5.1]{FKR}. This proves
\eqref{eq:innerdelta} on $C$.

The commuting diagonalizable operators $\mathcal H,\sigma$ preserve
the socle. Decompose $h$ into its finitely many
joint eigencomponents. The zero-mode map satisfies
\[
 [\mathcal H,h_{(0)}]=(\mathcal Hh-h)_{(0)},
 \qquad \sigma h_{(0)}\sigma^{-1}=(\sigma h)_{(0)}.
\]
Since $\delta$ commutes with both operators, only the component with
$\mathcal H$-weight $1$ and $\sigma$-eigenvalue $1$ can contribute to
$\delta$. This component gives \eqref{eq:hchoice}.
\end{proof}

\subsection{Pseudo-derivations and their filtration}

Let $C$ be a finite semisimple Lie conformal superalgebra whose
minimal ideals are noncurrent, and let $\mathcal H$ be a
diagonalizable grading operator on $C$. Write $\wt(a)=\Delta$
when $\mathcal Ha=\Delta a$. Thus $[\mathcal H,\dd]=\dd$ and, for
homogeneous states,
\[
 \wt(a_{(j)}b)=\wt(a)+\wt(b)-j-1.
\]
Let $\sigma$ be a diagonalizable automorphism commuting with $\mathcal H$.
Suppose that $h\in C$ is even, $\mathcal Hh=h$, $\sigma h=h$, and
$\delta=h_{(0)}$ satisfies $\delta^n=0$ for some $n\geq1$.
The construction does not require $\delta h=0$.

Let $J=\Soc C$. The conformal derivation algebra $M=\Cder C$ is finite and free over
$\C[\dd]$. The descriptions of the conformal derivations of the remaining minimal
ideals make $\Cder J$ finite over $\C[\dd]$, and the normalizer theorem
identifies $\Cder C$ with a $\C[\dd]$-submodule of $\Cder J$
\cite[Propositions~3.11, 4.6, and 5.6]{FKR}. Hence $M$ is finite; it is
torsion-free as a submodule of the conformal endomorphisms of the free
module $C$, and is therefore free. The grading induced on $M$ is
characterized by
\[
 [\mathcal H,\psi_{(j)}]=(w-j-1)\psi_{(j)}
 \qquad(\psi\in M\text{ of weight }w).
\]
The automorphism $\sigma$ acts on $M$ by conjugation:
\[
 (\sigma\psi)_\lambda a
   =\sigma\bigl(\psi_\lambda(\sigma^{-1}a)\bigr).
\]
Both $C$ and $M$ have finite parity-homogeneous $\C[\dd]$-bases of weight
vectors, which may also
be chosen to consist of $\sigma$-eigenvectors. Indeed, finite homogeneous
generating sets follow by taking the weight components of finite
generating sets. Their weights give a lower bound for the real parts of
the spectrum. Lift a joint eigenbasis of the finite-dimensional quotient
by $\dd$; the lifts generate by downward induction on weight, and their
number is the module rank. For $M$, the induced actions are locally finite
and diagonalizable because every conformal derivation is a finite
polynomial matrix in $\lambda$ and $\dd$ in a joint homogeneous basis of
$C$. Homogeneous components remain conformal derivations.

Let
\[
 R=\C[z,z^{-1},\zeta],\qquad
 d_z=\dd_z+z^{-1}\dd_\zeta,\qquad C_R=R\otimes C,
 \qquad \dd_R=d_z+\dd.
\]
In this subsection scalar coefficients in $R$ have weight zero: every
weight bound refers to the state in $C$. For $\psi\in M$ and $f\in R$ set
\begin{equation}\label{eq:pseudoder}
 \mathsf X(\psi,f)=\sum_{j\geq0}\frac{d_z^jf}{j!}\psi_{(j)}
 \quad\text{on }C_R.
\end{equation}
This is an $R$-linear operator, and its value on every state is a finite
sum. The scalar-extension construction, or the zero-mode calculation in
$M\ltimes C$, proves that it is an ordinary derivation of the extended
conformal algebra. In particular,
\begin{equation}\label{eq:partial-pseudo}
 [\mathsf X(\psi,f),\dd_R]=0,\qquad
 \mathsf X(\dd\psi,f)=-\mathsf X(\psi,d_zf).
\end{equation}
This agrees with the pseudo-derivation construction in
\cite[Proposition 3.11]{Li}. Its commutator identity is
\begin{equation}\label{eq:pseudo-bracket}
 [\mathsf X(\psi,f),\mathsf X(\chi,g)]
 =\sum_{j\geq0}\mathsf X\left(\psi_{(j)}\chi,
                      \frac{d_z^jf}{j!}g\right).
\end{equation}
This follows from the mode commutator formula and the Leibniz rule for $d_z$. The sum on the right is finite by
conformal locality in $M$.

Choose a joint homogeneous basis $\psi_1,\ldots,\psi_s$ of $M$, with
weights $w_1,\ldots,w_s$, and let $\Pcal$ consist of the finite sums of
operators in \eqref{eq:pseudoder}.

\begin{lemma}
\label{lem:pseudo-filtration}
Every element of $\Pcal$ has a unique expression
\begin{equation}\label{eq:pseudo-normal-form}
 Y=\sum_{i=1}^s\mathsf X(\psi_i,f_i),\qquad f_i\in R.
\end{equation}
For an integer $r\geq0$, the following conditions are equivalent:
\begin{enumerate}
\item $Y$ lowers the real part of state weight by at least $r$;
\item $f_i=0$ whenever $\operatorname{Re}w_i>1-r$.
\end{enumerate}
Write $\Pcal^r$ for this subspace. Then
\begin{equation}\label{eq:Pfiltration}
 [\Pcal^r,\Pcal^s]\subseteq\Pcal^{r+s},
 \qquad [\dd_\zeta,\Pcal^r]\subseteq\Pcal^r,
\end{equation}
and, for any integer $R_0\geq1$ with
$R_0>1-\min_i\operatorname{Re}w_i$,
\begin{equation}\label{eq:Pterminate}
 \Pcal^{R_0}=0.
\end{equation}
Moreover, every $B\in\Pcal^r$ has a primitive $K\in\Pcal^r$ satisfying
$[\dd_\zeta,K]=-B$. If $B$ is even and $\sigma$-invariant, $K$ may be
chosen with the same properties.
\end{lemma}
\begin{proof}
Existence of \eqref{eq:pseudo-normal-form} follows by moving powers of
$\dd$ from $\psi$ to $f$ with \eqref{eq:partial-pseudo}. The zero-mode
map
\[
 M/\dd M\longrightarrow\Der C,\qquad\psi+\dd M\longmapsto\psi_{(0)},
\]
is injective. If $\psi_{(0)}=0$, then $\psi_\lambda$ is divisible by
$\lambda$; the map $-\lambda^{-1}\psi_\lambda$ is again a conformal
derivation, so $\psi\in\dd M$. Consequently the operators
$(\psi_i)_{(0)}$ are linearly independent over $\C$.

In a relation of the form \eqref{eq:pseudo-normal-form}, take the maximal
real part of a weight $w_i$ for which $f_i\ne0$. The part of maximal
state-weight shift is
\[
 \sum_{\operatorname{Re}w_i\text{ maximal}}
 f_i(\psi_i)_{(0)}.
\]
Higher modes have strictly smaller real weight shift. Terms with
different imaginary parts of $w_i$ have different shifts and cannot
cancel. Within each shift, expanding the finitely many $f_i$ in Laurent
monomials in $z$ and monomials in $\zeta$ reduces cancellation to
$\C$-linear dependence of the zero modes. Thus the displayed part is
nonzero. This proves uniqueness and the stated equivalence with the
operator weight bound.

The commutator of operators lowering weight by $r$ and $s$ lowers it by
$r+s$. It belongs to $\Pcal$ by \eqref{eq:pseudo-bracket}, proving the
first assertion of \eqref{eq:Pfiltration}. The finite set of basis
weights then proves \eqref{eq:Pterminate}. Since
$[\dd_\zeta,d_z]=0$, we also have
\[
 [\dd_\zeta,\mathsf X(\psi,f)]
   =\mathsf X(\psi,\dd_\zeta f).
\]
In the unique normal form for $B$, replace each coefficient $f_i$ by
$-\int_0^\zeta f_i(z,u)\,du$. This gives the required primitive. Joint
eigenvectors in the chosen basis make this operation commute with
$\sigma$; it preserves parity and the set of basis weights.
\end{proof}

\subsection{Finite untwisting}
\begin{theorem}
\label{thm:untwist}
For $C,\mathcal H,\sigma,h,\delta$ as in the preceding subsection, there is a
$\sigma$-equivariant $R$-linear conformal automorphism $F$ of $C_R$ such
that
\begin{equation}\label{eq:Fidentities}
 F^{-1}(\dd_\zeta+\delta)F=\dd_\zeta,
 \qquad F\dd_R=\dd_RF.
\end{equation}
Both $F$ and $F^{-1}$ are finite sums on every state.
\end{theorem}
\begin{proof}
Set
\begin{equation}\label{eq:Xinitial}
 \mathcal X=\mathsf X(\ad h,\zeta)
 =\zeta\delta+Q,\qquad
 Q=\sum_{j\geq1}\frac{(-1)^{j-1}}jz^{-j}h_{(j)}.
\end{equation}
Here $\mathcal X\in\Pcal^0$ and $\delta=\mathsf X(\ad h,1)\in\Pcal^0$.
The operator $Q$ lowers state weight by at least one. The decomposition
is taken in $\End_R C_R$; no separate membership of $Q$ or
$\zeta\delta$ in $\Pcal$ is required.

The operator $\mathcal X$ is locally nilpotent. Indeed, the spectrum of
$C$ has a lower bound on its real part. In a word in $\delta$ and $Q$
applied to a fixed state, only a bounded number of occurrences of $Q$
can survive. Between them there can be fewer than $n$ consecutive
occurrences of $\delta$, since $\delta^n=0$. Thus sufficiently long
words vanish, and $F_0=e^{-\mathcal X}$ is a well-defined conformal
automorphism commuting with $\dd_R$ and $\sigma$.

On $\End_R C_R$, nilpotence of $\delta$ gives
\[
 (\ad\delta)^{2n-1}=0.
\]
Write $q=2n-1$. For $Y\in\Pcal^r$, expand
$(\ad\mathcal X)^qY$ in the full endomorphism algebra using
$\ad\mathcal X=\zeta\ad\delta+\ad Q$. The term with no $\ad Q$ is
$\zeta^q(\ad\delta)^qY=0$. Every remaining word lowers state weight
by at least $r+1$. The complete sum belongs to $\Pcal$ by
\eqref{eq:pseudo-bracket}, so Lemma~\ref{lem:pseudo-filtration} gives
\begin{equation}\label{eq:adjoint-filtration-drop}
 (\ad\mathcal X)^q\Pcal^r\subseteq\Pcal^{r+1}.
\end{equation}
In particular $(\ad\mathcal X)^{qR_0}\Pcal^0=0$.

Since $[\dd_\zeta,\mathcal X]=\delta$, finite conjugation now gives
\begin{align*}
 F_0^{-1}(\dd_\zeta+\delta)F_0
 &=\dd_\zeta+B_1,\\
 B_1&=\sum_{k\geq1}\frac{k}{(k+1)!}
                   (\ad\mathcal X)^k\delta.
\end{align*}
By \eqref{eq:adjoint-filtration-drop}, its terms vanish for
$k\geq qR_0$. Also
$[\mathcal X,\delta]=[Q,\delta]$ lowers state weight by at least one.
Hence $[\mathcal X,\delta]\in\Pcal^1$, and
$[\Pcal^0,\Pcal^1]\subseteq\Pcal^1$ implies $B_1\in\Pcal^1$.

Suppose the connection at some stage is $\dd_\zeta+B_r$, with
$B_r\in\Pcal^r$, where $r\geq1$. Choose an even, $\sigma$-invariant
$K_r\in\Pcal^r$ with $[\dd_\zeta,K_r]=-B_r$.
It strictly lowers state weight, so $e^{K_r}$ and $e^{-K_r}$ are locally
finite. The commutator filtration makes the conjugation formula finite
and gives
\[
 e^{-K_r}(\dd_\zeta+B_r)e^{K_r}
   =\dd_\zeta+B_{2r},\qquad B_{2r}\in\Pcal^{2r}.
\]
The terms linear in $B_r$ cancel; every remaining term contains a
commutator of $K_r$ with $B_r$. Iterating with $r=1,2,4,\ldots$ stops
when $r\geq R_0$. The finite product
$F=F_0e^{K_1}e^{K_2}e^{K_4}\cdots$ has the required properties.
\end{proof}

\subsection{The ordinary finite model}
For a $\C/\Z$-grading $L=\bigoplus_\alpha L^\alpha$, compatibility
means $[L^\alpha,L^\beta]\subseteq L^{\alpha+\beta}$.
A distribution of degree $\alpha$ has the form
$\sum_{m\in\alpha+\Z}a_mz^{-m-1}$ with $a_m\in L^\alpha$.
An \emph{ordinary regular realization} of $L$ consists of such a
grading, an even grading-preserving derivation $T$, and a finite
spanning family of homogeneous mutually local distributions satisfying
\begin{equation}\label{eq:ordinary-cov}
 T a(z)=\dd_z a(z).
\end{equation}
Under these hypotheses, \cite[Theorems 5.3 and 5.5]{Billig} give a finite
conformal algebra and a central extension from its twisted coefficient
algebra onto $L$.

\begin{theorem}\label{thm:ordinary-model}
Let $L$ be a simple superconformal Lie algebra. Then $L$ admits an
ordinary regular realization by a finite
semisimple Lie conformal superalgebra $C$, and
\[
 L\simeq\Alg(M^{\oplus r},\sigma)/Z\bigl(\Alg(M^{\oplus r},\sigma)\bigr),
\]
where $M$ is a noncurrent finite simple conformal superalgebra,
$\sigma$ is diagonalizable and permutes the $r$ constituents transitively.
The possibilities for $M$ are
\[
 W_N,\quad S_{N,0},\quad K_N\ (N\ne4),\quad K'_4,\quad CK_6.
\]
The construction may change the translation, the compatible
$\C/\Z$-grading, and the chosen Virasoro embedding.
\end{theorem}
\begin{proof}
The socle description \eqref{eq:socle-cartantype} and
Corollary~\ref{ft:cor:exclude} give $J\simeq M^{\oplus r}$ with the
stated possibilities for $M$, and $\sigma$ acts transitively on its
summands. Apply Theorem~\ref{thm:untwist} to $C$ and the state $h$
supplied by Theorem~\ref{thm:compatible-inner}.

For $\alpha\in\C/\Z$, let $C^\alpha$ be the $\sigma$-eigenspace with
eigenvalue $e^{-2\pi i\alpha}$. Choose a representative for each coset and
set
\[
 \Acal_0=\bigoplus_\alpha z^\alpha\C[z,z^{-1}]\otimes C^\alpha,
 \qquad \mathcal B=\Acal_0[\zeta],
 \qquad \Acal_\delta=e^{-\zeta\delta}\Acal_0.
\]
The choice of representatives does not change these spaces. Extend the
differential scalar construction to these sectors by
$d_z(z^\alpha)=\alpha z^{\alpha-1}$, with products given by
\eqref{eq:scalar-extension-products}. The eigenvalue condition implies
that $\mathcal B$ is closed under these products.
Moreover,
\[
 \Acal_0=\ker_{\mathcal B}\dd_\zeta,
 \qquad
 \Acal_\delta=\ker_{\mathcal B}(\dd_\zeta+\delta).
\]
The second identity follows by solving the polynomial differential
equation in $\zeta$, using nilpotence of $\delta$. Both kernels are
conformal subalgebras stable under $\dd_R$, because the corresponding
operators are ordinary derivations commuting with $\dd_R$.

The automorphism $F$ preserves $\mathcal B$: it commutes with $\sigma$
and has only integral Laurent powers of $z$. The same holds for its
inverse. Equation~\eqref{eq:Fidentities} therefore gives a conformal
isomorphism $F:\Acal_0\longrightarrow\Acal_\delta$ commuting with
$\dd_R$. On the horizontal generators one has
\[
 \dd_R(e^{-\zeta\delta}z^ma)
 =e^{-\zeta\delta}\bigl(z^m\dd a+z^{m-1}(m-\delta)a\bigr).
\]
Formula~\eqref{eq:scalar-extension-products} similarly gives the
binomial operator $\binom{m-\delta}{j}$ in the zero product. Thus the
two coefficient algebras are
\[
 \Alg(C,\sigma)=\Acal_0/\dd_R\Acal_0,
 \qquad
 \Gcal(C;\sigma,\delta)=\Acal_\delta/\dd_R\Acal_\delta,
\]
and $F$ induces an isomorphism of the full Lie superalgebras
\begin{equation}\label{eq:coeffiso}
 \overline F:\Alg(C,\sigma)\xrightarrow{\ \sim\ }
                       \Gcal(C;\sigma,\delta).
\end{equation}
In particular it identifies their full centers.

Write $G=e^{\zeta\delta}F$. The first identity in
\eqref{eq:Fidentities} gives $\dd_\zeta G=0$, and the same calculation
applies to $G^{-1}$. Hence $G$ and $G^{-1}$ are Laurent-polynomial on
each state. If
\[
 G(z)a=\sum_jz^jb_j,
\]
then the explicit coefficient isomorphism is
\[
 \overline F(a[m])=\sum_jb_j[m+j].
\]
Each $b_j$ belongs to the same $\sigma$-eigenspace as $a$, so every mode
on the right is allowed. Although $F$ is a conformal automorphism of the
differential scalar extension, $G$ need not be one; the factor
$e^{\zeta\delta}$ need not commute with $\dd_R$.

Let $Y$ denote the logarithmic realization and let
$X(a,z)=Y(a,z,\zeta)|_{\zeta=0}$.
Extend both maps by scalar multiplication:
\[
 Y_R(fa)=fY(a),\qquad X_R(fa)=fX(a).
\]
The products of scalar multiples of local fields satisfy
\eqref{eq:scalar-extension-products}; their Taylor expansion uses the
total derivative $D_z$. Thus $Y_R$ respects the conformal products and
translation. Define
\begin{equation}\label{eq:ordinary-transformed-fields}
 Y'(a,z)=Y_R(Fa)=X_R(G(z)a,z).
\end{equation}
The second equality follows from $Y(a)=e^{\zeta D}X(a)$ and
$Y(\delta a)=D Y(a)$, and shows that $Y'$ has no logarithms.
The composition $Y_R\circ F$ respects conformal products and
translation. Hence $Y'$ respects the ordinary conformal products and
satisfies $Y'(\dd a)=\dd_zY'(a)$.

The map $Y'$ is faithful. If $Y'(a)=0$, write
$G(z)a=\sum_jz^jb_j$ with distinct integral powers. Lemma~\ref{lem:separation}
applied to \eqref{eq:ordinary-transformed-fields} gives $b_j=0$ for all
$j$. Invertibility of $G$ then gives $a=0$. Its coefficients still span
$L$: one inclusion follows from the finite Laurent expression for $G$,
and the reverse inclusion follows by applying the expression for
$G^{-1}$ to $Y'$. The original logarithmic fields and their restrictions
$X$ have the same coefficient span because $\delta$ acts on $C$ and
$Y(a)=e^{\zeta D}X(a)$.

The realization also has the required direct grading.
The algebra $\Alg(C,\sigma)$ is graded by the mode coset in $\C/\Z$.
Its center is graded: the components of a central element commute
separately with every homogeneous element. Proposition~\ref{prop:central}
and \eqref{eq:coeffiso} show that
$\rho'=\rho\overline F$ has kernel equal to this full center. Transport
the resulting quotient grading to $L$. The original Virasoro state
$\ell$ is $\sigma$-fixed and satisfies $\ell_{(0)}=\dd$. Hence, with
\[
 t'=\rho'(\ell[0]),\qquad T'=\ad t',
\]
the ordinary mode commutator formula gives
\[
 [t',Y'(a,z)]=Y'(\dd a,z)=\dd_zY'(a,z).
\]
The element $t'$ has degree zero, so $T'$ preserves the quotient
grading. A finite $\C[\dd]$-generating set of $C$ now gives a finite
spanning family of mutually local $T'$-covariant distributions. Thus
$L$ is regular in the sense of \cite[Definition 5.1]{Billig}. This
construction permits the translation and the grading to differ from
the original physical-mode presentation.

Every conformal derivation of $C$ preserves its perfect minimal
ideals \cite[Lemma~1.5(3)]{FKR}. Hence every coefficient of such a
derivation preserves $J=\Soc C$, as do the pseudo-derivations and
finite exponentials defining $F$ and $F^{-1}$. The coefficient
comparison therefore restricts to an isomorphism
\[
 \overline F_J:\Alg(J,\sigma)\xrightarrow{\ \sim\ }
                  \Gcal(J;\sigma,\delta).
\]
Central recovery for the faithful fields of $J$ was established in
Section~\ref{sec:coeff}. Composing that recovery map with
$\overline F_J$, which identifies the full centers, gives
\begin{equation}\label{eq:ordinaryquotient}
 L\simeq\Alg(J,\sigma)/Z(\Alg(J,\sigma)),
 \qquad J\simeq M^{\oplus r}.
\end{equation}
The ambient Virasoro state gives ordinary regularity; the
central-quotient identification uses the faithful coefficient
realization of $J$.
\end{proof}

Sections~\ref{sec:endpoints}--\ref{sec:chengkac} identify these quotients
as underlying Lie superalgebras, completing the proof of
Theorem~\ref{thm:synthesis}.

\section{Twisted coefficient algebras}\label{sec:cyclic}
\subsection{Cyclic reduction}
The coefficient isomorphisms in this and the following family sections
concern the full parity-graded Lie superalgebras. They therefore identify
the full centers and induce isomorphisms of the central quotients.
They may change the distinguished Virasoro embedding. Diagonalizable
twists may have infinite order.

Let $C$ be a Lie conformal superalgebra with translation $\partial$,
and let $\mu\in\Aut C$ be diagonalizable. The mode $a[m]$ is allowed
when $\mu a=e^{-2\pi i m}a$. The ordinary coefficient relations are
\[
 (\partial a)[m]=-m a[m-1],\qquad
 [a[m],b[n]]=\sum_{j\ge0}\binom mj(a_{(j)}b)[m+n-j].
\]
This is the usual differential affinization followed by translation coinvariants; see \cite[Section 5]{Billig} and \cite[Sections 1--2]{KLP}.

Fix an integer $r\ge1$ and define $\sigma$ on $C^{\oplus r}$ by
\[
 \sigma(a_0,\ldots,a_{r-1})=(\mu a_{r-1},a_0,\ldots,a_{r-2}).
\]
The restriction of $\sigma^r$ to each summand is $\mu$, which we call
the return automorphism.
For $\mu a=e^{-2\pi i m}a$, set $\theta_m=e^{-2\pi i m/r}$ and
\[
 E_m(a)=(a,\theta_m^{-1}a,\ldots,\theta_m^{-(r-1)}a).
\]
Then $\sigma E_m(a)=\theta_m E_m(a)$ and
\[
 E_m(a)_{(j)}E_n(b)=E_{m+n}(a_{(j)}b).
\]
Consequently the Fourier identification $E_m(a)[m/r]\mapsto a\nn{m}$ gives
\begin{equation}\label{cy:eq:fourier}
 \Alg(C^{\oplus r},\sigma)\simeq\Alg_r(C,\mu),
\end{equation}
where the right-hand side has relations
\begin{align}
 (\partial a)\nn{m}&=-\frac mr a\nn{m-r},\label{cy:eq:rder}\\
 [a\nn{m},b\nn{n}]&=\sum_{j\ge0}\binom{m/r}{j}(a_{(j)}b)\nn{m+n-rj}.\label{cy:eq:rbr}
\end{align}
Every $\sigma$-eigenvector occurs in the Fourier identification, and it intertwines the derivative relations. Thus \eqref{cy:eq:fourier} is an isomorphism of the full coefficient algebras. In the scalar coordinate $t=s^r$, \eqref{cy:eq:rder}--\eqref{cy:eq:rbr} correspond to
\[
 d_r=\frac{1}{r s^{r-1}}\frac{d}{ds},\qquad
 d_r^j(s^m)=j!\binom{m/r}{j}s^{m-rj}.
\]
The normalization below combines the index change with a conformal-weight correction.

\subsection{Normalization in a primary quadratic presentation}
\begin{definition}
A primary quadratic presentation consists of a finite-dimensional superspace $A$, a semisimple even operator $P$ on $A$, and a Lie conformal superalgebra free as a module
\[
 C=\C[\partial]\otimes A,
\]
with an even Virasoro state $\ell\in A$ satisfying $P\ell=\ell$.
Writing $A_p=\ker(P-pI)$, we require
\begin{equation}\label{cy:eq:primary}
 [\ell_\lambda a]=(\partial+(p+1)\lambda)a\qquad(a\in A_p).
\end{equation}
The bracket of generators is required to have the form
\begin{equation}\label{cy:eq:quad}
 [a_\lambda b]=\partial u(a,b)+\lambda v(a,b)+w(a,b),
 \qquad u(a,b),v(a,b),w(a,b)\in A.
\end{equation}
Thus $p$ is conformal weight minus one. This definition does not prescribe an associative multiplication on $A$.
\end{definition}

For a primary quadratic presentation, write $H=\ell_{(1)}$ for the
conformal-weight operator on $C$. Thus $H|_A=P+I$.

\begin{lemma}\label{cy:lem:covariance}
For $a\in A_p$ and $b\in A_q$ in a primary quadratic presentation,
\[
 u(a,b),v(a,b)\in A_{p+q-1},\qquad w(a,b)\in A_{p+q},
\]
and
\begin{equation}\label{cy:eq:balance}
 (p+q)u(a,b)=p\,v(a,b).
\end{equation}
\end{lemma}
\begin{proof}
The coefficient of the first power of the Virasoro parameter in the conformal Jacobi identity gives
\[
 H(a_{(j)}b)=(Ha)_{(j)}b+a_{(j)}Hb-(j+1)a_{(j)}b,
 \qquad [H,\partial]=\partial.
\]
Since $C$ is free on $A$, comparison of independent powers of $\partial$ in \eqref{cy:eq:quad} gives the stated weights.

Apply the Jacobi identity to $\ell,a,b$, using $x$ for the Virasoro parameter and $\lambda$ for the bracket of $a,b$:
\[
 [\ell_x[a_\lambda b]]
 =(-\lambda+px)[a_{\lambda+x}b]
  +(\partial+\lambda+(q+1)x)[a_\lambda b].
\]
The coefficient of $x^2$ on the left is $(p+q)u(a,b)$, and on the right it is $p\,v(a,b)$. This proves \eqref{cy:eq:balance}. No division by $p$, $q$, or $p+q$ is used.
\end{proof}

\begin{theorem}\label{cy:thm:normalization}
Let $C=\C[\partial]\otimes A$ have a primary quadratic presentation. Let $\mu$ be a diagonalizable conformal automorphism satisfying
\begin{equation}\label{cy:eq:hyp}
 \mu\ell=\ell,\qquad \mu(A)=A.
\end{equation}
For $r\ge1$, define $\kappa_r$ on $A_p$ by
\[
 \kappa_r(a)=e^{2\pi i(r-1)p}a
\]
and extend it $\C[\partial]$-linearly. Then $\kappa_r$ is a conformal automorphism commuting with $\mu$, and $\tau=\mu\kappa_r$ is diagonalizable. There is an isomorphism
\begin{equation}\label{cy:eq:main}
 \Phi_r:\Alg_r(C,\mu)\longrightarrow\Alg(C,\tau)
\end{equation}
given, for a simultaneous eigenvector $a\in A_p$, by
\begin{equation}\label{cy:eq:map}
 \Phi_r(a\nn{m})=r^{-p}a[m-(r-1)p].
\end{equation}
Here $r^{-p}=\exp(-p\log r)$ with the real logarithm of the positive integer $r$. The map preserves parity and sends the center isomorphically onto the center.
\end{theorem}
\begin{proof}
The automorphism $\mu$ commutes with $H$ and therefore with $P$ on $A$. The weights in Lemma~\ref{cy:lem:covariance} show that $\kappa_r$ preserves \eqref{cy:eq:quad}: the discrepancy of one in the weights of $u,v$ contributes the factor $e^{-2\pi i(r-1)}=1$. It also commutes with $\partial$ by definition. Commuting semisimple operators on the finite-dimensional space $A$ are simultaneously diagonalizable, proving the assertion about $\tau$.

If $\mu a=e^{-2\pi i m}a$, then
\[
 \tau a=e^{-2\pi i(m-(r-1)p)}a.
\]
Thus the target index in \eqref{cy:eq:map} is allowed. Because $C=\C[\partial]\otimes A$ is free, the translation relations give
explicit normal forms. For $s\ge0$,
\[
 (\partial^s a)\nn{m}=(-1)^s r^{-s}
 \prod_{j=0}^{s-1}(m-rj)\,a\nn{m-rs},
\]
while in the ordinary coefficient algebra
\[
 (\partial^s a)[k]=(-1)^s
 \prod_{j=0}^{s-1}(k-j)\,a[k-s].
\]
Choose a simultaneous eigenbasis of $A$. Define the normal-form map from
the free span of all derivative-mode symbols to the free span of the
undifferentiated generator modes by the displayed formulas. It kills
every translation relation and is the identity on the undifferentiated
symbols. Hence the undifferentiated symbols form a basis of the translation
quotient. Therefore \eqref{cy:eq:map} is a bijection on
coefficient vector spaces. Its inverse on generator modes is
\[
 a[k]\longmapsto r^p a\nn{k+(r-1)p}.
\]
Derivative states on both sides are then determined uniquely by their
respective coefficient relations.

For $a\in A_p,b\in A_q$, denote the three coefficients in \eqref{cy:eq:quad} by $u,v,w$. In the domain,
\begin{equation}\label{cy:eq:src}
 [a\nn{m},b\nn{n}]
 =\frac1r\bigl(mv-(m+n)u\bigr)\nn{m+n-r}+w\nn{m+n}.
\end{equation}
Set $M=m-(r-1)p$ and $N=n-(r-1)q$. The target bracket of the two images is
\[
 r^{-p-q}\left(\bigl(Mv-(M+N)u\bigr)[M+N-1]+w[M+N]\right).
\]
By \eqref{cy:eq:balance},
\begin{align*}
 Mv-(M+N)u
 &=mv-(m+n)u-(r-1)\bigl(pv-(p+q)u\bigr)\\
 &=mv-(m+n)u.
\end{align*}
The source term involving $u,v$ has $P$-weight $p+q-1$, so its scalar after applying $\Phi_r$ is
\[
 r^{-1}r^{-(p+q-1)}=r^{-p-q},
\]
and its index is
\[
 m+n-r-(r-1)(p+q-1)=M+N-1.
\]
The $w$ term has weight $p+q$, so its scalar and index also agree. This proves bracket preservation. Center preservation follows from the isomorphism.
\end{proof}

\begin{corollary}\label{cy:cor:cyclic}
Under the hypotheses of Theorem~\ref{cy:thm:normalization},
\[
 \Alg(C^{\oplus r},\sigma)\simeq
 \Alg\bigl(C,\mu\exp(2\pi i(r-1)P)\bigr).
\]
This isomorphism induces an isomorphism of the quotients by the full
centers.
\end{corollary}

On the distinguished Virasoro modes, set
\[
 L_n^{\mathrm{diag}}=(\ell,\ldots,\ell)[n+1],\qquad L_n=\ell[n+1].
\]
The composite of \eqref{cy:eq:fourier} and \eqref{cy:eq:map} sends
\begin{equation}\label{cy:eq:virimage}
 L_n^{\mathrm{diag}}\longmapsto\frac1r L_{rn}.
\end{equation}
For $r>1$, the map does not identify the two marked Virasoro copies.

\subsection{Invariant primary spaces}
A state $a$ is \emph{primary of weight $h$} if
$[\ell_\lambda a]=(\partial+h\lambda)a$. Denote the space of these
states by $\Prim_h C$. In the following standard presentations,
the condition $\mu(A)=A$ follows from $\mu\ell=\ell$.

\begin{lemma}\label{cy:lem:prim}
Suppose $C=\C[\partial]\otimes A$ satisfies \eqref{cy:eq:primary}. Let $A^{(d)}$ be the space of generators of conformal weight $d$, so $A^{(d)}=A_{d-1}$. Then the primary states of weight $h$ are
\[
 \Prim_h C=A^{(h)}\quad(h\ne1),\qquad
 \Prim_1 C=A^{(1)}\oplus\partial A^{(0)}.
\]
In particular, if $A^{(0)}=0$, every conformal automorphism fixing $\ell$ preserves $A$.
\end{lemma}
\begin{proof}
A primary state is an eigenvector of $H$. Write it with respect to a
homogeneous basis of $A$. Each nonzero component is a scalar multiple
of $\partial^k a$, where $a$ has weight $d$ and $d+k=h$.
The primary condition becomes
\[
 (\partial+\lambda)^k(\partial+d\lambda)
 =\partial^{k+1}+(d+k)\lambda\partial^k.
\]
For $k\ge1$, the coefficient of $\lambda^{k+1}$ forces $d=0$. If $k\ge2$, the coefficient of $\partial\lambda^k$ then gives a contradiction. The only additional possibility is $k=1,d=0$, which gives weight one. Independence of the free basis prevents cancellation between distinct coordinates.
\end{proof}

Reduced quasi-Poisson presentations from \cite[Theorem 6.10(c)]{Billig} have
\begin{equation}\label{cy:eq:qp}
 [a_\lambda b]=\partial(P(a)b)
 +\lambda(P(a)b+aP(b))+\beta(a,b).
\end{equation}
Here juxtaposition and $\beta$ denote the two bilinear operations on
$A$ in that presentation.
For $P(a)=pa$, $P(b)=qb$, this is \eqref{cy:eq:quad} with $u=pab$, $v=(p+q)ab$, and $w=\beta(a,b)$. Thus the normalization theorem applies.

\subsection{A characteristic current ideal}\label{sec:normalization}
\begin{lemma}\label{ck:lem:currentcore}
Let $C$ be a Lie conformal superalgebra whose even part is
\[
 C_{\overline0}=\C[\partial]\ell\ltimes\Cur\mathfrak g,
 \qquad [\ell_\lambda a]=(\partial+\lambda)a\quad(a\in\mathfrak g),
\]
where $\mathfrak g$ is a nonzero finite-dimensional semisimple Lie algebra.
Then every parity-preserving conformal automorphism of $C$ fixes $\ell$.
\end{lemma}
\begin{proof}
The centralizer of $\Cur\mathfrak g$ in $C_{\overline0}$ is zero. If an
 element $p(\partial)\ell+b$ centralizes all constant $a\in\mathfrak g$,
the coefficient of $\partial a$ in its bracket is $p(-\lambda)$, whereas
$[b_\lambda a]$ contains no $\partial$ term. Thus $p=0$; then
centerlessness of $\mathfrak g$ gives $b=0$.

Every minimal ideal of $C_{\overline0}$ lies in $\Cur\mathfrak g$: otherwise
its intersection with that ideal is zero, so it centralizes it. Conversely,
the current algebras of the simple summands of $\mathfrak g$ are minimal
ideals. Therefore $\Cur\mathfrak g$ is the socle of the even algebra and is
characteristic.

A conformal automorphism of $\Cur\mathfrak g$ preserves its constant
subspace, by \cite[Theorem 3.4]{KLP}; the hypothesis that every ideal of
$\mathfrak g$ is centerless holds here. The quotient by the current ideal is
$\operatorname{Vir}$, whose conformal automorphism group is trivial: a free rank-one
module automorphism is multiplication by a scalar, and the bracket forces
that scalar to be $1$. Thus any automorphism $\varphi$ of $C$ satisfies
$\varphi(\ell)=\ell+b$, $b\in\Cur\mathfrak g$. For every constant $a$,
\[
 [\ell+b{}_\lambda\varphi(a)]=(\partial+\lambda)\varphi(a)
 =[\ell_\lambda\varphi(a)].
\]
Surjectivity on $\mathfrak g$ and its zero centralizer give $b=0$.
\end{proof}

The hypothesis specifies the entire even algebra. It applies in particular to $CK_6$, whose even
part is $\operatorname{Vir}\ltimes\Cur\mathfrak{so}_6$ \cite[Example 3.9]{FKR}.
It also applies to the standard $K_3$ presentation and to the usual
$S_{2,0}$ presentation with even current algebra $\mathfrak{sl}_2$;
see \cite[Section 3.4]{KLP} for the latter. Their automorphisms therefore fix the specified Virasoro state.

\subsection{Current transformations}\label{sec:current-transform}
The following lemma is a coefficient-algebra version of the logarithmic
pseudo-derivation construction in \cite[Propositions 3.11 and 3.19]{Li}.

\begin{lemma}\label{ck:lem:flow}
Let $C$ be a Lie conformal superalgebra with an even Virasoro state
$\ell$ satisfying $[\ell_\lambda\ell]=(\partial+2\lambda)\ell$ and
$\ell_{(0)}=\partial$. Suppose that $H=\ell_{(1)}$ is semisimple and
its weights are bounded below in each of finitely many integer cosets.
Let $h$ be even, primary of weight $1$, with
$[h_\lambda h]=0$, and suppose that $h_{(0)}$ is diagonalizable. If $\sigma$
is a diagonalizable automorphism fixing $h$, define
$z^{h_{(0)}}a=z^qa$ whenever $h_{(0)}a=qa$, and set
\[
 \tau=\sigma e^{-2\pi i h_{(0)}},\qquad
 A_h(z)=\sum_{j\ge1}\frac{(-1)^{j-1}}j z^{-j}h_{(j)},\qquad
 \Gamma_h(z)=z^{h_{(0)}}e^{A_h(z)}.
\]
Then $\Gamma_h$ induces a Lie superalgebra isomorphism
\[
 \Alg(C,\sigma)\simeq\Alg(C,\tau).
\]
The sums defining $A_h(z)a$ and $e^{A_h(z)}a$ are finite for each state.
The isomorphism preserves centers; its inverse is induced by
$\Gamma_{-h}$.
\end{lemma}
\begin{proof}
The operator $h_{(j)}$ lowers conformal weight by $j$. Locality makes the
sum finite on each state, and the lower bounds make $A_h$ locally
nilpotent. The conformal commutator formula and $[h_\lambda h]=0$ give
$[h_{(i)},h_{(j)}]=0$ for $i,j\ge0$.

Adjoin a formal logarithm $\zeta$ and all complex powers of $z$ generated
by the eigenvalues in question, with $d\zeta/dz=z^{-1}$. In the differential
scalar extension the zeroth product of $\zeta h$ is the derivation
\[
 \mathsf X_h=\zeta h_{(0)}+A_h(z).
\]
It commutes with the total translation $\partial+d/dz$. The two summands
commute. Its exponential is well defined after adjoining exponentials of
$\zeta$: the first factor is diagonal and the second is finite on each
state. The differential substitution $e^{q\zeta}\mapsto z^q$ gives precisely
$\Gamma_h$. Thus $\Gamma_h$ is an automorphism of the differential scalar
extension commuting with total translation. This can equivalently be
verified by exponentiating the derivation identity first as a formal
identity in an auxiliary scalar and then using the finite weight-lowering
part. The inverse is obtained by replacing $h$ by $-h$.

The operators $\sigma$ and $h_{(0)}$ commute. For a joint eigenvector with
$\sigma a=e^{-2\pi i m}a$ and $h_{(0)}a=qa$, every term of
$e^{A_h(z)}a$ has the same eigenvalues. Write
\[
 e^{A_h(z)}a=\sum_j z^{-j}a_j.
\]
The affinization map sends $z^ma$ to
$\sum_j z^{m+q-j}a_j$. Since
\[
 \tau a_j=e^{-2\pi i(m+q-j)}a_j,
\]
it sends exactly the $\sigma$-twisted exponent sectors to the
$\tau$-twisted sectors. Commutation with total translation lets it descend
to the coefficient quotients, where it preserves the bracket induced by
the zeroth product. Its inverse respects the reverse sectors.
\end{proof}

\section{The Witt family}\label{sec:endpoints}
\subsection{Automorphisms fixing the standard Virasoro state}
\begin{theorem}\label{cy:thm:witt}
Let $C=W_N$ in its reduced primary presentation, with standard Virasoro state $\ell=1$. If $\mu$ is a diagonalizable conformal automorphism fixing $1$, then
\[
 \Alg(W_N,\mu)\simeq W(1,N).
\]
Consequently, if $\sigma$ is the corresponding cyclic automorphism of $W_N^{\oplus r}$, then
\[
 \Alg(W_N^{\oplus r},\sigma)\simeq W(1,N).
\]
\end{theorem}
\begin{proof}
Write $B=\Lambda(\xi_1,\ldots,\xi_N)$ and $A=B\oplus\Der B$. The $P$-eigenspaces are $B$ with eigenvalue $1$ and $\Der B$ with eigenvalue $0$. Lemma~\ref{cy:lem:prim} and $\mu\ell=\ell$ imply that $\mu$ preserves both. Comparing the first product on two functions shows that $\phi=\mu|_B$ is a unital associative automorphism of $B$. Comparing the reduced mixed bracket gives
\[
 (\mu X)(\phi f)=\phi(X(f)),\qquad
 \mu X=\phi X\phi^{-1}.
\]
Conversely every even associative automorphism of $B$ extends in this way, preserving the presentation in \cite[Appendix A, type $W_N$]{Billig}.

Let $\mathfrak m$ be the augmentation ideal of $B$. It is characteristic. Since $\phi$ is semisimple, there is a $\phi$-invariant complement to $(\mathfrak m^2)_{\overline1}$ in $\mathfrak m_{\overline1}$. Choose an eigenbasis of this complement. These odd lifts freely generate $B$: the induced map on $\mathfrak m/\mathfrak m^2$ is invertible, and the induced map on the associated exterior graded algebra is invertible. After an exterior coordinate change, therefore,
\[
 \phi\xi_i=e^{-2\pi i b_i}\xi_i\qquad(b_i\in\C).
\]
For $I\subseteq\{1,\ldots,N\}$ write $b_I=\sum_{i\in I}b_i$ and $\xi_I=\prod_{i\in I}\xi_i$ in increasing order.

In a ring allowing the required complex powers of $t$, consider
\[
 \mathcal R=\bigoplus_I t^{b_I}\C[t,t^{-1}]\xi_I.
\]
With $\eta_i=t^{b_i}\xi_i$, this is the ordinary Laurent exterior algebra
\[
 \mathcal R=\C[t,t^{-1}]\otimes\Lambda(\eta_1,\ldots,\eta_N).
\]
The twisted coefficient generators have the vector-field realization
\[
 f[m]\longmapsto-t^m f\partial_t,\qquad
 X[m]\longmapsto t^mX.
\]
For $f=\xi_I$, its mode set is $b_I+\Z$; for $X=\xi_I\partial_i$, it is $b_I-b_i+\Z$. The displayed vector fields preserve $\mathcal R$, and their bracket agrees with the coefficient bracket by direct differentiation.

They span every derivation of $\mathcal R$. Indeed,
\[
 \partial_{\eta_i}=t^{-b_i}\partial_i,
 \qquad
 \left.\partial_t\right|_{\eta}
 =\left.\partial_t\right|_{\xi}-t^{-1}\sum_i b_i\xi_i\partial_i
\]
belong to their span. Their multiples by every Laurent exterior monomial also belong to that span, by the allowed mode cosets. Independence follows from the independent Laurent exterior monomials in these vector fields. Thus the realization is an isomorphism onto $\Der\mathcal R\simeq W(1,N)$.

Finally $P$ has integer eigenvalues, so $\kappa_r=\id$, and Corollary~\ref{cy:cor:cyclic} proves the cyclic assertion.
\end{proof}
In the last coordinate identification, the standard Witt generators become
\[
 -t^{n+1}\left.\partial_t\right|_{\xi}
 =-t^{n+1}\left.\partial_t\right|_{\eta}
  -t^n\sum_i b_i\eta_i\partial_{\eta_i}.
\]

\subsection{Diagonalizable automorphisms outside rank one}
\label{witt:sec:all}

The geometric automorphism theorem for the formal vector-field algebra
has an exceptional case $W(1,1)$, which we treat separately.
Write
\[
 B=\Lambda(\xi_1,\ldots,\xi_N),\qquad
 \mathcal O_N=\C[[t]]\otimes B,\qquad
 \widehat{\mathfrak W}_N=\operatorname{Der}_{\rm cont}\mathcal O_N.
\]
Here the hat distinguishes the completed annihilation algebra from the
Laurent vector-field algebra denoted by $W(1,N)$ in
Theorem~\ref{cy:thm:witt}.

\begin{lemma}\label{witt:lem:lift}
For $N\geq2$, every even conformal automorphism of $W_N$ is induced on
its completed annihilation algebra by a unique coordinate automorphism
\begin{equation}\label{witt:eq:coords}
 \Phi(t)=t+f(\xi),\qquad
 \Phi(\xi_i)=\phi(\xi_i),
 \quad f\in B_{\overline0}\cap\mathfrak m,\quad
 \phi\in\operatorname{Aut}_{\overline0} B,
\end{equation}
where $\mathfrak m$ is the augmentation ideal of $B$.
Conversely, every coordinate automorphism in
\eqref{witt:eq:coords} has an even conformal lift.  Consequently
\[
 \operatorname{Aut}_{\overline0}W_N
 \simeq(B_{\overline0}\cap\mathfrak m)\rtimes
             \operatorname{Aut}_{\overline0}B.
\]
\end{lemma}
\begin{proof}
Use the ambient convention of \cite[Example 3.4]{FK}, in which
$W_N=\C[\partial]\otimes(\operatorname{Der}B\oplus B)$ and
\[
 X[n]=t^nX,\qquad f[n]=t^nf\partial_t,
 \qquad \ell_{\rm std}=-1.
\]
The completed annihilation algebra is $\widehat{\mathfrak W}_N$.
A conformal automorphism induces a continuous automorphism commuting
with $\operatorname{ad}\partial_t$. The geometric automorphism result in the
$W(m,n)$ case of the proof of
\cite[Theorem 4.5]{CantariniKacAut} applies because $(1,N)\ne(1,1)$:
the induced automorphism is conjugation by a continuous coordinate
automorphism $\Phi$ of $\mathcal O_N$.

The center of $\widehat{\mathfrak W}_N$ is zero.  Commutation with
$\operatorname{ad}\partial_t$ therefore gives
$\Phi\partial_t\Phi^{-1}=\partial_t$.  On the coordinates this says
\[
 \partial_t\Phi(t)=1,\qquad \partial_t\Phi(\xi_i)=0.
\]
Thus $\Phi$ has form \eqref{witt:eq:coords}. Continuity
excludes a nonzero scalar term in $f$, and parity makes $f$ even.
The map induced on the odd linear quotient is invertible, so the
restriction $\phi$ is an exterior-algebra automorphism.  A coordinate
automorphism acting trivially on all vector fields is the identity:
commutation with the constant partial derivatives first restricts it
to a scalar translation, and commutation with $t\partial_t$ eliminates
that translation.  Hence the coordinate description is unique.

The exterior change $\phi$ has the conformal lift
\[
 g\longmapsto\phi(g),\qquad
 X\longmapsto\phi X\phi^{-1}.
\]
For a shift $t\mapsto t+q(\xi)$, its lift in this ambient convention is
\begin{equation}\label{witt:eq:shift}
 \exp(q_{(0)}).
\end{equation}
Indeed, $q[0]=q\partial_t$, so the induced annihilation action is
$\exp(\operatorname{ad}(q\partial_t))$.  This is precisely conjugation
by the coordinate shift.  The exponential is finite already on the
conformal algebra.  Give $g\in B$ its exterior degree, give
$g\partial_i$ degree $\deg g-1$, and give $\partial$ degree zero.
The ambient brackets give, for even $q$,
\[
 q_{(0)}g=-\partial(qg),\qquad
 q_{(0)}X=-X(q)-\partial(qX).
\]
Every positive-degree even component of $q$ has degree at least two;
each displayed term raises that auxiliary degree by at least two.
Only finitely many such degrees occur.  Thus $q_{(0)}$ is nilpotent,
commutes with $\partial$, and its finite exponential is a conformal
automorphism.  Composing these two types of lifts covers every
\eqref{witt:eq:coords}.

Finally the action of conformal automorphisms on the completed
annihilation algebra is faithful.  For
$v=\sum_{j=0}^d\partial^jv_j$, with
$v_j\in\operatorname{Der}B\oplus B$, one has, for $n\geq d$,
\[
 v[n]=\sum_{j=0}^d(-1)^j n(n-1)\cdots(n-j+1)v_j[n-j].
\]
In the vector-field realization the different $j$ have different
$t$-degrees, and the scalar factors are nonzero.  Vanishing of these
modes forces every $v_j=0$.  The group law of the coordinate
description is $(f,\phi)(g,\psi)=(f+\phi(g),\phi\psi)$, proving
the asserted semidirect product.
\end{proof}

\begin{theorem}\label{witt:thm:semisimple}
If $N\geq2$, every diagonalizable even conformal automorphism of $W_N$
is conformally conjugate to an exterior linear diagonal automorphism
fixing the standard Witt state.  In particular, for every such $\mu$,
and for every diagonalizable cyclic automorphism $\sigma$ of
$W_N^{\oplus r}$ with return automorphism $\mu$,
\[
 \Alg(W_N,\mu)\simeq W(1,N),\qquad
 \Alg(W_N^{\oplus r},\sigma)\simeq W(1,N).
\]
The same coefficient conclusions hold for $N=0$.
\end{theorem}
\begin{proof}
Because $W_N$ is finitely generated over $\C[\partial]$, a
diagonalizable conformal automorphism $\mu$ has only finitely many
eigenvalues: decompose finitely many module generators into
eigenvectors and use $\mu\partial=\partial\mu$.  Consequently
$p(\mu)=0$ for a squarefree polynomial $p$.
The induced annihilation automorphism satisfies the same polynomial
on polynomial vector fields, and therefore on their completion by
continuity.

Let $\Phi$ be given by Lemma~\ref{witt:lem:lift}.  The finite
dimensional space
\[
 E=(\C t\oplus B)\partial_t
\]
is invariant, since $\Phi$ fixes $\partial_t$, and its action is
$g\partial_t\mapsto\Phi(g)\partial_t$.  Its action is semisimple
by the squarefree polynomial just obtained.  On
$E/(B\partial_t)$ the action is the identity.  A semisimple
equivariant splitting gives an invariant even lift
\[
 \tau=t+q(\xi),\qquad \Phi(\tau)=\tau.
\]
Subtracting its scalar constant makes $q\in B_{\overline0}\cap
\mathfrak m$.  The same polynomial implies that
$\phi=\Phi|_B$ is semisimple.  Take odd eigenvector lifts
$\eta_i\in\mathfrak m_{\overline1}$ of an eigenbasis of
$\mathfrak m/\mathfrak m^2$.  They generate $B$, since their
linear terms are independent.  Thus $(\tau,\eta_1,\ldots,\eta_N)$
are coordinates, commuting with $\partial_t$, in which
\[
 \Phi(\tau)=\tau,\qquad \Phi(\eta_i)=\lambda_i\eta_i.
\]
The coordinate change has a conformal lift by
Lemma~\ref{witt:lem:lift}; faithfulness there proves that it
conjugates $\mu$ itself to the indicated diagonal automorphism.
This fixes the standard Witt state.  The change between the ambient
and reduced Witt conventions sends that state $-1$ to the state $1$
of Theorem~\ref{cy:thm:witt}.  That theorem and
Corollary~\ref{cy:cor:cyclic} now apply; the Witt cyclic phase is the
identity.  For $N=0$, $W_0=\operatorname{Vir}$ and its conformal
automorphism group is trivial, so the same conclusion follows
directly.
\end{proof}

\subsection{The rank-one exception and its two coefficient forms}
\label{witt:sec:rankone}

Use the reduced Witt presentation of $W_1$ in \cite[Theorem 6.10(c), Appendix A]{Billig}.
It is free over $\C[\partial]$ on even states $L,J$ and odd states $x,y$, with
\begin{align*}
 [L_\lambda L]&=(\partial+2\lambda)L,&
 [L_\lambda J]&=(\partial+\lambda)J,&[J_\lambda J]&=0,\\
 [L_\lambda x]&=(\partial+2\lambda)x,&
 [L_\lambda y]&=(\partial+\lambda)y,\\
 [J_\lambda x]&=x,&[J_\lambda y]&=-y,\\
 [x_\lambda y]&=L+(\partial+\lambda)J,&
 [y_\lambda x]&=L-\lambda J,&
 [x_\lambda x]&=[y_\lambda y]=0.
\end{align*}
The remaining brackets follow from super-skew-symmetry. In this convention,
$L$ corresponds to the function $1$, $x$ to $\xi$, $y$ to $\partial_\xi$, and
$J$ to $\xi\partial_\xi$.

\begin{proposition}\label{witt:prop:rankone-aut}
Every even conformal automorphism of $W_1$ is one of
\begin{align*}
 g_a:&\quad L\mapsto L,\quad J\mapsto J,\quad
          x\mapsto ax,\quad y\mapsto a^{-1}y,\\
 r_a:&\quad L\mapsto L+\partial J,\quad J\mapsto-J,\quad
          x\mapsto ay,\quad y\mapsto a^{-1}x,
 \qquad a\in\C^\times.
\end{align*}
In particular,
\[
 \operatorname{Aut}_{\overline0}W_1
       =\C^\times\rtimes\langle\rho\rangle,
 \qquad \rho=r_1,
\]
and every automorphism fixes $\ell=L+\tfrac12\partial J$.
All these automorphisms are diagonalizable.  Every $r_a$ is conjugate
to $\rho$.
\end{proposition}
\begin{proof}
The solvable radical of the even conformal algebra is
$\C[\partial]J$, and its quotient is $\operatorname{Vir}$.
The radical is characteristic, and an automorphism of its rank-one
free module is multiplication by a nonzero constant.  Every
conformal automorphism of $\operatorname{Vir}$ is the identity.
Thus an arbitrary $\mu$ satisfies
\[
 \mu J=cJ,\qquad \mu L=L+P(\partial)J,
 \qquad c\in\C^\times.
\]
On the odd module, $J_{(0)}$ has the two eigenspaces
$\C[\partial]x$ and $\C[\partial]y$, with eigenvalues $1,-1$.
Intertwining this action forces $c=1$ or $c=-1$.
For $c=1$ the two rank-one modules are preserved; for $c=-1$ they
are exchanged.  The nonzero diagonal or antidiagonal entries of the
resulting invertible matrix over $\C[\partial]$ must be units,
and hence constants.  Accordingly, either
$(\mu x,\mu y)=(ax,by)$ or $(ay,bx)$.

Apply $\mu$ to $[x_\lambda y]=L+(\partial+\lambda)J$.
In the first case its comparison with $ab[x_\lambda y]$ gives
$ab=1$ and $P=0$.  In the second case comparison with
$ab[y_\lambda x]=ab(L-\lambda J)$ gives $ab=1$ and
$P=\partial$.  Conversely the displayed bracket table verifies
both types.  These formulas give the stated group and the invariant
state $\ell$.  Each $g_a$ is diagonal and each $r_a$ squares to the
identity.  Finally $g_d^{-1}r_ag_d=r_{ad^2}$, so a square root of
$a^{-1}$ conjugates $r_a$ to $\rho$.
\end{proof}

\begin{proposition}\label{ck:prop:obstruction}
For every $c\in\C$, the state $L_c=L+c\partial J$ is Virasoro and
$(L_c)_{(0)}=\partial$. The states $L$ and $L_{1/2}$ are not conjugate under
conformal automorphisms. Moreover, the involution
\begin{equation}\label{ck:eq:involution}
 \rho(L)=L+\partial J,\qquad \rho(J)=-J,\qquad
 \rho(x)=y,\qquad \rho(y)=x
\end{equation}
is not conjugate to any automorphism fixing $L$.
\end{proposition}
\begin{proof}
Sesquilinearity gives
\[
 [L_\lambda\partial J]+[\partial J_\lambda L]
 =(\partial+\lambda)^2J-\lambda^2J
 =(\partial+2\lambda)\partial J.
\]
Since $[J_\lambda J]=0$, this proves the Virasoro identity for $L_c$.
A derivative in the first argument has zero zeroth product, so its zeroth
product is still $\partial$. On the odd generators,
\[
 [ (L_c)_\lambda x]=(\partial+(2-c)\lambda)x,\qquad
 [ (L_c)_\lambda y]=(\partial+(1+c)\lambda)y.
\]
For $c=0$, the spectrum of $(L_c)_{(1)}$ on the odd part is contained in
$\Z$, whereas for $c=1/2$ it is contained in $3/2+\Z$. A parity-preserving
conformal automorphism conjugating these states would conjugate these
endomorphisms, which is impossible.

By Proposition~\ref{witt:prop:rankone-aut}, the automorphisms fixing
$L$ are the $g_a$, whereas $\rho=r_1$. The scalar by which an
automorphism acts on $J$ is a character of the automorphism group,
and hence is invariant under conjugation. It is $1$ for every $g_a$
and $-1$ for $\rho$. Thus $\rho$ is not conjugate to an automorphism
fixing $L$.
\end{proof}

For $\ell=L+\tfrac12\partial J$, the generators have weights
$2,1,3/2,3/2$, respectively, and
\begin{equation}\label{witt:eq:k2table}
 [x_\lambda y]=\ell+(\tfrac12\partial+\lambda)J,
 \qquad [\ell_\lambda x]=(\partial+\tfrac32\lambda)x,
 \qquad [\ell_\lambda y]=(\partial+\tfrac32\lambda)y.
\end{equation}
This is the primary contact presentation $K_2$. Under the identification
$\ell\leftrightarrow1$ and $J\leftrightarrow xy$, the odd generators are
$x,y$, the algebra $\Lambda(x,y)$ has
$\{x,y\}=1$, $\{x,x\}=\{y,y\}=0$, and $P=1-E/2$.
The standard contact bracket is
\[
 [a_\lambda b]
 =\partial\bigl((Pa)b\bigr)
  +\lambda\bigl((Pa)b+a(Pb)\bigr)+\{a,b\}.
\]
It gives the entire table, including
$[J_\lambda x]=x$ and $[J_\lambda y]=-y$.
The involution $\rho$ exchanges the isotropic odd coordinates
$x,y$, and is an orthogonal reflection in this two-dimensional
contact presentation.

\begin{theorem}\label{witt:thm:rankone-forms}
The full coefficient algebras of diagonalizable twists of $W_1$
have exactly the following two parity-preserving isomorphism types:
\[
 \mathfrak k_2^{\rm id}:=\Alg(K_2,\operatorname{id}),
 \qquad
 \mathfrak k_2^{\rm ref}:=\Alg(K_2,\rho).
\]
They are nonisomorphic, both have zero full center, and the same
two possibilities exhaust cyclic orbits $W_1^{\oplus r}$.
The first is the usual Laurent $W(1,1)$ coefficient algebra.
The second is the reflection-twisted contact form.
\end{theorem}
\begin{proof}
Proposition~\ref{witt:prop:rankone-aut} reduces a single-factor twist
to $g_c$ or $\rho$.  Write $c=e^{2\pi i a}$.
The current $h=aJ$ is even and primary of weight one, has
$[h_\lambda h]=0$, and has semisimple zeroth product.  It is fixed
by $g_c$, and $g_c=e^{2\pi i h_{(0)}}$.
Lemma~\ref{ck:lem:flow} therefore removes the $g_c$ twist.
The induced mode map in the $\ell$ presentation is
\begin{equation}\label{witt:eq:k2flow}
 x[m]\mapsto x[m+a],\quad y[n]\mapsto y[n-a],\quad
 J[k]\mapsto J[k],\quad
 \ell[k]\mapsto\ell[k]+aJ[k-1].
\end{equation}
In the source, the allowed $x$ modes are $-a+\Z$, the $y$ modes
are $a+\Z$, and the even modes are integral.  The target modes
in \eqref{witt:eq:k2flow} are all integral.  For example,
\[
 [x[m],y[n]]
  =\ell[m+n]+\tfrac12(m-n)J[m+n-1],
\]
and the two mode shifts give precisely the extra $aJ[m+n-1]$
in the image of $\ell[m+n]$.  The other brackets follow either
from the lemma or directly from the same table. The identity twist is also the untwisted
coefficient algebra of $W_1$, hence $W(1,1)$.

To distinguish the two outcomes, consider the center of the even
part, an invariant of the parity-graded Lie superalgebra.
Set $L_n=\ell[n+1]$ and $J_m=J[m]$. In either case the even
part is generated by these integral Witt modes and the $J$ modes,
with
\[
 [L_n,J_m]=-mJ_{n+m},\qquad [J_m,J_k]=0.
\]
For the identity twist $m\in\Z$; for the reflection twist
$m\in\tfrac12+\Z$.  The Witt quotient is centerless, so a
central element of the even part lies in the current ideal.
The displayed brackets then give
\[
 Z((\mathfrak k_2^{\rm id})_{\overline0})=\C J_0,
 \qquad Z((\mathfrak k_2^{\rm ref})_{\overline0})=0.
\]
Thus the two superalgebras are not isomorphic.  In the identity
case $J_0$ acts with nonzero opposite eigenvalues on the odd
generators and hence is not central in the full superalgebra;
its action also excludes an odd center.  In the reflection case
set $u=x+y$ and $v=x-y$.  Then $u$ has integral modes, $v$ has
half-integral modes, and any allowed $J_s$ maps $u[m]$ to
$v[s+m]$ and $v[n]$ to $u[s+n]$. This is injective on their
direct sum, excluding an odd center there as well.

Finally this primary presentation is invariant under every
automorphism.  Its cyclic phase is
$\kappa_r=g_{(-1)^{r-1}}$, since $P$ has eigenvalues
$1,0,1/2,1/2$ on $\ell,J,x,y$.
The cyclic coefficient theorem therefore yields the twist
$\mu\kappa_r$.  This lies in the same scaling or reflection
component as $\mu$.  The preceding analysis gives exactly the
same two forms for every cycle length.
\end{proof}

\section{The special family}\label{sec:special}

\subsection{Primary presentations}
Throughout this section $S_N=S_{N,0}$, $N\ge2$, in the vector-field
convention of Subsection~\ref{ft:sec:translation}. The distinguished Virasoro
state is
\[
 \ell_*=-1+\frac1N\partial E,
 \qquad E=\sum_{i=1}^N\xi_i\partial_i.
\]

\begin{lemma}[An invariant primary quadratic generating space]
\label{sp:lem:primary}
Let $S(N)$ be the exterior divergence-zero subalgebra of $W(N)$. For a
homogeneous $f\in\Lambda_N$ of degree $d<N$, set
\[
 e_f=f-\frac{\partial(fE)}{N-d}.
\]
Then
\[
 A_N=S(N)\oplus\Span\{e_f:\deg f<N\},\qquad
 S_N=\C[\partial]\otimes A_N
\]
is a $\GL_N$-invariant primary quadratic presentation for $\ell_*=-e_1$.
If $X\in S(N)$ has vector-field degree $r$, the conformal weights are
\[
 \Delta(X)=1-\frac rN,\qquad
 \Delta(e_f)=2-\frac dN.
\]
All these weights are positive.
\end{lemma}
\begin{proof}
The exterior divergence satisfies $\operatorname{div}(fE)=(d-N)f$.
This gives the asserted free decomposition of $\ker\Div$. It is
invariant under linear coordinate changes because $E$ is invariant.
A direct calculation gives
\[
 [\ell_*{}_{\lambda}X]
   =\left(\partial+\left(1-\frac rN\right)\lambda\right)X,
 \qquad
 [\ell_*{}_{\lambda}e_f]
   =\left(\partial+\left(2-\frac dN\right)\lambda\right)e_f.
\]
In the second identity, the two apparent $\lambda^2fE$ terms cancel.
The maximum degree of a nonzero field in $S(N)$ is $N-2$;
therefore its weights are at least $2/N$.

The generator brackets are as follows. If $\deg f=d$, $\deg g=e$, and
$d+e<N$, then
\begin{align*}
 [e_f{}_{\lambda}e_g]
 ={}&-\frac{N-d-e}{N-e}\partial e_{fg}\\
 &-\frac{(N-d-e)(2N-d-e)}{(N-d)(N-e)}\lambda e_{fg}.
\end{align*}
If $d+e\ge N$, this bracket is zero. At $d+e=N$ this follows from
$fgE=0$ and cancellation of the function coefficient; no division by
$N-d-e$ is made.

For $X\in S(N)$ of degree $r$ and $d+r<N$, define
\[
 Z=(-1)^{p(X)p(f)}fX+\frac{X(f)E}{N-d-r}\in S(N).
\]
Then
\[
 [X_\lambda e_f]
 =e_{X(f)}+\frac r{N-d}\partial Z
       -\frac{N-d-r}{N-d}\lambda Z.
\]
Here $e_0=0$. If $d+r\ge N$, both $fX$ and $X(f)$ vanish;
in the equality case, $X(f)=0$ follows from the exterior integration
identity for divergence-zero fields. Finally
$[X_\lambda Y]=[X,Y]$ for $X,Y\in S(N)$. These formulas show that
$[A_N{}_{\lambda}A_N]\subseteq\partial A_N+\lambda A_N+A_N$.
\end{proof}

\begin{lemma}
\label{sp:lem:smalln4}
Assume $N=2$ and write
\[
 q_i=\partial_i,\qquad
 p_i=\xi_i-\partial(\xi_iE),\qquad
 J_{ij}=\xi_i\partial_j-\frac12\delta_{ij}E,
 \qquad
 \ell_*=-1+\frac12\partial E.
\]
Then $J_{11}+J_{22}=0$, and
\[
 \ell_*,\ J_{11},\ J_{12},\ J_{21},\ q_1,\ q_2,\ p_1,\ p_2
\]
form a $\C[\partial]$-basis of $S_2$. Their nonzero
brackets include
\begin{align*}
 [\ell_*{}_{\lambda}J_{ij}]&=(\partial+\lambda)J_{ij},&
 [\ell_*{}_{\lambda}q_i]&=\left(\partial+\frac32\lambda\right)q_i,\\
 [\ell_*{}_{\lambda}p_i]&=\left(\partial+\frac32\lambda\right)p_i,&
 [J_{ij}{}_{\lambda}J_{kl}]&=\delta_{jk}J_{il}-\delta_{il}J_{kj},\\
 [J_{ij}{}_{\lambda}q_k]&=-\delta_{ik}q_j+\frac12\delta_{ij}q_k,&
 [J_{ij}{}_{\lambda}p_k]&=\delta_{jk}p_i-\frac12\delta_{ij}p_k,\\
 [q_i{}_{\lambda}p_j]&=-\delta_{ij}\ell_*+(\partial+2\lambda)J_{ji},&
 [q_i{}_{\lambda}q_j]&=[p_i{}_{\lambda}p_j]=0.
\end{align*}
After rescaling the two odd doublets and choosing the standard basis of
the traceless $J_{ij}$, these are the constant small-$N=4$ relations in
\cite[Section~3.4]{KLP}, with the state denoted there by $L$ equal to
$\ell_*$.
\end{lemma}
\begin{proof}
The displayed generators form the $N=2$ case of the generating space in
Lemma~\ref{sp:lem:primary}. The formulas follow
directly from the $W_2$ brackets and conformal sesquilinearity. For the
mixed odd bracket, for example,
\[
 [\partial_i{}_{\lambda}\xi_j]=\delta_{ij}+\lambda\xi_j\partial_i,
 \qquad
 [\partial_i{}_{\lambda}(\xi_jE)]=\delta_{ij}E-\xi_j\partial_i,
\]
which gives the stated expression. Replacing $q_i$ by $\sqrt2q_i$ and
$p_i$ by $-\sqrt2p_i$ gives the normalization of
\cite[equations~(3.39)--(3.47)]{KLP}. The identification of the same
Virasoro state is also consistent with \cite[Proposition~6.4]{FKR}.
\end{proof}

\subsection{Diagonalizable automorphisms}
\begin{lemma}\label{sp:lem:coordinate-lift}
Let $g$ be a continuous coordinate automorphism of
$\C[[t]]\otimes\Lambda_N$ of the form
\[
 g(t)=t+F(\xi),\qquad g(\xi_i)=\phi(\xi_i),
\]
where $F$ is even with zero constant term and $\phi\in\Aut\Lambda_N$.
Write $g_*D=gDg^{-1}$ and $\phi_*X=\phi X\phi^{-1}$.
For $X\in\Der\Lambda_N$ and $f\in\Lambda_N$, let
$a=X+f\partial_t$, whose coefficients are independent of $t$, and define
\[
 a^g=\phi_*X+
 \bigl(\phi(f)-(\phi_*X)(F)\bigr)\partial_t.
\]
Identify $X+f\partial_t$ with the corresponding free generator $X+f$
of $W_N$. The formula
\begin{equation}\label{sp:eq:coordinate-lift}
 \widehat g(a)=\sum_{j\ge0}\frac{(-1)^j}{j!}
                  \partial^j(F^j a^g)
\end{equation}
extends to a conformal automorphism of $W_N$. If $g$ multiplies the
Berezin volume by a nonzero scalar, it restricts to a conformal
automorphism of $S_N$.
\end{lemma}
\begin{proof}
The sum is finite because $F$ is nilpotent. In the vector-field coefficient
realization $a[m]=t^ma$, the image of its $m$th mode is
\[
 \sum_j\binom mj t^{m-j}F^ja^g
       =(t+F)^ma^g=g_*(t^ma).
\]
Here $(\partial^jb)[m]=(-1)^j m(m-1)\cdots(m-j+1)b[m-j]$.
These identities also hold for every nonnegative mode in the
annihilation algebra. The coefficient formulas determine the conformal
brackets, so the Lie identity for $g_*$ proves that $\widehat g$ preserves
them. The inverse coordinate change gives the inverse conformal map.
If $g$ multiplies the volume by a scalar, it preserves divergence-zero
vector fields. Thus all nonnegative modes of
$\Div(\widehat g(a))$, for $a\in S_N$, vanish. The mode realization of
the free module $\C[\partial]\otimes\Lambda_N$ is faithful, so
$\Div(\widehat g(a))=0$. Apply the same argument to $g^{-1}$.
\end{proof}

\begin{theorem}\label{sp:thm:semisimple}
Every parity-preserving diagonalizable conformal automorphism of $S_N$,
$N\ge2$, is conjugate by a conformal automorphism to a linear change of
the odd coordinates. In particular, it is conjugate to an automorphism
\[
 \xi_i\longmapsto c_i\xi_i\qquad(c_i\in\C^\times).
\]
In these conjugated coordinates it fixes $\ell_*$ and preserves the
primary quadratic generating space $A_N$ of Lemma~\ref{sp:lem:primary}.
\end{theorem}
\begin{proof}
First assume $N>2$. By \cite[Example 3.5]{FK}, the annihilation algebra
of $S_N$ is the \emph{derived} formal special algebra
$\mathfrak s=S(1,N)'$, with distinguished derivation
$-\ad\partial_t$. The full divergence-zero algebra has the additional
one-dimensional direction $\C\nu\partial_t$, where
$\nu=\xi_1\cdots\xi_N$; it is not substituted for $\mathfrak s$ here.

A conformal automorphism induces a continuous automorphism of
$\mathfrak s$ commuting with this distinguished derivation. The geometric
automorphism theorem \cite[Theorem 4.5]{CantariniKacAut} realizes it by a
continuous coordinate change multiplying the Berezin volume by a
nonzero scalar. The centralizer of $\mathfrak s$ in
$\widehat{\mathfrak W}_N$ is zero:
commuting with $\partial_t,\partial_1,\ldots,\partial_N$ forces an even
centralizing vector field to be a scalar multiple of $\partial_t$, and
commuting with $t\partial_t+E/N$ makes that scalar zero. Therefore the
coordinate change commutes with $\partial_t$ itself. It has exactly the
form in Lemma~\ref{sp:lem:coordinate-lift}.

Let $G$ be the group of these coordinate changes with constant
Berezinian. It acts faithfully on the finite-dimensional vector space
$\C t\oplus\Lambda_N$, so is an algebraic group.
Let $\mathfrak m$ be the augmentation ideal of $\Lambda_N$.
The linear part on $\mathfrak m/\mathfrak m^2$ gives a split epimorphism
\[
 G\longrightarrow\GL_N(\C).
\]
Its kernel $U$ is unipotent: it acts as the identity on the associated
graded space for the exterior-degree filtration, assigning degree one
to $t$. Its remaining matrix entries strictly raise this degree.
The constant-Berezinian condition is a closed algebraic condition.
Consequently $G=U\rtimes\GL_N(\C)$.

Diagonalizability of the original conformal automorphism implies
diagonalizability of its annihilation action. Indeed, a finite set of
$\C[\partial]$-generators is contained in finitely many eigenspaces,
so one squarefree polynomial annihilates the automorphism on $S_N$,
and therefore on its annihilation algebra. The subspace of
$\mathfrak s$ consisting of vector fields of $t$-degree at most one is
a finite-dimensional faithful $G$-module: it is $G$-stable and contains
all constant partial derivatives. Thus the corresponding element of
$G$ is semisimple.

A semisimple element of $U\rtimes\GL_N$ is conjugate into the displayed
$\GL_N$. Choose a finite
$\GL_N$-stable filtration
\[
 U=U_1\supset U_2\supset\cdots\supset U_s=1
\]
whose quotients $V_j=U_j/U_{j+1}$ are vector groups. Suppose inductively
that the element has the form $u_jg_0$ with $u_j\in U_j$ and semisimple
linear part $g_0\in\GL_N$. In the quotient algebraic group
\[
 H_j=V_j\rtimes\GL_N
\]
write
\[
 \bar u_j=\bar u_j^0+\bar u_j^1,
 \qquad
 \bar u_j^0\in\ker(1-\operatorname{Ad}g_0),
 \quad
 \bar u_j^1\in\im(1-\operatorname{Ad}g_0).
\]
Conjugation by an element of $U_j$ removes $\bar u_j^1$. The remaining
$\bar u_j^0$ commutes with $g_0$ in $H_j$. Since the image of the
original element in $H_j$ is semisimple, its commuting unipotent Jordan
factor must be trivial; hence $\bar u_j^0=0$. The conjugated element
therefore lies in $U_{j+1}g_0$. Induction gives conjugacy into
$\GL_N$. Lemma~\ref{sp:lem:coordinate-lift} lifts every conjugating
element to the conformal algebra, and a further linear conjugation
diagonalizes the resulting matrix.

For $N=2$, let $\mu$ be the given automorphism.
Lemma~\ref{sp:lem:smalln4} identifies our presentation,
including its specified Virasoro state, with the constant small-$N=4$
conformal algebra. Extend $\mu$ to the differential coefficient ring
$\widehat S=\C[t^q\mid q\in\mathbb Q]$, with derivation $d/dt$.
By \cite[Proposition~3.63]{KLP}, this extension is represented by
$Y\in\mathrm{SL}_2(\widehat S)$ and
$\begin{pmatrix}c&d\\ e&f\end{pmatrix}\in\mathrm{SL}_2(\C)$.
The images of the odd generators lie both in the constant algebra and
in their free $\widehat S$-span, so their coefficients are constant.
The two projections in \cite[equations~(3.66)--(3.67)]{KLP} therefore show that
$c(Y^{-1})^T$ and $dY$ have constant entries. Since $cf-de=1$,
at least one of $c,d$ is nonzero, and hence
$Y\in\mathrm{SL}_2(\C)$. Thus $Y'=0$, and
\cite[equation~(3.64)]{KLP} gives $\mu(\ell_*)=\ell_*$.
Conversely, constant pairs preserve the constant algebra, so
\cite[Proposition~3.69 and equation~(3.70)]{KLP} yields
\[
 \Aut S_2\simeq
 (\mathrm{SL}_2(\C)\times\mathrm{SL}_2(\C))/\{\pm(I,I)\}.
\]
This group acts faithfully on the finite-dimensional standard primary
generating space, so diagonalizability of $\mu$ makes its group element
semisimple. Every semisimple element lifts to a semisimple element of
the product, since the covering has finite central kernel. Conjugating
both factors puts it in a maximal torus. The diagonal coordinate
subgroup acts by
\[
 p_i\longmapsto c_i p_i,\qquad
 q_i\longmapsto c_i^{-1}q_i\qquad(i=1,2),
 \qquad c_1,c_2\in\C^\times.
\]
It is a two-dimensional torus, and the displayed automorphism group
has rank two, so this torus is maximal. Conjugacy of maximal tori now
gives the required conjugation into the diagonal coordinate subgroup.

Linear coordinate changes fix $1$ and $E$, and thus fix $\ell_*$.
They preserve $A_N$ by Lemma~\ref{sp:lem:primary}.
\end{proof}

\subsection{Coefficient algebras}
Let
\[
 B_N=\C[t,t^{-1}]\otimes\Lambda(\eta_1,\ldots,\eta_N),
 \qquad \nu_\eta=\eta_1\cdots\eta_N.
\]
For $B\in\C$, define
\[
 \mathfrak S_B=
 \left\{D\in\Der B_N:
       \operatorname{div}D+B t^{-1}D(t)=0\right\},
 \qquad \mathfrak S'_B=[\mathfrak S_B,\mathfrak S_B].
\]
Equivalently this is the condition $\operatorname{div}(t^BD)=0$,
interpreted formally even when $t^B$ is not an element of $B_N$.
We write $S'_B(1,N)=\mathfrak S'_B$, as in
\cite[Appendix A, item~3, pp.~29--30]{Billig}.

\begin{theorem}\label{sp:thm:coefficients}
Suppose $\mu\in\Aut S_N$ is diagonalizable. After the conjugation of
Theorem~\ref{sp:thm:semisimple}, write
\[
 \mu(\xi_i)=e^{-2\pi i b_i}\xi_i,\qquad B=\sum_i b_i.
\]
Then there is an isomorphism of full Lie superalgebras
\[
 \Alg(S_N,\mu)\simeq\mathfrak S'_B.
\]
If $B\notin\Z$, $\mathfrak S'_B=\mathfrak S_B$. If $B\in\Z$,
\[
 \mathfrak S_B=\mathfrak S'_B\oplus
       \C t^{-B}\nu_\eta\partial_t
\]
as vector spaces. The algebras $\mathfrak S'_B$ are centerless. The extra direction
is not a coefficient mode of $S_N$ and is not a central quotient of
that coefficient algebra.
\end{theorem}
\begin{proof}
For the diagonal linear automorphism, its extension to $W_N$ has the
usual twisted vector-field realization. In its coefficient ring set
\[
 \eta_i=t^{b_i}\xi_i.
\]
These variables and $t^{\pm1}$ generate the ring $B_N$. A monomial
coefficient of twisted degree $b_I+\Z$ becomes an ordinary Laurent
coefficient after this substitution. Likewise, the allowed modes of
$\xi_I\partial_j$ have degree $b_I-b_j+\Z$. The resulting change of
coordinates identifies the full twisted coefficient algebra of $W_N$
with $\Der B_N$. Its Berezinian transforms the original volume into
\[
 t^B\,dt\,d\eta_1\cdots d\eta_N.
\]
Thus the divergence-zero condition becomes the defining equation of
$\mathfrak S_B$.

Let $\nu=\xi_1\cdots\xi_N$, and let $\Lambda_N^{<N}$ be the span of
the exterior monomials of degree less than $N$.
The image of conformal divergence is the free module
\[
 I=\Div W_N
   =\bigl(\C[\partial]\otimes\Lambda_N^{<N}\bigr)
        \oplus\C[\partial]c,
 \qquad c=\partial\nu.
\]
The sequence $0\to S_N\to W_N\to I\to0$ is equivariant.
On the twisted scalar extension of the free module $I$, the total
translation $\partial+d/dt$ is injective: the highest $\partial$-degree
in a finite expression cannot cancel. The exact sequence of
coinvariants therefore embeds $\Alg(S_N,\mu)$ into
$\Alg(W_N,\mu)$, with image the kernel of the induced map to
$\Alg(I,\mu)$.

The natural map from the latter coefficient module to twisted
functions sends
\[
 f[m]\longmapsto t^mf\quad(\deg f<N),\qquad
 c[m]\longmapsto -m t^{m-1}\nu,
       \quad m\in B+\Z.
\]
For $B\notin\Z$, this map is injective, so the coefficient image is all
of $\mathfrak S_B$. For $B\in\Z$, its kernel is exactly $\C c[0]$.
The element $\nu\partial_t$ of the original twisted vector-field
algebra maps to $-c[0]$, and after the coordinate substitution it is
$t^{-B}\nu_\eta\partial_t$. The coefficient image is consequently a
codimension-one Lie subsuperalgebra of $\mathfrak S_B$ complementary to
that line.

Finally, the coefficient algebra of a perfect conformal algebra under
a diagonalizable twist is perfect. For homogeneous $a,b$, hold the
sum of their mode indices fixed and vary the first index in its
integer coset. The mode commutator is a finite polynomial in that
index in the binomial basis. Inverting a finite Vandermonde matrix
places every coefficient of every product $a_{(j)}b$ in the ordinary
derived algebra. The conformal derived algebra spans $S_N$, so this
proves the assertion. Therefore the image above is perfect.
A perfect codimension-one Lie subsuperalgebra is an ideal: its induced
action on the one-dimensional quotient superspace has abelian image,
and thus vanishes on its derived algebra, which is the subalgebra
itself. The quotient of $\mathfrak S_B$ by the coefficient image is
therefore zero or one-dimensional and abelian. Perfectness gives the
reverse inclusion in the derived algebra, so the coefficient image is
exactly $\mathfrak S'_B$.

For centerlessness, the image contains every
$t^k\partial_{\eta_i}$ and every $\eta_i\partial_{\eta_j}$ with
$i\ne j$. A central vector field commuting with all constant odd
partials has coefficients depending only on $t$. Its brackets with
$t\partial_{\eta_i}$ force its $\partial_t$ coefficient to vanish;
its brackets with the off-diagonal linear fields force all remaining
coefficients to vanish. The latter fields are available because
$N\ge2$. This excludes both an even and an odd center.
\end{proof}

\begin{corollary}\label{sp:cor:cyclic}
For a cycle of length $r$ with diagonalizable return automorphism $\mu$
on $S_N$, the full cyclic coefficient algebra has an isomorphism
\[
 \Alg_r(S_N,\mu)\simeq S'_{B+r-1}(1,N).
\]
Here $\mu$ is first conjugated as in
Theorem~\ref{sp:thm:semisimple}, and $B$ is the sum of chosen
logarithms of its linear-coordinate eigenvalues.
\end{corollary}
\begin{proof}
Apply Theorem~\ref{cy:thm:normalization} to the invariant presentation
of Lemma~\ref{sp:lem:primary}. On that presentation
$P=\ell_{*(1)}-I$, and
$\kappa_r=e^{2\pi i(r-1)P}$ is the scalar linear change
\[
 \xi_i\longmapsto e^{-2\pi i(r-1)/N}\xi_i.
\]
It adds $(r-1)/N$ to every $b_i$, hence adds $r-1$ to $B$.
Theorem~\ref{sp:thm:coefficients} gives the conclusion.
\end{proof}

\subsection{Parameter changes and the \texorpdfstring{$\widetilde S$}{tilde S}-realization}
\begin{proposition}
\label{sp:prop:parameters}
For every $B\in\C$ and $k\in\Z$,
\[
 S'_B(1,N)\simeq S'_{B+k}(1,N),\qquad
 S'_B(1,N)\simeq S'_{-B-2}(1,N).
\]
In particular $B'\in B+\Z$ or $B'\in-B+\Z$ is sufficient for an
isomorphism of the underlying Lie superalgebras.
\end{proposition}
\begin{proof}
In the density $t^Bdt\,d\eta_1\cdots d\eta_N$, make the Laurent
coordinate change $\eta'_1=t^k\eta_1$, leaving the other coordinates
fixed. The odd Berezinian gives the density
$t^{B+k}dt\,d\eta'_1\,d\eta_2\cdots d\eta_N$.
The resulting coordinate isomorphism carries the full divergence
stabilizers to one another and therefore carries their derived
algebras to one another. For the second assertion use $t=s^{-1}$
and leave the odd variables fixed. The density becomes
$-s^{-B-2}ds\,d\eta_1\cdots d\eta_N$; its nonzero constant factor
does not affect the stabilizer. This proves the two isomorphisms.
The calculations are identities of logarithmic density derivatives,
so they also apply when the complex powers are only formal symbols.
\end{proof}

\begin{remark}
Proposition~\ref{sp:prop:parameters} gives sufficient parameter
identifications. The case $N=2$, including these transformations,
is discussed in \cite[\S1.1 and \S\S1.7--1.8]{GLS}. The parameter
$\lambda$ there is the present $B$ for the full density stabilizer
$\mathfrak S_B$. By Theorem~\ref{sp:thm:coefficients}, this equals
$S'_B(1,N)$ for $B\notin\Z$; at integral parameters one passes to
the derived ideal. A necessary-and-sufficient isomorphism criterion
for arbitrary $N$ is not proved here.
\end{remark}

\subsubsection*{The \texorpdfstring{$\widetilde S$}{tilde S}-realization}
The conformal algebra $\widetilde S_N$ is not needed in the reduction of
Theorem~\ref{thm:ordinary-model}, but it appears in the finite conformal
classification and in modern statements of the conjectural list. The
following observation explains why it contributes no additional abstract
simple Lie-superalgebra type.

\begin{proposition}
\label{sp:prop:tilde-derived}
Let $N\ge2$ be even, set
\[
 \nu=\xi_1\cdots\xi_N,\qquad F=1+\nu,
\]
and define the Laurent density stabilizers
\[
 \widetilde{\mathfrak S}
   =\{D\in W(1,N):\diver(FD)=0\},\qquad
 \mathfrak S_0=\{D\in W(1,N):\diver D=0\}.
\]
Then $\widetilde{\mathfrak S}\simeq\mathfrak S_0$ and
\[
 [\widetilde{\mathfrak S},\widetilde{\mathfrak S}]
   \simeq [\mathfrak S_0,\mathfrak S_0]
   =S'_0(1,N).
\]
Moreover, there is a vector-space decomposition
\[
 \widetilde{\mathfrak S}
   =[\widetilde{\mathfrak S},\widetilde{\mathfrak S}]
     \oplus\C\nu\partial_t.
\]
Thus the quotient by the derived algebra is one-dimensional,
even and abelian.
\end{proposition}
\begin{proof}
The Laurent coordinate change
\[
 t'=t(1+\nu),\qquad \xi'_i=\xi_i
\]
has inverse $t=t'(1-\nu)$ and Berezinian $F$. Hence it carries
$F\,dt\,d\xi_1\cdots d\xi_N$ to the standard volume form in the
primed coordinates. Conjugation of vector fields gives
$\widetilde{\mathfrak S}\simeq\mathfrak S_0$ and sends
\[
 \nu\partial_t\longmapsto\nu F\partial_{t'}
 =\nu\partial_{t'}.
\]
The assertions about the derived algebra and its complement now
follow from Theorem~\ref{sp:thm:coefficients} with $B=0$.
\end{proof}

The definition in \cite[Appendix A, item~4, pp.~30--31]{Billig}
gives the full stabilizer $\widetilde{\mathfrak S}$. In this
realization the simple algebra is its derived algebra
$S'_0(1,N)$; the full stabilizer has the one-dimensional quotient
above. The finite conformal algebra $\widetilde S_N$ has
annihilation algebra $S(1,N)'$ with a different distinguished
translation \cite[Example~3.7]{FK}. The isomorphism of derived
Laurent algebras therefore does not identify the conformal
presentations or their distinguished translations.

\section{The contact families}\label{sec:contact}
\subsection{Diagonalizable automorphisms}
\label{ct:sec:arbitrary}

Let $\widehat{\mathfrak K}_N$ denote the formal contact algebra over
$\C[[t]]\otimes\Lambda(\xi_1,\ldots,\xi_N)$.
By \cite[Theorem 4.5]{CantariniKacAut}, every continuous automorphism
of $\widehat{\mathfrak K}_N$ is induced by a formal coordinate change
preserving the contact form up to an invertible factor; none of the
exceptions in that theorem applies here.
The algebra $K_0$ is the Virasoro conformal algebra. For $N\ge1$,
\cite[Example 3.8]{FK} identifies the annihilation algebra of $K_N$
with $\widehat{\mathfrak K}_N$ and its distinguished translation
with $-\partial_t$. This includes $K_4$, although $K_4$ itself
is not simple as a conformal algebra.

Write $B=\Lambda(\xi_1,\ldots,\xi_N)$ and equip it with the Poisson
bracket $\beta(\xi_i,\xi_j)=\delta_{ij}$. Let $d$ be the exterior
differential on super differential forms and set
$\alpha=\sum_i\xi_i\,d\xi_i$. A Poisson automorphism of $B$ means an
even associative-algebra automorphism preserving $\beta$;
equivalently, it preserves $d\alpha$.
In the contact-function realization the
form is $\omega=dt-c\alpha$, where $c\ne0$ is a fixed scalar determined
by the normalization of the Poisson bracket. A fixed rescaling of the
odd variables changes $c$ and has no effect on the argument below.
A change of coordinates is called strict contact if it preserves $\omega$.

In the contact convention of \cite[Appendix A]{Billig},
\[
 A=\Lambda(\xi_1,\ldots,\xi_N),\qquad
 P=I-\tfrac12 E,\qquad \ell=1,
\]
where $E$ is exterior degree. The bracket is \eqref{cy:eq:qp} with the nondegenerate Poisson bracket $\beta$.

For $N=4$, write $\nu=\xi_1\xi_2\xi_3\xi_4$.

\begin{lemma}\label{ct:lem:lift}
Every conformal automorphism of $K_N$ is induced on its annihilation
algebra by a unique strict-contact change
\begin{equation}\label{ct:eq:coordinates}
 t\longmapsto t+f(\xi),\qquad
 \xi\longmapsto\varphi(\xi),
\end{equation}
where $\varphi$ is an even Poisson automorphism of $B$ and
$f\in\bigoplus_{d\ge4}\Lambda^d(\C^N)$ is even. Conversely every Poisson
automorphism $\varphi$ has a unique such lift. Its conformal action on
the free generating space $B$ is
\begin{equation}\label{ct:eq:conformal-lift}
 a\longmapsto\sum_{j\ge0}\frac{(-\partial)^j}{j!}
                         \bigl(f^j\varphi(a)\bigr).
\end{equation}
The sum is finite, and the inverse has the same form for the inverse
coordinate change.
\end{lemma}
\begin{proof}
A conformal automorphism induces a continuous automorphism of the
annihilation algebra commuting with $\operatorname{ad}\partial_t$.
Since $\widehat{\mathfrak K}_N$ is centerless, it fixes $\partial_t$ itself. Apply the
geometric theorem. Commutation with $\partial_t$ forces the coordinate
change to have the form \eqref{ct:eq:coordinates}. Its contact multiplier
is one, by comparison of the $dt$ coefficient. Therefore
\begin{equation}\label{ct:eq:primitive}
 df=c\bigl(\varphi(\alpha)-\alpha\bigr),\qquad
 \varphi(d\alpha)=d\alpha.
\end{equation}
The second identity says precisely that $\varphi$ is Poisson.

Conversely, for any Poisson $\varphi$, the one-form
$c(\varphi(\alpha)-\alpha)$ is closed. The linear part of $\varphi$ is
orthogonal, so its degree-two part is zero. Euler contraction gives a
primitive $f$ with zero constant term; it has exterior degree at least
four. More explicitly, on a form homogeneous of total degree $q>0$,
counting both $\xi_i$ and $d\xi_i$ with degree one, the homotopy is
$q^{-1}\iota_E$, where $E=\sum_i\xi_i\partial_i$.
The primitive is unique because the kernel of $d$ on $B$ consists of
constants. A nonzero constant translation of $t$ is excluded by
continuity on $B[[t]]$. Thus the strict-contact lift exists and is
unique. Uniqueness also makes lifting compatible with composition.

The induced conformal maps are finite. On contact functions,
\[
 t^ma\longmapsto(t+f)^m\varphi(a)
   =\sum_{j\ge0}\binom mj t^{m-j}f^j\varphi(a).
\]
Nilpotence of $f$ makes this a finite identity, including for formal
power series by continuity. The mode relation
$(\partial^jb)[m]=(-1)^jm(m-1)\cdots(m-j+1)b[m-j]$
now gives \eqref{ct:eq:conformal-lift}. The contact change preserves
the mode brackets, hence the conformal products. To see the last
implication directly, apply the mode commutator formula for all
sufficiently large nonnegative indices and compare its finitely many
polynomial coefficients. A state all of whose sufficiently large modes
vanish is zero in the free contact presentation. The inverse coordinate
change gives the inverse conformal map by the same argument.
\end{proof}

\begin{lemma}
\label{ct:lem:darboux}
Every diagonalizable even Poisson automorphism $\varphi$ of $B$ is
conjugate, by a Poisson automorphism with identity linear part, to its
orthogonal linear part $g$.
\end{lemma}
\begin{proof}
Let $\mathfrak m$ be the augmentation ideal. Since $\varphi$ is
diagonalizable on the finite-dimensional space $B$, choose an equivariant
section of
$\mathfrak m_{\overline1}\longrightarrow
 (\mathfrak m/\mathfrak m^2)_{\overline1}$.
The resulting odd coordinate map $\theta$ has identity linear part and
satisfies $\varphi\theta=\theta g$. It is an algebra automorphism,
because it induces the identity on the associated exterior graded
algebra. It need not be Poisson.

Set $\Omega_0=d\alpha$ and $\Omega=\theta^{-1}\Omega_0$.
Then $g\Omega=\Omega$, and $\Omega-\Omega_0$ is a closed two-form
whose coefficients have exterior degree at least two. Euler homotopy
gives a $g$-invariant one-form $\eta$ with
$d\eta=\Omega-\Omega_0$, whose coefficients have degree at least three.
For $s$ an ordinary auxiliary parameter, the form
\[
 \Omega_s=\Omega_0+s(\Omega-\Omega_0)
\]
is nondegenerate. Its inverse is polynomial in $s$ and the exterior
variables, because its difference from the invertible constant matrix
has nilpotent coefficients. Solve
$\iota_{V_s}\Omega_s=-\eta$.
The resulting even vector field is $g$-invariant and raises exterior
degree by at least two. Its time-dependent flow is a finite sum:
every sufficiently long composition of the $V_s$ vanishes on $B$.
The coefficients are obtained by integrating polynomials in $s$.
Cartan's formula proves that the flow pulls $\Omega_s$ back to
$\Omega_0$. Taking the inverse of its value at $s=1$ gives an
automorphism $\kappa$ commuting with $g$ and satisfying
$\kappa\Omega_0=\Omega$.

Thus $\theta\kappa$ preserves $\Omega_0$ and satisfies
$\varphi\theta\kappa=\theta\kappa g$, as required.
\end{proof}

\begin{theorem}\label{ct:thm:normalization}
Every diagonalizable automorphism of $K_N$ is conjugate by a conformal
automorphism to an element of $\OO_N(\C)$ acting linearly on the odd
generators. In these conjugated coordinates it fixes the standard contact
Virasoro state and preserves the primary quadratic generating space
$\Lambda(\C^N)$.
The same conclusion holds for every diagonalizable automorphism of
$K'_4$; in these conjugated coordinates it fixes the standard contact
Virasoro state and preserves the primary generating space
$\Lambda^{\le3}(\C^4)\oplus\C\partial\nu$.
\end{theorem}
\begin{proof}
For $K_N$, the action modulo $\partial K_N$ is the Poisson
automorphism $\varphi$ in Lemma~\ref{ct:lem:lift}. It is diagonalizable,
being a quotient of a diagonalizable operator. Lemma~\ref{ct:lem:darboux}
and uniqueness of the strict-contact lift conjugate its coordinate
action to $(t,\xi)\mapsto(t,g\xi)$. Formula~\eqref{ct:eq:conformal-lift}
lifts this conjugation to the conformal algebra. The resulting linear
orthogonal automorphism fixes $1$ and preserves $\Lambda(\C^N)$.

For $K'_4$, use the identification $K_4\simeq\Cder K'_4$
established in the proof of Proposition~\ref{prop:inner-S0-K4}.
Conjugation on conformal derivations extends every automorphism of
$K'_4$ uniquely to $K_4$; under the identification, the inner conformal
derivations are exactly $K'_4$. Conversely, $K'_4$ is the derived
conformal algebra of $K_4$ and is therefore characteristic. Restriction
and extension are inverse, because $K_4$ acts faithfully on $K'_4$.

The extension of a diagonalizable automorphism is diagonalizable.
Indeed, a diagonalizable $\C[\partial]$-module automorphism of a finite
module has only finitely many eigenvalues: decompose a finite generating
set into eigenvectors and use commutation with $\partial$.
Conjugation on
$\operatorname{Hom}_{\C}(K'_4,\C[\lambda]\otimes K'_4)$
therefore decomposes into the finitely many eigenvalue ratios. Its
invariant subspace of conformal derivations is diagonalizable.
Apply the result for $K_4$ and restrict the conjugating automorphism to
its characteristic derived algebra. The resulting orthogonal action
preserves the stated primary generating space.
\end{proof}

\subsection{Primary presentations and the cyclic phase}
There may be generators of conformal weight zero, namely exterior degree four. Lemma~\ref{cy:lem:prim} nevertheless limits possible derivative mixing to
\[
 \Lambda^2\longrightarrow\Lambda^2+\partial\Lambda^4.
\]
An automorphism fixing $1$ sends each $\xi_i$ to a linear combination of the $\xi_j$, since these are exactly the primary states of weight $3/2$. But
\[
 (\xi_i)_{(1)}\xi_j=\xi_i\xi_j.
\]
Therefore every quadratic generator maps into $\Lambda^2$, with no
$\partial\Lambda^4$ component. All other exterior degrees are preserved
by Lemma~\ref{cy:lem:prim}. Thus $A$ is invariant.

For a degree-$d$ generator, $p=1-d/2$, and hence
\[
 \kappa_r(a)=e^{2\pi i(r-1)(1-d/2)}a=(-1)^{(r-1)d}a.
\]
Let $\varepsilon$ be the parity automorphism,
$\varepsilon(a)=(-1)^{p(a)}a$ for homogeneous $a$. Then
\begin{equation}\label{cy:eq:parity}
 \tau=\mu\varepsilon^{r-1}.
\end{equation}
Thus the cyclic phase is $\varepsilon^{r-1}$. For the analogous distinction in permutation-twisted vertex-superalgebra modules, see \cite{BVW}.

\subsubsection*{The derived algebra \texorpdfstring{$K'_4$}{K4-prime}}
Let $\nu=\xi_1\xi_2\xi_3\xi_4$ and $c_4=\partial\nu$. The derived conformal algebra has a free generating space
\[
 A'=\bigoplus_{d=0}^3\Lambda^d\ \oplus\C c_4;
\]
see \cite[Example 3.8]{FKR}. Its generators have weights $2-d/2$ on $\Lambda^d$ and weight $1$ on $c_4$. They are primary for $\ell=1$. None has weight zero, so every automorphism fixing this $\ell$ preserves $A'$.

Its brackets in this basis are first order in $\partial,\lambda$. Indeed, a top-degree exterior term from the product term of \eqref{cy:eq:qp} can only contribute a multiple of $c_4$; its $\lambda$ coefficient vanishes because the $P$-weights of the two generators sum to zero. A top-degree term from the Poisson bracket of two generators of degree at most three is zero: the only possible degrees are three and three, and two triples in four variables share at least two indices, leaving a repeated odd variable after differentiation. Brackets with $c_4$ follow by sesquilinearity. For example,
\[
 [(c_4)_\lambda 1]=\lambda c_4,\qquad
 [\xi_i{}_{\lambda}c_4]=(\partial+\lambda)\beta(\xi_i,\nu),\qquad
 [(c_4)_\lambda c_4]=0.
\]
Thus Theorem~\ref{cy:thm:normalization} applies directly to $K'_4$, with $p_{c_4}=0$:
\[
 f\nn{m}\longmapsto r^{-(1-d/2)}f[m-(r-1)(1-d/2)],\qquad
 c_4\nn{m}\longmapsto c_4[m].
\]
This calculation in the free basis of $K'_4$ gives the automorphism
change \eqref{cy:eq:parity}. In particular, the cyclic map fixes
$c_4[0]$ whenever that mode is allowed.

\subsection{Contact coefficient algebras}
\label{ct:sec:coefficients}

Let $\rho$ denote an orthogonal reflection in a nonisotropic line when
$N\ge1$.

\begin{theorem}\label{ct:thm:twists}
Let $M=K_N$ with $N\ne4$, or let $M=K'_4$, and let $\mu$ be an
arbitrary diagonalizable automorphism of $M$. After the conformal
conjugacy in Theorem~\ref{ct:thm:normalization}, write $g\in\OO_N$
for its orthogonal representative. Then
\[
 \Alg(M,\mu)\simeq
 \begin{cases}
 \Alg(M,\id),&\det g=1,\\
 \Alg(M,\rho),&\det g=-1.
 \end{cases}
\]
These are isomorphisms of the full coefficient algebras. For a cyclic
orbit of length $r$, replace $\det g$ by
\begin{equation}\label{ct:eq:cyclic-sign}
 \det g\,(-1)^{N(r-1)}.
\end{equation}
For $N=0$ only the identity case occurs. For $K'_4$ the cyclic factor
in \eqref{ct:eq:cyclic-sign} is always one.
\end{theorem}
\begin{proof}
Conformal conjugacy induces a coefficient isomorphism. If $\det g=1$,
choose a semisimple $X\in\mathfrak{so}_N$ with $e^{2\pi iX}=g$
by taking opposite logarithms on reciprocal eigenvalue pairs. The
$-1$ eigenspace has even dimension, so it can also be paired.
Let $h$ be the quadratic current whose zeroth product acts as $X$
on the odd linear generators. Then $h$ has weight one,
$[h_\lambda h]=0$, and $g$ fixes $h$.
Lemma~\ref{ck:lem:flow} removes the twist.

If $\det g=-1$, its $-1$ eigenspace is nondegenerate and nonzero.
All reflections in nonisotropic lines are orthogonally conjugate,
so we may choose $\rho$ to reflect a line in this eigenspace.
Then $g\rho\in\SO_N$ commutes with $\rho$. Choose semisimple
$X\in\mathfrak{so}_N$ with $e^{2\pi iX}=g\rho$, commuting with
$g$ and $\rho$ and vanishing on the reflected line. The corresponding
quadratic current changes the twist from $g$ to $\rho$ by
Lemma~\ref{ck:lem:flow}.
The cyclic coefficient theorem applies to the
invariant presentation just obtained. Its phase is parity to the power
$r-1$, induced by $(-I_N)^{r-1}$, proving
\eqref{ct:eq:cyclic-sign}.
\end{proof}

For $K_N$ with $N\ne4$, the two representatives have an explicit
contact interpretation. The identity gives the usual contact algebra
over $\C[t,t^{-1}]\otimes\Lambda(\xi_1,\ldots,\xi_N)$.
For $K'_4$, the derived image and the possible central kernel are
described in Proposition~\ref{ct:prop:center}.
For a reflection acting by $-1$ on $\xi_N$, use the algebra of functions
spanned by $t^m a$, where $m$ is integral for monomials not containing
$\xi_N$ and half-integral for monomials containing $\xi_N$.
This is the ordinary order-two contact twist. With
$\eta=t^{1/2}\xi_N$, it is a contact algebra on the Laurent
superalgebra with form
\[
 dt-c\sum_{i<N}\xi_i\,d\xi_i-c\,t^{-1}\eta\,d\eta.
\]
The term involving $dt$ in $\xi_Nd\xi_N$ vanishes because $\eta^2=0$.
The further substitution $\eta=t\theta$ gives
$t^{-1}\eta\,d\eta=t\theta\,d\theta$, since $\theta^2=0$.
This is the contact form with one odd exponent equal to one in
\cite[Section~1]{KacVanDeLeur}. Thus both representatives are explicit
contact forms on the punctured superline.

\begin{proposition}
\label{ct:prop:center}
For $N\ne4$ and every diagonalizable $g\in\OO_N(\C)$,
\[
 Z(\Alg(K_N,g))=0.
\]
For a diagonalizable $g\in\OO_4(\C)$, set
$\nu=\xi_1\xi_2\xi_3\xi_4$ and $c_4=\partial\nu$. Then
\[
 Z(\Alg(K'_4,g))=
 \begin{cases}
 \C c_4[0],&\det g=1,\\
 0,&\det g=-1.
 \end{cases}
\]
The conjugacy, current transformations, and cyclic normalization above
preserve this central coefficient whenever it is present.
\end{proposition}
\begin{proof}
By Theorem~\ref{ct:thm:twists}, it suffices to take $g=\id$ or
$\rho$, with only the identity case when $N=0$. For $K_N$, $N\ne4$,
the free coefficient basis identifies the algebra with the corresponding
contact-function realization. For $K'_4$, the natural map to the
coefficient algebra of $K_4$ sends
\[
 c_4[m]\longmapsto-m\nu[m-1].
\]
On all other free generators it is the identity. The free-module
coefficient bases show that its kernel is exactly $\C c_4[0]$ when
$m=0$ is allowed, and zero otherwise. Since $g\nu=(\det g)\nu$,
zero is allowed exactly in the determinant-one case. In the identity
case the image is the derived contact algebra, omitting precisely the
top coefficient $\nu[-1]$; compare \cite[Example 3.9]{FK}.
In the reflection case all allowed top coefficients occur, so the image
is the full reflection twisted contact algebra.
The possible kernel is central: sesquilinearity in the ambient $K_4$
gives $(c_4)_{(0)}b=(\partial\nu)_{(0)}b=0$ for every $b\in K'_4$,
and the coefficient commutator with $c_4[0]$ has only this zeroth term.

All these contact-function algebras are centerless. A
central function commutes with $\partial_t$, hence is independent of
$t$. Commuting with the contact field of $t$ then restricts it to
exterior degree two. Its brackets with the allowed linear odd modes
force all its odd partial derivatives to vanish, so it is zero.
These fields also belong to the derived contact image. This proves
the zero-center assertion for $K_N$, $N\ne4$, and shows that the stated
kernel is the entire center for $K'_4$.

Every strict-contact change sends the top monomial to its linear
determinant times itself: higher exterior terms vanish and $f\nu=0$.
The conjugator in Lemma~\ref{ct:lem:darboux} has identity linear part,
so it fixes $\nu$ and $c_4$. Quadratic currents commute conformally with
$\nu$, and hence with $c_4$, so the current transformations fix
$c_4[0]$. The cyclic formula also sends $c_4[0]$ to $c_4[0]$.
\end{proof}

\begin{remark}
The two parity representatives are stated to be nonisomorphic in
\cite[Sections~1 and~4]{KacVanDeLeur}, while the isomorphisms within
each parity class are given in \cite[Theorem~4.1(a)]{KacVanDeLeur}.
Our exhaustion theorem uses only the reduction to these representatives.
In rank four, the even-parity representative is the derived contact
algebra, as described in Proposition~\ref{ct:prop:center}.
\end{remark}

\section{The Cheng--Kac family}\label{sec:chengkac}
Let $V=\C^6$ with its standard nondegenerate symmetric form and
orthonormal basis $\xi_1,\ldots,\xi_6$. Choose an orientation and denote
the two Hodge eigenspaces in $\Lambda^3V$ by $\Lambda^3_+V$ and
$\Lambda^3_-V$. Let $\operatorname{pr}_+$ be the projection onto
$\Lambda^3_+V$ along $\Lambda^3_-V$. Set
\[
 A=\Lambda^0V\oplus\Lambda^1V\oplus\Lambda^2V\oplus\Lambda^3_+V,
 \qquad C=\C[\partial]\otimes A.
\]
This is the primary model in
\cite[Theorem 6.10(c), Appendix A, item 5]{Billig}.
Let $\beta$ be the inherited Poisson bracket with
$\beta(\xi_i,\xi_j)=\delta_{ij}$. Let $\pi_A:\Lambda V\to A$ be the projection with kernel
$\Lambda^3_-V+\Lambda^{\ge4}V$, and let $T_A$ act by $2$ on
$\Lambda^3_+V$ and by $1$ in the other degrees. Define
\[
 a\cdot b=T_A\pi_A\bigl(T_A^{-1}a\wedge T_A^{-1}b\bigr).
\]
Let $E$ be the exterior-degree operator, so $Ea=da$ for
$a\in\Lambda^dV$. Set
\[
 P=I-\tfrac12 E,\qquad
 [a_\lambda b]=\partial(P(a)\cdot b)
   +\lambda(P(a)\cdot b+a\cdot P(b))+\beta(a,b).
\]
Its standard Virasoro state is $\ell=1$. A generator of exterior degree
$d$ has weight $2-d/2$; the four weights are $2,3/2,1,1/2$.

\begin{theorem}\label{ck:thm:autck}
Every conformal automorphism of $C=CK_6$ fixes $1$, preserves the primary
generating space $A$, and is induced by a unique element of $\SO(V)$.
Consequently $\Aut C\simeq\SO_6(\C)$.
\end{theorem}
\begin{proof}
By Lemma~\ref{ck:lem:currentcore}, every conformal automorphism fixes $1$.
The primary generator
weights are $2,3/2,1,1/2$, all nonzero.
Lemma~\ref{cy:lem:prim} therefore shows that every automorphism
preserves $A$ and its weight spaces. Since these weights distinguish
the exterior degrees, each exterior-degree subspace is invariant.

Let $g$ be its restriction to $V$. For $v\in V$,
\[
 v_{(0)}v=(v,v)1,
\]
since $v\cdot v=0$. Preservation of this bracket gives $g\in\OO(V)$.
The identity $v_{(1)}w=v\wedge w$ determines the action on $\Lambda^2V$
as $\Lambda^2g$. For $v\in V$, $b\in\Lambda^2V$,
\[
 v_{(1)}b=\operatorname{pr}_+(v\wedge b).
\]
These products span $\Lambda^3_+V$. Therefore their relations require
$\Lambda^3g$ to preserve the kernel $\Lambda^3_-V$ of the projection.
An orthogonal map of determinant $-1$ interchanges the two Hodge
eigenspaces, while one of determinant $1$ preserves them. Hence
$\det g=1$.

Conversely, every element of $\SO(V)$ preserves $P$, $\beta$, the two
Hodge eigenspaces, and the rescaled truncated product. It therefore defines
a conformal automorphism. The displayed products also prove uniqueness.
\end{proof}

\begin{theorem}\label{ck:thm:cktwist}
For every diagonalizable conformal automorphism $\mu$ of $CK_6$,
\[
 \Alg(CK_6,\mu)\simeq\Alg(CK_6,\id)\simeq CK(6).
\]
For the associated cyclic automorphism $\sigma$ of $CK_6^{\oplus r}$,
one also has
\[
 \Alg(CK_6^{\oplus r},\sigma)\simeq CK(6).
\]
All these coefficient algebras have zero center.
\end{theorem}
\begin{proof}
By Theorem~\ref{ck:thm:autck}, $\mu$ is induced by a semisimple
$g\in\SO_6(\C)$. Choose a hyperbolic basis in which $g$ has reciprocal
eigenvalue pairs $\lambda_i,\lambda_i^{-1}$. The eigenvalue $-1$ has even
multiplicity, so it can also be paired. Choose $q_i\in\C$ with
$e^{2\pi i q_i}=\lambda_i$. The operator with eigenvalues
$q_i,-q_i$ on each pair lies in $\mathfrak{so}_6$ and exponentiates to $g$.

The current space $\Lambda^2V$ acts on $V$ as $\mathfrak{so}_6$ by
$h_{(0)}v=\beta(h,v)$. Choose $h$ representing that semisimple operator.
It is primary of weight one and $[h_\lambda h]=0$. The products generating
$A$ show that
\[
 \mu=e^{2\pi i h_{(0)}}
\]
on all of $C$. Lemma~\ref{ck:lem:flow} with $\sigma=\mu$ has $\tau=\id$.

To relate this ordinary coefficient algebra to the usual superconformal
notation, \cite[Theorem 6.10(c), Appendix A, item 5]{Billig} identifies
$CK(6)$ with $\Alg(CK_6,\varepsilon)$, where
$\varepsilon=e^{2\pi iP}$ is parity in the present normalization.
The element $-I_V$ lies in $\SO_6$, so the preceding argument also removes
$\varepsilon$. Thus $\Alg(CK_6,\id)\simeq CK(6)$. The latter is simple;
all centers in these identified coefficient algebras are zero.

For a cycle of length $r$, the phase $\kappa_r=\varepsilon^{r-1}$ is
induced by $(-I_V)^{r-1}\in\SO_6$ and commutes with every automorphism.
Corollary~\ref{cy:cor:cyclic} reduces the cyclic coefficient algebra to
another diagonalizable twist of $CK_6$, so the first assertion applies.
\end{proof}

On a free primary generator $a$ with $h_{(0)}a=qa$, the
coefficient isomorphism takes the form
\begin{equation}\label{ck:eq:explicitflow}
 a[m]\longmapsto a[m+q]+(h_{(1)}a)[m+q-1].
\end{equation}
Indeed, $h_{(j)}a=0$ for $j\ge2$ and $(h_{(1)})^2a=0$ on $A$.
The latter follows by exterior degrees: $h_{(1)}1=h$,
$h_{(1)}V\subseteq\Lambda^3_+V$, and $h_{(1)}$ kills
$\Lambda^2V+\Lambda^3_+V$.
Extend this formula to derivative states by the coefficient relations.
For $L_n=1[n+1]$, the map gives $L_n\mapsto L_n+h[n]$.

\section{Completion of the proof and polynomial growth}
\label{sec:final-reduction}

\begin{proof}[Proof of Theorem~\ref{thm:synthesis}]
By Theorem~\ref{thm:ordinary-model},
\[
 L\simeq\mathcal B/Z(\mathcal B),\qquad
 \mathcal B=\Alg(M^{\oplus r},\sigma),
\]
where $M$ is a noncurrent finite simple conformal superalgebra and
$\sigma$ is diagonalizable and permutes its $r$ copies transitively.
Its return automorphism $\mu$ is diagonalizable.
Identifying successive constituents by the component maps of $\sigma$
puts the cycle in the form used in Section~\ref{sec:cyclic}.
A conformal conjugacy of $\mu$ extends componentwise to this cycle.

The family results give the coefficient algebras and full centers in
the table below. The case $K_0\simeq W_0$ is included in the first row.
For the contact rows, let $g$ be the orthogonal representative of the
return automorphism, let $\rho$ be a reflection, and set
\[
 \rho_0=\id,\quad \rho_1=\rho,\qquad
 (-1)^\epsilon=\det(g)(-1)^{N(r-1)},\quad\epsilon\in\{0,1\}.
\]
In rank four, write $\nu=\xi_1\xi_2\xi_3\xi_4$.

\begin{center}
\renewcommand{\arraystretch}{1.15}
\begin{tabularx}{\textwidth}{@{}p{0.16\textwidth}Xp{0.26\textwidth}@{}}
\toprule
Constituent & Coefficient algebra and full center & Justification\\
\midrule
$W_N$\newline $(N\ne1)$ &
$W(1,N)$; center $0$ &
Theorem~\ref{witt:thm:semisimple}\\[4pt]
$W_1$ &
$\Alg(K_2,\id)$ or $\Alg(K_2,\rho)$; two nonisomorphic forms, each
with center $0$ &
Proposition~\ref{witt:prop:rankone-aut}\newline
Theorem~\ref{witt:thm:rankone-forms}\\[4pt]
$S_{N,0}$ &
$S'_{B+r-1}(1,N)$, $B\in\C$; center $0$ &
Theorem~\ref{sp:thm:coefficients}\newline
Corollary~\ref{sp:cor:cyclic}\\[4pt]
$K_N$\newline $N\ge1$, $N\ne4$ &
$\Alg(K_N,\rho_\epsilon)$; center $0$ &
Theorem~\ref{ct:thm:twists}\newline
Proposition~\ref{ct:prop:center}\\[4pt]
$K'_4$ &
$\Alg(K'_4,\rho_\epsilon)$; center $\C(\partial\nu)[0]$ if
$\epsilon=0$, and $0$ if $\epsilon=1$ &
Theorem~\ref{ct:thm:twists}\newline
Proposition~\ref{ct:prop:center}\\[4pt]
$CK_6$ &
$CK(6)$; center $0$ &
Theorem~\ref{ck:thm:cktwist}\\
\bottomrule
\end{tabularx}
\end{center}

Since $B+r-1\in\C$, the special
row lies in the stated family. Quotienting the contact rows by their
full centers gives $\mathfrak K_N^\epsilon$. Theorem~\ref{thm:ordinary-model}
allows no other constituents, so the list is exhaustive.
\end{proof}

\subsection{Polynomial growth and bounded weight spaces}\label{sec:growth}
For a finite subset $F$, let $V_j(F)$ denote the span of iterated
brackets in elements of $F$ involving at most $j$ factors. Highest weights below
are $L_0$-eigenvalues; thus a highest weight $h$ occurs in the physical
space $L[-h]$.

\begin{lemma}\label{gr:lem:highest}
Let $L$ contain the centerless Virasoro algebra. If the adjoint Virasoro
module of $L$ has a nontrivial simple highest- or lowest-weight
subquotient, then there is a finite subset $F\subset L$ such that
\begin{equation}\label{gr:eq:partition-lower}
 \dim V_{q+1}(F)\geq \sum_{n=0}^{q}p(n)\qquad(q\geq0),
\end{equation}
where $p(n)$ is the partition number. In particular, for this fixed
$F$, the function $\dim V_j(F)$ is not bounded by any polynomial in $j$.
\end{lemma}

\begin{proof}
Suppose that $N/P$ is a nontrivial simple highest-weight subquotient,
with highest vector $v$, so that $L_nv=0$ for $n>0$ and $L_0v=hv$.
Its central charge is zero. The unique simple highest-weight module at
central charge and highest weight zero is trivial. Thus $h\ne0$.

Choose a positive integer $m>24|h|$ and define
\[
 E_n=m^{-1}L_{mn}.
\]
Then $[E_r,E_s]=(r-s)E_{r+s}$, and $v$ is a highest vector of weight
$h/m$ for this copy of the centerless Virasoro algebra. At central
charge zero, the reducibility weights are
\[
 h_{r,s}=\frac{(3r-2s)^2-1}{24},\qquad r,s\geq1;
\]
see \cite[Theorem~2.8 and (2.9)--(2.10)]{BFT}. Every nonzero value has
absolute value at least $1/24$. Hence the Verma module of weight $h/m$
is simple, and its nonzero map to the module generated by $v$ is
injective.

For a partition $\lambda=(\lambda_1,\ldots,\lambda_s)$ of $n$, the PBW
vectors
\[
 E_{-\lambda_1}\cdots E_{-\lambda_s}v
\]
are therefore linearly independent. Choose a lift $\widetilde v\in N$
and let
\[
 F=\{E_{-1},E_{-2},\widetilde v\}.
\]
For $r\geq2$,
\[
 E_{-r}=\frac{1}{(r-2)!}(\ad E_{-1})^{r-2}E_{-2}.
\]
Thus a level-$n$ PBW vector has a lift obtained from at most $n+1$
factors from $F$. Taking all levels at most $q$ proves
\eqref{gr:eq:partition-lower}.

Suppose, for contradiction, that $\dim V_j(F)\leq Cj^d$. Choose an integer $r>d$. For $0\leq a_i\leq j$,
the vectors
\[
 E_{-1}^{a_1}E_{-2}^{a_2}\cdots E_{-r}^{a_r}v
\]
are linearly independent. Their lifts lie in
\[
 V_{1+jr(r+1)/2}(F),
\]
while their number is $(j+1)^r$. This contradicts the assumed bound
for large $j$. The lowest-weight case follows from the automorphism
$L_n\mapsto-L_{-n}$.
\end{proof}

\begin{proof}[Proof of Theorem~\ref{thm:growth-rigidity}]
Uniform boundedness implies linear commutator growth.
If $\dim L[k]\leq d$ for all $k$ and $F$ is finite, choose $a\geq0$
so that every weight component of every element of $F$ has index in
$[-a,a]$. Then
\[
 V_j(F)\subseteq\bigoplus_{|k|\leq aj}L[k],\qquad
 \dim V_j(F)\leq d(2aj+1)\leq d(2a+1)j.
\]
Linear growth is a special case of polynomial growth.

Conversely, assume polynomial commutator growth.
Lemma~\ref{gr:lem:highest} excludes every nontrivial simple highest- or
lowest-weight subquotient of the adjoint Virasoro module.
A weight $\Vir$-module with finite-dimensional weight spaces is
Harish--Chandra. Submodules and quotients remain weight modules because
polynomials in $L_0$ separate the finitely many weight components of
each vector. By \cite[Theorem~1, p.~225]{Mathieu}, a simple
Harish--Chandra module is highest weight, lowest weight, or of the
intermediate series. A nontrivial simple intermediate-series module
with integral support has dimension one at every nonzero integral
weight.

Fix nonzero integers $k,l$ and a parity $\bar p$. We prove by induction
on $\dim W$ that
\begin{equation}\label{gr:eq:induction}
 \dim W\leq \dim (U(\Vir)W)[l]
 \qquad\text{for every }W\subseteq L[k]_{\bar p}.
\end{equation}
For $W\ne0$, the finitely generated module $M=U(\Vir)W$ has a maximal
proper submodule $N$: the union of a chain of proper submodules cannot
contain all of a finite generating set. The simple quotient $Q=M/N$
is generated by the nonzero image of $W$ in $Q[k]$. Since $k\ne0$,
the quotient is nontrivial. Lemma~\ref{gr:lem:highest} excludes the
highest- and lowest-weight cases, so $Q$ is of intermediate-series type.
Hence
\[
 \dim Q[k]=\dim Q[l]=1.
\]
Let $W'=W\cap N$ and $M'=U(\Vir)W'$. Then
$\dim W'=\dim W-1$ and $M'\subseteq N$. By induction and exactness on
weight spaces,
\[
 \dim W=1+\dim W'
 \leq1+\dim M'[l]
 \leq1+\dim N[l]=\dim M[l].
\]
Taking $W=L[k]_{\bar p}$ and interchanging $k,l$ proves equality of all
nonzero weight dimensions in each parity. Hence
\[
 \dim L[k]=\dim L[1]\qquad(k\ne0).
\]
The remaining space $L[0]$ is finite-dimensional by hypothesis.
\end{proof}

\begin{corollary}\label{gr:cor:original}
Let $L$ be a simple Lie superalgebra containing the centerless Virasoro
algebra, with diagonal integral grading and finite-dimensional weight
spaces. If $L$ has polynomial commutator growth, then it is isomorphic,
as a parity-graded Lie superalgebra, to a member of the list in
Theorem~\ref{thm:synthesis}. In particular, the conclusion holds for
algebras satisfying the original conditions SC1--SC3 and the stated
grading and finiteness hypotheses.
\end{corollary}
\begin{proof}
Theorem~\ref{thm:growth-rigidity} gives the boundedness required in
Definition~\ref{def:superconformal}, so Theorem~\ref{thm:synthesis}
applies. The original condition SC3 is linear commutator growth.
\end{proof}

\section{Beyond polynomial growth}\label{sec:unbounded}
Under the hypotheses of Theorem~\ref{thm:growth-rigidity}, polynomial
growth forces uniformly bounded weight spaces. Finite-dimensional
weight spaces alone do not imply this bound, as the following
construction shows. Since an
ordinary Lie algebra is a Lie superalgebra with zero odd part, the
example also lies in the ambient super setting.

\subsection{A local Lie algebra}\label{sec:unbounded-local}
Let
\[
 \Vir_+=\bigoplus_{n>0}\C L_n,
 \qquad
 \Vir_-=\bigoplus_{n<0}\C L_n.
\]
Let $M=M^+(3)$ be the highest-weight Verma module generated by $v$ with
\[
 L_nv=0\quad(n>0),\qquad L_0v=3v,
\]
and let $N=M^-(-3)$ be the lowest-weight Verma module generated by $w$
with
\[
 L_nw=0\quad(n<0),\qquad L_0w=-3w.
\]
Both modules have central charge zero.

\begin{lemma}\label{unb:lem:irreducible}
The modules $M$ and $N$ are irreducible. Their level-$n$ weight spaces
have dimension $p(n)$.
\end{lemma}
\begin{proof}
At central charge zero, the reducibility weights are
\[
 \frac{(3r-2s)^2-1}{24}\qquad(r,s\geq1).
\]
The value $3$ would require $(3r-2s)^2=73$, which is impossible. Thus
$M$ is irreducible by the same criterion used in
Lemma~\ref{gr:lem:highest}. The automorphism $L_n\mapsto-L_{-n}$ gives
the assertion for $N$. The PBW theorem gives $p(n)$ basis vectors at
level $n$.
\end{proof}

The construction uses the following induced-module identity. For highest and
lowest Verma modules of weights $h$ and $h'$,
\begin{equation}\label{unb:eq:tensor-induced}
 M^+(h)\otimes M^-(h')
 \simeq U(\Vir)\otimes_{U(\C L_0)}\C_{h+h'},
\end{equation}
where $L_0$ acts on $\C_{h+h'}$ by $h+h'$.
Let $\Theta$ be the natural map from the induced module to the tensor
product, sending the inducing vector to the tensor of the two extremal
vectors. Give the source its quotient PBW filtration and the
target the total PBW filtration from $U(\Vir_-)$ and $U(\Vir_+)$.
A negative mode raises PBW degree by one on the first factor and
preserves the filtration on the second; a positive mode behaves
symmetrically. Thus $\operatorname{gr}\Theta$ is the canonical
isomorphism
\[
 S(\Vir_-\oplus\Vir_+)\simeq S(\Vir_-)\otimes S(\Vir_+).
\]
Both filtrations are exhaustive and bounded below, so $\Theta$ is an
isomorphism.

Let
\[
 \mathfrak a_0=\Vir\oplus\C d,
\]
where $d$ centralizes $\Vir$. Extend the $\Vir$-actions by letting $d$
act as $1$ on $M$ and as $-1$ on $N$. Since the two extremal weights
sum to zero, \eqref{unb:eq:tensor-induced} and Frobenius reciprocity
give a unique $\Vir$-homomorphism
\[
 B:M\otimes N\longrightarrow\mathfrak a_0,
 \qquad B(v\otimes w)=L_0+d.
\]
Here $L_0+d$ is annihilated by $\ad L_0$. Since $d$ acts trivially on
both $M\otimes N$ and $\mathfrak a_0$, the map $B$ is also
$\mathfrak a_0$-equivariant. Moreover, it is onto: for $n\ne0$,
\[
 L_n\cdot(L_0+d)=nL_n,
\]
so its image contains every $L_n$, then $L_0$, and finally $d$.

This defines the local Lie algebra
\[
 \mathfrak a_{-1}=N,
 \qquad \mathfrak a_0=\Vir\oplus\C d,
 \qquad \mathfrak a_1=M,
\]
with mixed bracket $[u,f]=B(u\otimes f)$ for $u\in M$ and $f\in N$.
The local bracket is prescribed only when the sum of the two degrees
belongs to $\{-1,0,1\}$.
The equivariance of $B$ is precisely the local Jacobi compatibility.
By the maximal construction for local Lie algebras
\cite[Proposition~2.5]{JLW}, which applies to arbitrary local parts,
there is a maximal graded Lie algebra
\[
 \mathfrak g=\bigoplus_{r\in\Z}\mathfrak g_r
\]
in which the local part embeds, with
\[
 \mathfrak g_0=\mathfrak a_0,\qquad
 \mathfrak g_{>0}=\mathcal L(M),
 \qquad
 \mathfrak g_{<0}=\mathcal L(N),
\]
where $\mathcal L$ denotes the free Lie algebra. The element $d$ gives
the auxiliary grading,
\begin{equation}\label{unb:eq:dgrading}
 [d,x]=rx\qquad(x\in\mathfrak g_r).
\end{equation}

\subsection{The simple quotient}\label{sec:unbounded-simple}
\begin{lemma}\label{unb:lem:proper-ideal}
Every proper ideal of $\mathfrak g$ has zero intersection with
$\mathfrak g_{-1}\oplus\mathfrak g_0\oplus\mathfrak g_1$.
\end{lemma}
\begin{proof}
Equation~\eqref{unb:eq:dgrading} shows that every ideal is graded by the
auxiliary degree: polynomials in $\ad d$ separate the finitely many
homogeneous components of each element.

Let $I$ be an ideal. If $I\cap M\ne0$, irreducibility gives $M\subseteq I$.
Since
\[
 [M,N]=\mathfrak a_0,
\]
we obtain $\mathfrak a_0\subseteq I$, and $[d,N]=-N$ then gives
$N\subseteq I$. The local part generates $\mathfrak g$, so
$I=\mathfrak g$. The same argument applies if $I\cap N\ne0$.

Suppose now that $I\cap\mathfrak a_0\ne0$. The action of
$\mathfrak a_0$ on $M$ is faithful. Indeed, if $a+cd$ acts trivially,
with $a\in\Vir$, then $a$ acts on $M$ as a scalar. Hence $[a,\Vir]$ acts
trivially. The nontrivial $\Vir$-module $M$ is faithful because $\Vir$
is simple. Thus $a$ centralizes $\Vir$, so $a=0$, and then $c=0$.
Consequently a nonzero element of $I\cap\mathfrak a_0$ has a nonzero
bracket with some element of $M$, reducing to the preceding case.
\end{proof}

Let $R$ be the sum of all proper ideals of $\mathfrak g$. By
Lemma~\ref{unb:lem:proper-ideal}, the three local components of $R$ are
zero, so $R$ is proper. Since it contains every proper ideal, the quotient
\[
 L=\mathfrak g/R
\]
is simple. The local part embeds in $L$.

\subsection{Finite but unbounded weight spaces}\label{sec:unbounded-weights}
Consider the physical grading defined by $L_0$, rather than the auxiliary
grading defined by $d$. Every $L_0$-weight in $M$ is of the form $3+n$
with $n\geq0$. Hence every weight in the auxiliary component
$\mathcal L_r(M)$ is at least $3r$. Moreover,
\[
 \mathcal L_r(M)\subseteq M^{\otimes r}.
\]
For a fixed $L_0$-eigenvalue $k>0$, the relevant tensor weight space has
dimension at most
\[
 \sum_{\substack{n_1+\cdots+n_r=k-3r\\ n_i\geq0}}
 p(n_1)\cdots p(n_r),
\]
which is finite, and only $1\leq r\leq k/3$ can occur. The negative free
part has the symmetric property. The degree-zero algebra
$\Vir\oplus\C d$ has one-dimensional nonzero $L_0$-weight spaces and
two-dimensional zero-weight space. Thus every $L_0$-weight space of
$\mathfrak g$ is finite-dimensional. The action of $\ad L_0$ on
$\mathfrak g$ is diagonal with integral spectrum. Since $R$ is stable
under this action, polynomials in $\ad L_0$ show that $R$ is a direct
sum of its intersections with the eigenspaces. The quotient $L$ therefore
also has diagonal $L_0$-action and finite-dimensional integral weight
spaces.

The positive and negative free parts have no zero eigenvalue, and the
local part survives in $L$. At eigenvalue $3+n$, the local summands
$\Vir$ and $M$ contribute $1+p(n)$ independent vectors; the summands
$\Vir$ and $N$ do the same at eigenvalue $-3-n$. Thus, in the physical
notation \eqref{eq:weight},
\[
 \dim L[0]=2,
 \qquad
 \dim L[\pm(3+n)]\geq1+p(n)\qquad(n\geq0).
\]
The lower bounds prove the unboundedness asserted in
Theorem~\ref{thm:unbounded-example}.

The algebra is also finitely generated. The four elements
$L_{-2},L_{-1},L_1,L_2$ generate $\Vir$; the vectors $v$ and $w$ generate
$M$ and $N$ as $\Vir$-modules; and
\[
 d=[v,w]-L_0.
\]
Hence the images of these six elements generate $L$.

The algebra $L$ is not abstractly isomorphic to $\Vir$. The
centralizer of $L_0$ in $L$ contains the independent elements $L_0$ and
$d$. On the other hand, the centralizer of every nonzero element of the
Witt algebra is one-dimensional: if
$f(t)\partial_t$ and $g(t)\partial_t$ commute, then
$(g/f)'=0$, so $g$ is a scalar multiple of $f$.

\subsection{Exponential growth}\label{sec:unbounded-growth}

\begin{proposition}\label{unb:prop:exponential}
The elements $L_{-1}$ and $v$ generate a free Lie algebra of rank two
in $L$. In particular, $L$ has exponential commutator growth.
\end{proposition}
\begin{proof}
Use the finite generating set
$S=\{L_{-2},L_{-1},L_1,L_2,v,w\}$.
Define
\[
\begin{array}{lll}
 e_1=L_{-1},& f_1=L_1,& h_1=-2L_0,\\[1mm]
 e_2=v,& f_2=-6w,& h_2=-6(L_0+d).
\end{array}
\]
Set
\[
 A=(a_{ij})_{1\leq i,j\leq2}
   =\begin{pmatrix}-2&-6\\-6&-24\end{pmatrix}.
\]
The brackets in the embedded local part give
\begin{gather*}
 [e_i,f_j]=\delta_{ij}h_i,\qquad [h_i,h_j]=0,\\
 [h_i,e_j]=a_{ij}e_j,\qquad [h_i,f_j]=-a_{ij}f_j.
\end{gather*}
In particular, the off-diagonal mixed brackets vanish because
$L_1v=0$ and $L_{-1}w=0$.

Let $\mathfrak g(A)$ be the complex generalized Kac--Moody algebra
with this symmetric matrix presentation; see
\cite[\S1]{Borcherds}. Both diagonal entries are negative and the
off-diagonal entries are nonzero, so the presentation imposes no
additional Serre or commuting relations among the positive generators.
The displayed identities therefore define a homomorphism
\[
 \Phi:\mathfrak g(A)\longrightarrow L.
\]
Let $\alpha_i$ be the root degree of $e_i$ in $\mathfrak g(A)$.
On the component of degree $a\alpha_1+b\alpha_2$, the two Cartan
eigenvalues are $(-2a-6b,-6a-24b)$. Since $\det A=12$, distinct root
degrees have distinct joint weights. Hence $I=\ker\Phi$ is root-graded,
and its degree-zero component is zero because $L_0$ and $d$ remain
independent in $L$.

Let $\overline I$ denote the conjugate of $I$ with respect to the
real form of $\mathfrak g(A)$ defined by its generators.
The ideal $I+\overline I$ is root-graded, stable under conjugation,
and has zero degree-zero component. Its real fixed part is thus a
graded ideal disjoint from the real Cartan subalgebra.
By \cite[Corollary~2.6]{Borcherds}, this real ideal is zero, so
$I+\overline I=0$ and $\Phi$ is injective.
For this matrix, the positive part is the free Lie algebra on $e_1,e_2$;
see \cite[Proposition~2.1 and \S5.2]{Jurisich}.

The degree-$j$ part of the free Lie algebra on two generators has
dimension asymptotic to $2^j/j$ by the Witt formula
\cite[Proposition~5.2, equation~(13)]{Jurisich}.
By the Jacobi identity, left-normed words span all Lie words in $S$,
so the words of length at most $j$ span a space of dimension at most
$\sum_{r=1}^j6^r$. Thus, for all sufficiently large $j$,
\[
 \frac{2^{j-1}}j\leq\dim V_j(S)\leq\frac{6^{j+1}-6}{5},
\]
which proves the exponential growth assertion.
\end{proof}

This completes the proof of Theorem~\ref{thm:unbounded-example}.
Finite-dimensional weight spaces alone therefore cannot replace the
growth hypothesis in Theorem~\ref{thm:growth-rigidity}.

\subsection{Further questions}\label{sec:further-questions}

\begin{questionnote}{question:growth-threshold}{Growth beyond polynomial}
Under the integral grading and finite-dimensionality assumptions of
Theorem~\ref{thm:growth-rigidity}, bounded weight multiplicities also
follow from the condition
\[
 \log\bigl(1+\dim V_j(F)\bigr)=o(\sqrt j)
 \qquad(j\to\infty)
\]
for every finite $F\subset L$.
This condition is weaker than the requirement of polynomial growth.
The partition bound in Lemma~\ref{gr:lem:highest} excludes nontrivial
simple highest- and lowest-weight subquotients, and the rest of the
proof of Theorem~\ref{thm:growth-rigidity} applies unchanged.
What is the weakest growth condition with this consequence, especially
for simple Lie algebras?
\end{questionnote}

\begin{questionnote}{question:endpoint}{Isomorphisms between representatives}
Which representatives are isomorphic? For the special family,
Proposition~\ref{sp:prop:parameters} gives the sufficient condition
$B'\in\pm B+\Z$. This paper does not settle necessity or isomorphisms
between different families beyond the low-rank cases above.
\end{questionnote}

\begin{questionnote}{question:marking}{Prescribed Virasoro data}
Which isomorphisms preserve a prescribed Virasoro subalgebra, its
generators, or the physical grading? The constructions above may change
the translation, grading, and Virasoro embedding; they leave this
question unresolved.
\end{questionnote}

\section*{Acknowledgments}
Yong Yang acknowledges support from the Simons Foundation, Award No.~918096.

OpenAI's ChatGPT was used for mathematical exploration and editorial assistance. The authors checked the final arguments and are responsible for the contents of the paper. H. Chen thanks Professor Yucai Su for helpful discussions.

\end{document}